\documentclass[12pt,leqno]{amsart}
\usepackage{amssymb}
\usepackage{amsmath,amssymb}
\renewcommand{\baselinestretch}{1.1}

\newtheorem{theorem}{Theorem}[section]

\newtheorem{lemma}{Lemma}[section]
\newtheorem{proposition}{Proposition}[section]
\newtheorem{condition}{Condition}[section]

\begin{document}
\theoremstyle{plain}
\newtheorem{MainThm}{Theorem}
\newtheorem{thm}{Theorem}[section]
\newtheorem{clry}[thm]{Corollary}
\newtheorem{prop}[thm]{Proposition}
\newtheorem{lem}[thm]{Lemma}
\newtheorem{deft}[thm]{Definition}
\newtheorem{hyp}{Assumption}
\newtheorem*{ThmLeU}{Theorem (J.~Lee, G.~Uhlmann)}

\theoremstyle{definition}
\newtheorem{rem}[thm]{Remark}
\newtheorem*{acknow}{Acknowledgments}
\numberwithin{equation}{section}
\newcommand{\eps}{{\varphi}repsilon}
\renewcommand{\d}{\partial}
\newcommand{\re}{\mathop{\rm Re} }
\newcommand{\im}{\mathop{\rm Im}}
\newcommand{\R}{\mathbf{R}}
\newcommand{\C}{\mathbf{C}}
\newcommand{\N}{\mathbf{N}}
\newcommand{\D}{C^{\infty}_0}
\renewcommand{\O}{\mathcal{O}}
\newcommand{\dbar}{\overline{\d}}
\newcommand{\supp}{\mathop{\rm supp}}
\newcommand{\abs}[1]{\lvert #1 \rvert}
\newcommand{\csubset}{\Subset}
\newcommand{\detg}{\lvert g \rvert}
\newcommand{\dd}{\mbox{div}\thinspace}
\newcommand{\www}{\widetilde}
\newcommand{\ggggg}{\mbox{\bf g}}
\newcommand{\ep}{\varepsilon}
\newcommand{\la}{\lambda}
\newcommand{\va}{\varphi}
\newcommand{\ppp}{\partial}
\newcommand{\ooo}{\overline}
\newcommand{\wwwx}{\widetilde{x}_0}
\newcommand{\sumkj}{\sum_{k,j=1}^n}
\newcommand{\walpha}{\widetilde{\alpha}}
\newcommand{\wbeta}{\widetilde{\beta}}
\newcommand{\weight}{e^{2s\va}}
\newcommand{\fdif}{\partial_t^{\alpha}}
\newcommand{\OOO}{\Omega}
\newcommand{\LLL}{L_{\lambda,\mu}}
\newcommand{\LLLL}{L_{\widetilde{\la},\widetilde{\mu}}}
\newcommand{\ppdif}[2]{\frac{\partial^2 #1}{{\partial #2}^2}}
\newcommand{\uu}{\mbox{\bf u}}
\newcommand{\vv}{\mbox{\bf v}}
\newcommand{\ww}{\mbox{\bf w}}
\renewcommand{\v}{\mathbf{v}}
\newcommand{\y}{\mathbf{y}}
\newcommand{\ddd}{\mbox{div}\thinspace}
\newcommand{\rrr}{\mbox{rot}\thinspace}
\newcommand{\Y}{\mathbf{Y}}
\newcommand{\w}{\mathbf{w}}
\newcommand{\z}{\mathbf{z}}
\newcommand{\G}{\mathbf{G}}
\newcommand{\f}{\mathbf{f}}
\newcommand{\F}{\mathbf{F}}
\newcommand{\dddx}{\frac{d}{dx_0}}
\newcommand{\CC}{_{0}C^{\infty}(0,T)}
\newcommand{\HH}{_{0}H^{\alpha}(0,T)}
\newcommand{\llll}{L^{\infty}(\Omega\times (0,t_1))}
\renewcommand{\baselinestretch}{1.5}
\renewcommand{\div}{\mathrm{div}\,}  %div
\newcommand{\grad}{\mathrm{grad}\,}  %grad
\newcommand{\rot}{\mathrm{rot}\,}  %rot

\title
[]
{Carleman estimate and an inverse source problem for the Kelvin-Voigt model for
viscoelasticity}

%\pagestyle{myheadings}
%\markboth{O.~Yu.~Imanuvilov and \, M.~Yamamoto}
%{Calder\'on problem for Maxwell's equations}

\author{
O.~Yu.~Imanuvilov and \, M.~Yamamoto }

\thanks{
Department of Mathematics, Colorado State
University, 101 Weber Building, Fort Collins, CO 80523-1874, U.S.A.
e-mail: {\tt oleg@math.colostate.edu}
Partially supported by NSF grant DMS 1312900}\,
\thanks{ Department of Mathematical Sciences, The University
of Tokyo, Komaba, Meguro, Tokyo 153, Japan,\\
PeoplesfFriendship University of Russia (RUDN University)
6 Miklukho-Maklaya St, Moscow, 117198, Russian Federation
e-mail: {\tt myama@ms.u-tokyo.ac.jp}
Supported by Grant-in-Aid for Scientific Research (S) 15H05740 of
Japan Society for the Promotion of Science and
prepared with the support of the "RUDN University Program 5-100".
}

%\author{
%O.~Yu.~Imanuvilov, Department of Mathematics, Colorado
%State University, 101 Weber Building, Fort Collins, CO 80523-1874,
%U.S.A. E-mail: {\tt oleg@math.colostate.edu}},\,
%M.~Yamamoto  Department of Mathematical Sciences, The University
%of Tokyo, Komaba, Meguro, Tokyo 153, Japan e-mail:
%{myama@ms.u-tokyo.ac.jp}\,}

\date{}

\maketitle

\begin{abstract}
We consider the Kelvin-Voigt model for the
viscoelasticity, and prove a Carleman estimate for functions without
compact supports.
Then we apply the Carleman estimate to prove the Lipschitz stability
in determining a spatial varying function in an external source term
of Kelvin-Voigt model by a single measurement. Finally as a related
system, we consider an isothermal compressible fluid system and
apply the Carleman estimate to establish the Lipschitz stability for
an inverse source problem for the compressible fluid system.
\end{abstract}

\section{Introduction and main results}
Let $T$ be a positive constant, $x' = (x_1,..., x_n) \in \Bbb R^n$,
and $\Omega$ be a bounded domain in $\Bbb R^n$ with $\partial\Omega\in
C^\infty$, let $\vec\nu =\vec\nu(x')$ be the unit outward normal vector at
$x'$ to $\ppp\Omega$.
Here we understand $x' \in \Bbb R^n$ as the spatial variable and
$x_0$ as the time variable.  We set
$x = (x_0, x') = (x_0, x_1, ..., x_n)$, and
$$
Q:= (-T,T) \times \OOO, \quad \Sigma := (-T,T)\times\partial\Omega.
$$

Here and henceforth $i = \sqrt{-1}$ and
$\cdot^T$ denotes the transposes of vectors and matrices
under consideration, and $D= (D_0, D')$, $D_0=\frac{1}{i}
\partial_{x_0}$, $D'=(\frac{1}{i}\partial_{x_1},\dots,
\frac{1}{{i}}\partial_{x_n})$, $\nabla' = (\ppp_{x_1}, \dots,
\ppp_{x_n})$, $\nabla = (\ppp_{x_0}, \nabla')$.

In the cylinder domain $Q$, we consider the Kelvin-Voigt model:
\begin{equation}\label{1.1}
\rho\partial_{x_0}^2\mbox{\bf w}
= L_{\lambda,\mu}(x,D')\ppp_{x_0}\mbox{\bf w}
+ L_{\www{\lambda},\www{\mu}}(x,D') \ww + \F \quad \mbox{in $Q$},
\end{equation}
\begin{equation}\label{1.2}
\ww\vert_{\Sigma}=0,
\end{equation}
\begin{equation}\label{1.3}
\ww(0,x') = 0, \quad x'\in \Omega.
\end{equation}
Here $\mbox{\bf w}(x) =(w_1(x),\dots , w_n(x))^T$ is the displacement and
$\mbox{\bf F}(x)=(F_1(x),\dots, F_n(x))^T$ is an external force.
For Lam\'e coefficients $\lambda(x)$ and $\mu(x)$,
the partial differential operator $L_{\lambda,\mu}(x,D')$ is defined by
\begin{eqnarray}\nonumber
L_{\lambda,\mu}(x,D')\mbox{\bf w} = \mu(x)\Delta\mbox{\bf w}
+ (\mu(x)+ \lambda(x))\nabla'\mbox{div}\,\mbox{\bf w}\nonumber\\
+ (\mbox{div}\,\mbox{\bf w})\nabla'\lambda + (\nabla'\mbox{\bf w}
+ (\nabla'\mbox{\bf w})^T) \nabla'\mu, \qquad x \in Q
\end{eqnarray}
and $L_{\www{\lambda},\www{\mu}}(x,D')$ is defined similarly.

The equation (\ref{1.1}) is called the Kelvin-Voigt model and is one
model equation for the viscoelasticity.
The viscoelasticity indicates a mixed physical property
of the viscosity and the elasticity, and is frequently observed in human
tissues.  Thus the viscoelasticity is important for example for the medical
diagnosis.  In the diagnosis, the main task is to detect some anomaly of
spatially varying coefficients in (\ref{1.1}) which may indicate
some organizational abnormalities such as tumors.  Hence inverse problems of
determining coefficients $\la(x')$, $\mu(x')$, $\www{\la}(x')$,
$\www{\mu}(x')$ by available boundary data, are demanded from medical
points of view.  Such inverse problems of determining coefficients
can be solved by the corresponding inverse source problems which can be the
linearized problems of inverse coefficient problems,
but we do not here discuss details for the inverse coefficient problems.

As for medical applications related to the Kelvin-Voigt model, see
Catheline, Gennisson, Delon, Fink, Sinkus, Abouelkaram and Culiolic \cite{CGD},
and Royston, Mansy and Sandler \cite{Roy}.
As for applications of other viscoelasticity models to the diagnosis,
we refer to de Buhan \cite{dB},
Sinkus, Tanter, Xydeas, Catheline, Bercoff and Fink \cite{S}
and the references therein.
As monographs on the viscoelasticity, the readers can consult
Lakes \cite{L}, and Renardi, Hrusa and Nohel \cite{RHN}.

For the inverse problems related to the medical diagnosis, one
of theoretical fundamentals is the Carleman estimate.
The purpose of this paper is to establish a Carleman estimate for
(\ref{1.1}) and apply it to an inverse source problem of determining
a spatially varying factor of the external source $\F(x)$.
Also as is seen later,
the principal part of (\ref{1.1}) is a strongly coupled
parabolic Lam\'e system $\rho\ppp_{x_0} - L_{\la,\mu}(x,D')$, which
causes a serious difficulty for proving the Carleman estimate for (\ref{1.1}).
In addition to the Kelvin-Voigt model, we consider an isothermal
compressible fluid dynamics system which mathematically,
modulo nonlinear terms, is a first-order partial differential equation coupled
with the parabolic Lam\'e system
and apply the Carleman estimate to establish the Lipschitz stability for
some inverse source problem.
\\
\vspace{0.2cm}

Now we reduce (\ref{1.1}) to an integro-parabolic system.  Setting
\begin{equation}\label{1.5}
\uu = \ppp_{x_0}\ww,
\end{equation}
from (\ref{1.1}) - (\ref{1.3}), we readily obtain
\begin{equation}\label{gora1A}
P(x,D)\mbox{\bf u} \equiv \rho\partial_{x_0}\mbox{\bf u}
- L_{\lambda,\mu}(x,D')\mbox{\bf u}
- \int_0^{x_0}L_{\widetilde \lambda, \widetilde\mu}
(x, \widetilde x_0, D')\mbox{\bf u}(\widetilde x_0, x')d\widetilde x_0
= \mbox{\bf F}\quad\mbox{in}\,\,Q,
\end{equation}
\begin{equation}\label{gora2A}
\mbox{\bf u}\vert_{\Sigma}=0.
\end{equation}
Here and henceforth, for the sake of more generality, assuming
in (\ref{gora1A})
that $\www{\lambda}, \www{\mu}$ depend on $x=(x_0,x') \in Q$ and
$\www{x}_0 \in (-T,T)$, we set
\begin{eqnarray}
L_{\www\lambda,\www\mu}(x, \wwwx,D')\mbox{\bf u}
= \www\mu(x,\wwwx)\Delta\mbox{\bf u}
+ (\www\mu(x,\wwwx)+ \www\lambda(x,\wwwx))\nabla'\mbox{div}\,\mbox{\bf u}
                                    \nonumber\\
+ (\mbox{div}\,\mbox{\bf u})\nabla'\widetilde\lambda + (\nabla'\mbox{\bf u}
+ (\nabla'\mbox{\bf u})^T) \nabla'\widetilde\mu, \qquad (x,\widetilde x_0) \in
Q\times [-T,T],  \nonumber
\end{eqnarray}

We mainly consider system (\ref{gora1A}) with boundary (\ref{gora2A}).
Henceforth the coefficients  $\rho,$ $\lambda$, $\mu$, $\www\la$, $\www\mu$
are assumed to satisfy
\begin{equation}\label{pinok0}
\rho, \lambda, \mu \in C^2(\overline{Q}), \quad
 \rho(x)>0,\, \mu(x) > 0, \thinspace \lambda(x) + \mu(x) > 0 \quad \mbox{for}\,
x \in \overline{Q}
\end{equation}
and
\begin{equation}\label{pinok1}
\www\lambda, \www\mu\in C^2([-T,T] \times \ooo{\OOO}\times [-T,T]).
\end{equation}

More precisely, the first main purpose of this paper is to establish
\begin{enumerate}
\item
a Carleman estimate for functions without compact supports;
\item
the Lipschitz stability in an inverse source problem of determining
spatially varying factor of the external force ${\bf F}$
\end{enumerate}
for (\ref{gora1A}).

A Carleman estimate is an $L^2$-weighted estimate of solution $\uu$ to
(\ref{gora1A}) which holds uniformly in large parameter.
Carleman estimates have been well
studied for single equations (e.g., H\"ormander \cite{H}, Isakov \cite{Is}).
A Carleman estimate yields several important results such as
the unique continuation, the energy estimate called an observability
inequality and the stability in inverse problems.
However for systems whose principal part is coupled, for example even for
isotropic Lam\'e system (that is, (\ref{gora1A}) with $\www\lambda
= \www\mu \equiv 0$), the Carleman estimate is difficult to be proved for
functions $\mathbf{u}$ whose supports are not compact in $Q$.
In particular, for (\ref{gora1A}), no such Carleman estimates are not known.

In order to prove the stability global in $\Omega$ for the inverse problem
with lateral data limited on some subboundary, we need a Carleman estimate for
functions without compact supports.  Otherwise we have to observe
the extra data on the whole lateral subbboundary $(0,T)\times \ppp\Omega$.

In establishing a Carleman estimate for our system (\ref{gora1A}) for
non-compactly supported $\uu$, we should emphasize
the two main difficulties:
\begin{itemize}
\item
The principal part $\rho\ppp_{x_0} - L_{\lambda,\mu}(x,D')$ is
strongly coupled.
\item
The Lam\'e operator $L_{\www\la,\www\mu}(x,\wwwx,D')$ appears as an integral.
\end{itemize}

As for Carleman estimates for functions without compact supports and
applications to inverse problems for the Lam\'e system
with $L_{\la,\mu}(x,D') = 0$,
we refer to  Bellassoued, Imanuvilov and
Yamamoto \cite{BIY}, Bellassoued and Yamamoto \cite{BY},
Imanuvilov, Isakov and Yamamoto \cite{IIY},
Imanuvilov and Yamamoto \cite{IY4} - \cite{IY7}.
In this paper, we modify the arguments in those papers
and establish a Carleman estimate for system (\ref{gora1A}) where
$\mathbf{u}$ does not
have compact support.  Then we apply the Carleman estimate to an inverse
source problem by modifying the method in Imanuvilov and Yamamoto \cite{IY1} -
\cite{IY3} and Beilina, Cristofol, Li and M. Yamamoto \cite{BCLY}
discussing scalar hyperbolic equations.
As for the methodology for applying Carleman estimates to inverse
problems, we refer to a pioneering paper Bukhgeim and Klibanov
\cite{BK}, and also Beilina and Klibanov \cite{BeiKl},
Bellassoued and Yamamoto \cite{BY2017}, Klibanov \cite{Kl1, Kl2},
Klibanov and Timonov \cite{KlT}.
\vspace{0.2cm}

For the statement of the Carleman estimate for (\ref{gora1A}), we need to
introduce notations and a definition.

Let $\overline{\xi}$ denote the complex conjugate of $\xi \in \C$.
We set $<a,b> = \sum^n_{k=0} a_k\overline b_k$ for
$a = (a_0, \dots, a_n), b = (b_0, \dots , b_n) \in \Bbb C^n$,
$\xi=(\xi_0,\dots,\xi_n), \xi'=(\xi_1,\dots, \xi_n), \widetilde \xi=(\xi_0,
\dots,
\xi_{n-1}), \www\nabla=(\partial_{y_0},\dots \partial_{y_{n-1}})$,
$\zeta=(\xi_0,\dots, \xi_{n-1},\widetilde s),
\widetilde D=(D_0,\dots, D_{n-1})$.
\par
For $\beta \in C^2(\ooo Q)$, we introduce the symbol:
$$
p_{\rho,\beta}(x,\xi)=i\rho(x)\xi_0 + \beta(x)\vert \xi'\vert^2.
$$

Let $\Gamma_0$ be some relatively open subset on $\partial\Omega$.
We set $\www\Gamma=\partial\Omega\setminus\Gamma_0$ and
$\Sigma_0=(-T,T)\times \Gamma_0,$ $\www\Sigma=(-T,T)\times \www\Gamma.$

In order to prove the Carleman estimate for the viscoelastic Lam\'e  system,
we assume the existence of a real-valued function $\psi$
which is pseudoconvex
with respect to the symbols $p_{\rho,\mu}(x,\xi)$ and
$p_{\rho,\la+2\mu}(x,\xi)$.  More precisely, we can state as follows.

For functions $f(x,\xi)$ and $g(x,\xi)$, we introduce the Poisson bracket
$$
\{f,g\}=\sum_{j=0}^n\left(\partial_{\xi_j} f\partial_{x_j} g
-\partial_{\xi_j} g \partial_{x_j} f\right).
$$

Denote
\begin{equation}\label{pauk}
\widetilde \varphi (x_0)=\frac{1}{(x_0+T)^3(T-x_0)^3}.
\end{equation}
We introduce
\begin{condition}\label{A1}  {\it We say that a function}
$\varphi\in C^{0,2}(\overline Q)$ with $\partial_{x_0} \varphi\in L^\infty(Q)$,
{\it is pseudoconvex with respect to the symbol} $p_{\rho,\beta}(x,\xi)$
{\it if there exists a constant} $C_1>0$ {\it such that}
\begin{equation}\label{ma}
\frac{\mbox{Im}\{\overline
p_{\rho,\beta}(x,\xi_0,\overline{\widetilde\zeta}),{p_{\rho,\beta}(x,\xi_0,\widetilde\zeta)}\}}{\vert s\vert}>  C_1\widetilde\varphi(x_0)M(\xi, \widetilde\varphi(x_0) s)^2
\quad \forall (x,\xi,s)\in \mathcal S,
\end{equation}
{\it where}
$$
\mathcal S=\{(x,\xi,s);\, x\in\overline{Q},\,\, M(\xi,\widetilde \varphi(x_0) s)=1,
p_{\rho,\beta}(x,\xi_0,\widetilde\zeta)=0\},
$$
$$
\widetilde\zeta=(\xi_1+i\vert s\vert\partial_{x_1}
\varphi,\dots,\xi_n+i\vert s\vert\partial_{x_n} \varphi)
$$
and
\begin{equation}\label{boykot}
M(\xi,s)=(\xi_0^2+\sum_{i=1}^n\xi_i^4+s^4)^\frac 14.
\end{equation}
\end{condition}
\noindent
We assume that there exists a positive constant $C_2$ such that
\begin{equation}\label{pinok220}
\partial_{\vec\nu}\varphi(x)\vert_{\overline{\Sigma_0}}<0
\quad \mbox{and}
\quad \vert\partial_{\vec\nu}\varphi(x)\vert \ge C_2 \widetilde \varphi(x_0)  \quad \forall x\in \overline{\Sigma_0}
\end{equation}
and
\begin{equation}\label{rep}
-\partial_{\vec \nu}\varphi(x)>\frac {1}{\root\of 2}\root\of{\frac{\mu(x)}{(\lambda+2\mu)(x)}} \vert \partial_{\vec\tau}\varphi(x)\vert \quad \forall x\in \Sigma_0, \quad \forall \vec\tau\in T(\Gamma_0), \,\, \vert \vec\tau\vert=1.
\end{equation}

If the pair $(\widetilde\mu,\widetilde \lambda)$ is not identically equal
to zero in $Q$, then we assume that
\begin{equation}\label{pinok22}
\partial_{x_0}\varphi(x)<0 \quad \mbox{on}\,\,(0,T]\times \overline \Omega
\quad \mbox{and}\quad \partial_{x_0}\varphi(x)>0 \quad \mbox{on}\,\,[-T,0)
\times \overline\Omega.
\end{equation}

Let us assume
\begin{equation}\label{giga1}
\varphi(x)<0\quad \mbox{on}\,\, \overline Q,\quad
\nabla' \varphi(x) \ne  0\quad \mbox{on}
\,\,\overline Q.
% \quad \partial^2_{x_0}\varphi(0,x')>0\quad\mbox {on}\,\, \overline\Omega.
\end{equation}
Let $a_\beta(x,\xi')=\beta(x)\vert \xi'\vert^2$.  Furthemore we introduce
\begin{condition}\label{A2}
{\it We say that a function} $\varphi(x)\in C^{0,2}(\overline Q)$ {\it is
pseudoconvex with respect to the symbol} $a_\beta(x,\xi')$ {\it if there exists
a constant} $C_3>0$ independents of $x_0, s, \xi'$ {\it such that}
$$
\frac{\mbox{Im}\{\overline
a_\beta(x,\xi'-i\vert s\vert\nabla'\varphi),{a_\beta(x,\xi'+i\vert s\vert\nabla'\varphi)}\}}{\vert  s\vert}>C_3\widetilde\varphi(x_0)\vert (\xi', \widetilde\varphi(x_0) s)\vert^2
\quad \forall (x,\xi',s)\in
\mathcal K,
$$
{\it where}
$$
\mathcal K=\{(x,\xi',s); \, x\in\overline{Q},\,\, \vert(\xi',s)\vert=1,
a_\beta(x,\xi'+i\vert s\vert\nabla'\varphi)=0\}.
$$
\end{condition}

Let $\Bbb N = \{ 1,2,3,... \}$,
$\alpha=(\alpha_0,\dots,\alpha_n) =: (\alpha_0, \alpha')$ with
$\alpha_0, ..., \alpha_n \in \Bbb N \cup \{ 0\}$ and
$\vert \alpha\vert=2\alpha_0+\sum_{j=1}^n\alpha_j$,
$\ppp_x^{\alpha} = \ppp_{x_0}^{\alpha_0}\ppp_{x_1}^{\alpha_1} \cdots
\ppp_{x_n}^{\alpha_n}$.

Finally we assume that
\begin{equation}\label{1}
\lim_{x_0\rightarrow -T+0} \varphi(x_0,x')=\lim_{x_0\rightarrow T-0} \varphi(x_0,x')=-\infty,
\end{equation}
and
\begin{equation}\label{lomka1}
\frac{C_4}{(x_0+T)^3(T-x_0)^3}\le\sum_{\vert (0,\alpha_1,..., \alpha_n)\vert
\le 2} \vert \partial_x^\alpha\varphi(x)\vert
\le \frac{C_5}{(x_0+T)^3(T-x_0)^3}
\end{equation}
for $x \in Q$,
and
\begin{equation}\label{lomka}\vert \partial_{x_0} \varphi(x)\vert\le \frac{C_6}{(x_0+T)^{4}(T-x_0)^{4}} \quad x\in ((0,T)\times \overline\Omega)\cup ((-T,0)\times
\overline\Omega).
\end{equation}

%Using the function $\psi$, we introduce the function $\varphi(x)$ by
%\begin{equation}\label{giga3}
%\varphi(x)=e^{\tau\psi(x)},\quad \tau >1,
%\end{equation}
%where the parameter $\tau>0$ will be fixed below.

For any function $\mbox{\bf f}=(f_1,\dots ,f_n)$, we introduce the differential
form $\omega_{\mbox{\bf f}}=\sum_{j=1}^n f_jdx_j.$ Then $d\omega_{\mbox{\bf f}}
= \sum_{k<j}^n(\partial_{x_j} f_k
- \partial_{x_k} f_j) dx_k\wedge dx_j.$

We identify the differential form $d\omega_{\mbox{\bf f}}$ with the
vector-function:
$$
d\omega_{\mbox{\bf f}}=\left (\partial_{x_2} f_1
-\partial_{x_1} f_2,\dots , \partial_{x_n} f_1
-\partial_{x_1}f_n,\partial_{x_3} f_2
-\partial_{x_2} f_3,\dots , \partial_{x_n} f_2
-\partial_{x_2} f_n, \dots,\partial_{x_n}f_{n-1}
-\partial_{x_{n-1}} f_n\right).
$$

Denote
\begin{eqnarray}\label{gopak}
\Vert \mbox{\bf u}\Vert^2_{\mathcal B(\varphi,s,Q)}=
\int_Q\Biggl(\sum^3_{\vert\alpha\vert=0, \alpha_0\le 1}(s\widetilde\varphi)^{4-2\vert\alpha\vert}
\vert
\partial^{\alpha}\text{\bf u}\vert^2
+ s\widetilde\varphi\vert\nabla' d\omega_{\mbox{\bf u}}\vert^2
\nonumber\\+(s\widetilde\varphi)^3\vert d\omega_{\mbox{\bf u}}\vert^2
+ s\widetilde\varphi\vert\nabla' \mbox{div}\,\text{\bf u}\vert^2
+ (s\widetilde\varphi)^3\vert\mbox{div}\,\text{\bf u}\vert^2\Biggr)e^{2s\varphi}
dx,
\end{eqnarray}
$$
\Vert \mbox{\bf u}\Vert^2_{\mathcal X(\varphi,s,\Omega)}
= \int_\Omega\Biggl(\sum^2_{\vert\alpha'\vert=0}(s\widetilde\varphi)
^{4-2\vert\alpha\vert}\vert
\partial^{\alpha'}\text{\bf u}\vert^2
+ s\widetilde\varphi\vert\nabla' d\omega_{\mbox{\bf u}}\vert^2
$$
$$
+(s\widetilde\varphi)^3\vert d\omega_{\mbox{\bf u}} \vert^2
+ s\widetilde\varphi\vert\nabla' \mbox{div}\,\text{\bf u}\vert^2
+ (s\widetilde\varphi)^3\vert\mbox{div}\,\text{\bf u}\vert^2\Biggr)e^{2s\varphi}
dx'.
$$

Finally we introduce  norms
$$\Vert \text{\bf F}
e^{s\varphi}\Vert^2_{\mathcal Y(\varphi,s,Q)}=\Vert\mbox{div}\, \text{\bf F}
e^{s\varphi}\Vert^2_{L^2(Q)}
+ \Vert d\omega_{\mbox{\bf F}}e^{s\varphi}\Vert^2_{L^2(Q)}+\Vert
(s\widetilde\varphi)^\frac 12\text{\bf F}
e^{s\varphi}\Vert^2_{H^{\frac 14,\frac 12,\widetilde s}(\Sigma)}+\Vert \text{\bf F}
e^{s\varphi}\Vert^2_{L^2(Q)},
$$ and for any positive $p$ we introduce the norm
$$
\Vert \mbox{\bf u}\Vert_{H^{\frac p2, p, \widetilde s}(\Sigma)}=(\Vert \mbox{\bf u}\Vert^2_{H^{\frac p2, p}(\Sigma)}+\Vert (s\widetilde\varphi)^p \mbox{\bf u}\Vert^2_{L^2(\Sigma)})^\frac 12.
$$

Now we are ready to state our first main result, a Carleman estimate
as follows.

\begin{theorem}\label{opa3}
Let $\mbox{\bf F}, \mbox{div}\, \F, d\omega_{\F}\in L^2(Q)$, $\mbox{\bf u}\in  L^2(0,T; H^3(\Omega)),\partial_{x_0}\mbox{\bf u}\in L^2(0,T; H^1(\Omega))$ satisfy (\ref{gora1A}), (\ref{gora2A}).
Let (\ref{pinok0}) and (\ref{pinok1}) hold true, and let a function
$\varphi$ satisfy (\ref{pinok22}) - (\ref{lomka}), Conditions \ref{A1}
and \ref{A2} with $\beta = \mu$ and $\beta=\lambda+2\mu$.
Then there exists  $s_0>0$ such
that  for any  $s > s_0$  the following
estimate holds true:
\begin{eqnarray}\label{2.9'}
\Vert \mbox{\bf u}\Vert_{\mathcal B(\varphi,s,Q)}
+ \Vert\partial^2_{\vec \nu} \text{\bf u}
e^{s\varphi}\Vert_{H^{\frac 14,\frac 12,\widetilde s}( \Sigma_0)}
+\Vert \partial_{\vec \nu}
\text{\bf u}
e^{s\varphi}\Vert_{H^{\frac 34,\frac 32,\widetilde s}( \Sigma_0)}\le C_7(\Vert \text{\bf F}
e^{s\varphi}\Vert_{\mathcal Y(\varphi,s,Q)}
+\\
\Vert \partial_{\vec \nu} \text{\bf u}
e^{s\varphi}\Vert_{H^{\frac 34,\frac 32,\widetilde s}(\widetilde\Sigma)}+ \Vert \partial^2_{\vec \nu} \text{\bf u}
e^{s\varphi}\Vert_{H^{\frac 14,\frac 12,\widetilde s}(\widetilde\Sigma)}+\Vert\widetilde\varphi \partial_{x_0}\partial_{\vec \nu} \text{\bf u}
e^{s\varphi}\Vert_{ L^2(\widetilde\Sigma)}).
\nonumber
\end{eqnarray}
Here the constant $C_7>0$ is independent of $s.$
\end{theorem}

{\bf Example of a function $\varphi$ which satisfies the conditions
of Theorem \ref{opa3}.}
Let $\psi\in C^3(\overline\Omega)$ satisfy
$\psi\vert_{\partial\Omega\setminus\widetilde\Gamma}=0,$ $\vert \nabla '\psi(x)\vert >0$  on $ x'\in\overline\Omega$, $\partial_{\vec \nu}\psi< 0$ on
$\overline{\partial\Omega\setminus\widetilde{\Gamma}}$.
We set
$\varphi(x)=\frac{e^{\lambda \psi(x')}-e^{2\lambda\Vert \psi\Vert_{C^0(\overline\Omega)}}}{\widetilde\ell(x_0)}$, $\widetilde\ell(x_0)>0 $ on $[-T,T]$, $\partial_{x_0}\widetilde\ell(x_0)<0$ on $[0,T)$ and $\partial_{x_0}\widetilde\ell(x_0)>0$ on $[-T,0)$, $\partial^j_{x_0}\widetilde \ell(\pm T)=0$ for all $j\in\{0,1,2\},$
 $\partial^3_{x_0}\widetilde\ell(\pm T)\ne 0.$

Provided that parameter $\lambda>0$ is sufficiently large,  we can prove
that Conditions \ref{A1} and \ref{A2} hold true (e.g.,
H\"ormander \cite{H}, Imanuvilov, Puel and Yamamoto \cite{IPY}).
The normal derivative of the function $\varphi$ on $\Sigma_0$ is strictly
negative and so (\ref{rep}) holds true.
Inequality (\ref{giga1}) follows from the fact that the function
$\psi$ does not have critical points on $\overline\Omega$.
The properties of the function $\widetilde\ell$ imply (\ref{1})-(\ref{lomka}).
%The Carleman estimate (\ref{2.9'}) does not require data on the whole lateral
%boundary $\www\Sigma = (-T,T) \times \www\Gamma$.

The proof  of Theorem \ref{opa3} is given in Sections \ref{QQ1}-\ref{Q2}.
\\
\vspace{0.2cm}

Next we apply Carleman estimate (\ref{2.9'}) to an inverse source problem of
determining a spatially varying factor of source term of the form
$\mathbf{F}(x):= R(x)\f(x')$.
Now we assume that $\rho, \la, \mu$ are independent of $x_0$:
$\rho(x) = \rho(x'), \la(x) = \la(x'), \mu(x) = \mu(x')$ for $x\in (0,T)\times \Omega$.
Let $\OOO \subset \Bbb R^n$ be a smooth bounded domain and $\eta\in (0,T)$ be a fixed time moment.
We consider
\begin{equation}\label{(A1.1)}
\rho(x')\ppp^2_{x_0}\mbox{\bf u}
= L_{\lambda,\mu}(x',D') \ppp_{x_0}\mbox{\bf u} + L_{\www{\lambda},\www{\mu}}
(x,D')\mbox{\bf u}(x)  + R(x)\f(x') \quad
\mbox{in $(0,T)\times \OOO$},
\end{equation}
\begin{equation}\label{(A1.222)}
\mbox{\bf u}(\eta,x') = \mbox{\bf a}(x'), \quad\ppp_{x_0} \mbox{\bf u}(\eta,x') = \mbox{\bf b}(x')\quad
x' \in \OOO
\end{equation}
and
\begin{equation}\label{(A1.3)}
\mbox{\bf u}\vert_{(0,T)\times \ppp\OOO} = 0.
\end{equation}
Here $R(x)$ is an $n\times n$ matrix function and $\f(x')$ is an
$\Bbb R^n$-valued function.

We further assume
\begin{equation}\label{ogonek}
\www{\lambda}, \www{\mu}
\in C^2([0,T]\times \ooo{\OOO}).
\end{equation}

We consider\\
{\bf Inverse source problem.} {\it
Let functions $\mbox{\bf a},\mbox{\bf b},$ $R$ be given and $\widetilde\Gamma$ is an arbitrary fixed open subset of  $\partial\OOO.$
Determine function $\f(x')$ by $\ppp_{\vec\nu}\mbox{\bf u}\vert
_{(0,T)\times \widetilde\Gamma}$ and $\mbox{\bf u}(\eta, \cdot),
\ppp_{x_0}\mbox{\bf u}(\eta,\cdot)$ in $\Omega$.}
\\

We state our main result on the inverse source problem.
\\
\begin{theorem}\label{theorem 1.2}{\it
Let $\mbox{\bf a},\mbox{\bf b}\in H^3(\Omega),$  $\f \in H^1(\OOO),$ $\mbox{\bf u}\in  L^2(0,T; H^3(\Omega)),\partial^k_{x_0}\mbox{\bf u}\in L^2(0,T; H^1(\Omega)), k\in \{1,2,3\}$ satisfy (\ref{(A1.1)}) and (\ref{(A1.222)}).
We assume that there exists a constant $\delta_0>0$ such that
\begin{equation}\label{(1.18)}
\vert \mbox{det} \thinspace R(\eta,x') \vert \ge \delta_0 > 0,
 \quad x' \in \ooo{\OOO},\quad R\in C^1([0,T]\times\overline\Omega).
\end{equation}
Moreover we assume that the Lam\'e coefficients $\lambda,\mu,\rho,
\widetilde \mu,
\widetilde \lambda$ satisfy (\ref{ogonek}) and (\ref{pinok0}).
Then there exists a constant $C_8>0$ such that
\begin{eqnarray}\label{(1.19)}
\Vert \f\Vert_{H^1(\OOO)}
\le C_8\biggl( \Vert \mbox{\bf a}\Vert_{H^3(\Omega)}+\Vert \mbox{\bf b}\Vert_{H^3(\Omega)}+\sum_{j=0}^2 (\Vert \partial^j_{x_0}(d \omega_{\mbox {\bf u}},\mbox{div}\,\mbox {\bf u})
\Vert_{H^{\frac 34,\frac 32}((0,T)\times\widetilde\Gamma)}\nonumber\\
+ \Vert \partial^j_{x_0}\partial_{\vec \nu}(d \omega_{\mbox {\bf u}},\mbox{div}\,\mbox {\bf u})
\Vert_{H^{\frac 14,\frac 12}((0,T)\times\widetilde\Gamma)}+\Vert\partial^{j+1}_{x_0}\partial_{\vec \nu} \text{\bf u}
\Vert_{ L^2((0,T)\times\widetilde\Gamma)} \biggr).
\end{eqnarray}
}
\end{theorem}

The proof of Theorem \ref{theorem 1.2} is provided in Section \ref{Q1}.

There are other works on inverse problems related to the viscoelasticity
and we refer to Cavaterra, Lorenzi and Yamamoto \cite{CLY},
de Buhan \cite{dB}, de Buhan and Osses \cite{dBO},
Grasselli \cite{G}, Imanuvilov and Yamamoto \cite{IY2018},
Janno \cite{J}, Janno and von Wolfersdorf \cite{JW},
Lorenzi, Messina and Romanov \cite{LMR}, Lorenzi and Romanov \cite{LR},
Loreti, Sforza and Yamamoto \cite{LSY},
Romanov and Yamamoto \cite{RY}, von Wolfersdorf \cite{W}.
\\
\vspace{0.3cm}

Finally, as application of Carleman estimate Theorem \ref{opa3} in a special
case of $\www{\la} = \www{\mu} \equiv 0$, we consider an inverse source problem
for an isothermal compressible fluid equations (e.g.,  Landau and Lifshitz
\cite{LL}), which can be simplified as follows:
\begin{equation}\label{Z1A}
\partial_{x_0}\rho+\mbox{div}(\mbox{\bf v}\rho)=0\quad \mbox{in} \,Q=(-T,T)
\times \Omega,
\end{equation}
\begin{equation}\label{Z2A}
\rho\partial_{x_0}\mbox{\bf v}
= L_{\lambda,\mu}(x',D')\mbox{\bf v} - \rho(\mbox{\bf v},\nabla')\mbox{\bf v}
- h(\rho)\nabla' \rho + R\mbox{\bf f}\quad \mbox{in} \,Q,
\end{equation}
\begin{equation}\label{Z3A}
\mbox{\bf v}\vert_{\Sigma}=0.
\end{equation}

Here $\rho(x)$ and $\mbox{\bf v}(x)$ describe the density and the
velocity field at $x \in Q$ respectively.

Consider the following inverse problem:
{\it  Let $\rho(0,\cdot), \mbox{\bf v}(0,\cdot)$ in $\OOO$ and
$\partial_{\vec\nu}\mbox{\bf v}$ on $(-T,T)\times \Gamma$ be
given.  Then determine a function $\mbox{\bf f}(x').$}

Then we have

\begin{theorem}\label{theorem 1.3}
{\it Let $M_1>0$ and $M_2>0$ be constants, $h\in C^2(\Bbb R^1)$ and
for some constant $\delta_0 > 0$, let the function $R$ satisfy
\begin{equation}\label{grabesz}
R\in C^1(\overline Q), \quad \vert \mbox{det} \thinspace R(0,x') \vert \ge
\delta_0 > 0, \quad x' \in \ooo{\OOO}.
\end{equation}
Moreover for $k=1,2$, we assume that there exist constants $M_1>0$ and
$M_2>0$ such that the functions  $(\rho_k, \mbox{\bf v}_k)
\in C^3(\ooo{Q})\times C^3(\ooo{Q})$ satisfy (\ref{Z1A})- (\ref{Z3A})  with function $\f_k\in H^1_0(\Omega)$ and
\begin{equation}\label{compot1}
\Vert (\rho_k, \mbox{\bf v}_k)\Vert_{C^3(\overline Q)\times C^3(\overline Q)}
\le M_1
\end{equation}
and
\begin{equation}\label{compot}
\rho_k(0,x')\ge M_2>0\quad \forall x'\in \overline\Omega,\,\, \forall
k\in \{1,2\}.
\end{equation}
Furthermore let $\la = \la(x')$ and $\mu=\mu(x')$ satisfy
$\la, \mu \in C^2(\overline{\Omega})$ and
$\la + \mu, \mu > 0$ on $\overline{\OOO}$.

Then we can choose a constant $C_9$ independent of $(\rho_k , \mbox{\bf v}_k)$
such that
\begin{eqnarray}\label{(1.19)}
\Vert \f_1-\f_2\Vert_{H^1(\OOO)}
\le C_9( \Vert  (\rho_1-\rho_2)(0,\cdot)\Vert_{H^3(\Omega)}+\Vert  (\mbox{\bf v}_1-\mbox{\bf v}_2)(0,\cdot)\Vert_{H^3(\Omega)}\\+\sum_{j=0}^1 (\Vert \partial^j_{x_0}(d \omega_{\mbox {\bf v}_1-\mbox {\bf v}_2},\mbox{div}\,(\mbox {\bf v}_1-\mbox {\bf v}_2))
\Vert_{H^{\frac 34,\frac 32}(\widetilde\Sigma)}\nonumber\\
+ \Vert \partial^j_{x_0}\partial_{\vec \nu}(d \omega_{\mbox {\bf v}_1-\mbox {\bf v}_2},\mbox{div}\,(\mbox {\bf v}_1-\mbox {\bf v}_2))
\Vert_{H^{\frac 14,\frac 12}(\widetilde\Sigma)}+\Vert\partial^{j+1}_{x_0}\partial_{\vec \nu}( \text{\bf v}_1-\mbox {\bf v}_2)
\Vert_{ L^2(\widetilde\Sigma)} )).\nonumber
\end{eqnarray}                         }
\end{theorem}
\section{Proof of Theorem \ref{opa3}}\label{QQ1}
Sections \ref{QQ1}-\ref{Q1} are devoted to the proof of Theorem \ref{opa3} and in Section \ref{Q!8}
we collect necessary lemmata for the proof.

%Next we note that instead of (\ref{2.9'}), it suffices to prove a simpler
%estimate
%\begin{eqnarray}\label{3.1}
%\int_Q( s\tau^2\varphi\vert\nabla d\omega_{\mbox{\bf u}}\vert^2
%+s^3\tau^4\varphi^3\vert d\omega_{\mbox{\bf u}}\vert^2
%+s\tau^2\varphi\vert\nabla\text{div}\,\mbox{\bf u}\vert^2
%+s^3\tau^4\varphi^3\vert\text{div}\,\mbox{\bf u}\vert^2)
%e^{2s\varphi} dx\nonumber\\
%+s\tau\Vert\varphi^\frac 12 \nabla\partial_{\vec \nu} \text{\bf u}
%e^{s\varphi}\Vert_{{L}^2(\Sigma)}^2
%+ s\tau\Vert \varphi^\frac 12\partial_{\vec \nu}^2  \text{\bf u}e^{s\varphi}
%\Vert^2_{{L}^2(\Sigma)}\\
%\le C_7(\Vert \text{\bf F}
%e^{s\varphi}\Vert^2_{\mathcal Y(\psi,s,\tau,Q)}+ s\tau\Vert\varphi^\frac 12
%\nabla\partial_{\vec \nu} \text{\bf u}
%e^{s\varphi}\Vert_{{ L}^2(\widetilde\Sigma)}^2
%+ s\tau\Vert \varphi^\frac 12\partial_{\vec \nu}^2  \text{\bf u}e^{s\varphi}
%\Vert^2_{{ L}^2(\widetilde\Sigma)}+s^3\tau^5\Vert\varphi^\frac 52  \text{\bf u}
%e^{s\varphi}\Vert_{{L}^2(\widetilde\Sigma)}^2)\nonumber.
%\end{eqnarray}
%This fact is the simple corollary of the following proposition:

Consider the following boundary value problem
\begin{equation}\label{poko}
\mbox{\bf P}(x,D) \mbox{\bf v}\triangleq  (
P_{\rho,\mu} (x,D)\mbox{\bf v}_1, P_{\rho,\lambda+2\mu} (x,D) v_2)=\mbox{\bf q} \quad \mbox{in}\,\, Q.
\end{equation}

Here $ P_{\rho,\beta}(x,D)=\rho\partial^2_{x_0}-\beta \Delta,$ $\mbox{\bf q}$ is a given vector-valued function from $L^2(Q)$
with $\frac{n^2-n}{2}+1$ components, and a vector-valued function
$\mbox{\bf v}_1$ with $(n^2-n)/2$ components is defined by
$$
\mbox{\bf v}_1=(v_{1,2}, \dots  v_{n-1,n})
= (v_{k,j})_{1\le k<j\le n}.
$$
We set $\mbox{\bf v} = (\mbox{\bf v}_1, v_2)$.
Function $ \mbox{\bf v}$ satisfies the boundary conditions
\begin{equation}\label{poko1}
\mathcal B(x,D')\mbox{\bf v}=\mbox{\bf g}\quad \mbox{on}\,\Sigma.
\end{equation}
where $\mbox{\bf g}=(g_1,\dots, g_{\frac{n^2+n}{2}+1})$ is
a given vector-valued function and the boundary operator $\mathcal B(x,D')$ is
constructed in the following way
$$
\mathcal B(x,D')= (B_1(x,D'), B_2(x')),\quad
B_1(x,D')=(b_{1}(x,D'),\dots ,b_{n}(x,D'))
$$
and
$$
b_{k}(x,D')\mbox{\bf v}= -\sum_{j=1, j\ne k}^n\mbox{sign}(k-j)
\partial_{x_j} v_{j,k}-\frac{(\lambda+2\mu)(x)}{\mu(x)}\partial_{x_{k}}v_2
$$
for $1\le k \le n$, and $B_2(x')$ is the smooth matrix constructed in the
following way: Consider an $n\times n$ matrix such that on the main diagonal
we have $\nu_n(x')$, the $n$th row is $(\nu_1(x'),\dots, \nu_n(x'))$,
the first $n-1$ elements of the last column are $-\nu_1(x'),\dots,
-\nu_{n-1}(x')$ and all the rest elements are zero:
$$
\left (\begin{matrix}
\nu_n(x') &0&\dots& 0&\dots & -\nu_1(x')\\
0& \nu_n(x')&\dots &\dots\\
0& 0&\dots &\nu_n(x')&\dots&-\nu_j(x')\\
\dots&\dots&\dots &\dots\\ \nu_1(x')&\dots &\dots&\dots& \nu_{n-1}(x')
&\nu_n(x') \end{matrix}\right).
$$
If $\nu_n(x')\ne 0 $, then the determinant of such a matrix is not equal
to zero.
Denote the inverse to this matrix by $B_3(x')$
and set $\mbox{\bf w}:= (w_1, ..., w_n) = B_3\widetilde{ \mbox{\bf v}}$ where
$\widetilde{\mbox{\bf v}}= (\www{v}_1, ..., \www{v}_n)$,
$\www{v}_j = v_{j,n}$ for $j \in \{1,..., n-1\}$ and $\www{v}_n
= v_2$.
Then $B_2(x')\mbox{\bf v}=\mbox{\bf v}-(\nu_2w_1-\nu_1 w_2, \dots ,
\nu_n w_1-\nu_1 w_n, \dots , \nu_n w_{n-1}-\nu_{n-1}w_n,
\sum_{j=1}^n \nu_jw_j)$.

We are not investigating the existence of solution for problem
(\ref{poko}) - (\ref{poko1}), but assume the existence of  function $\mbox{\bf v}$ which  satisfies
(\ref{poko}) - (\ref{poko1}) with some $\mbox{\bf g}$.
Our current goal is to establish the Carleman estimate for solution of this
problem.

Without loss of generality  we can assume that
\begin{equation}\label{normal}
\vec \nu(0)=-\vec e_n=(0,\dots,0,-1).
\end{equation}

We have

\begin{proposition}\label{zoopa} Let $\mbox{\bf v}\in H^1(Q)$ satisfy (\ref{poko}), (\ref{poko1}).
There exists  $s_0>1$  such that
 for all
$s\ge s_0$ the following estimate  holds true
\begin{eqnarray}\label{3.2'1'}
\sum_{\vert\alpha\vert\le 2}\Vert (s\widetilde\varphi)^\frac {(3-2\vert\alpha\vert)}{2} \partial^\alpha  \mbox{\bf v}\,
e^{s\varphi}\Vert_{L^2(Q)}  +\Vert \mbox{\bf v} e^{s\varphi}\Vert_{H^{\frac 34,\frac 32,\widetilde s}(\Sigma)}\\
\le C_1(\Vert \mbox{\bf P}(x,D)\mbox{\bf v} e^{s\varphi}\Vert_{L^2(Q)}
+\Vert \mbox{\bf g}e^{s\varphi}\Vert_{H^{\frac 14,\frac 12,\widetilde s}(\Sigma)}
+\Vert \mbox{\bf v} e^{s\varphi}\Vert_{H^{\frac 34,\frac 32,\widetilde s}(\widetilde\Sigma)}+\Vert \partial_{\vec \nu}\mbox{\bf v} e^{s\varphi}\Vert_{H^{\frac 14,\frac 12,\widetilde s}(\widetilde\Sigma)}),\nonumber
\end{eqnarray}
where $C_1$ is independent of $s.$
\end{proposition}
The end of this section, sections \ref{QQ1}-\ref{QQ5}, and part of section \ref{Q2} are devoted to the proof of Proposition \ref{zoopa}.

First, by an argument based on the partition of unity (e.g., Lemma 8.3.1
in \cite{H}), it suffices to prove the inequality (\ref{3.2'1'})
locally, by assuming that
\begin{equation}\label{3.55}
\text{supp}\, \mbox{\bf v}\subset B (y^*,\delta),
\end{equation}
where $B(y^*,\delta)$ is the ball of the radius $\delta>0$
centered at some point $y^*$.

Otherwise, without the loss  of generality we may assume that $y^*=(y^*_0,0,
\dots, 0).$
Let $\theta\in C^\infty_0(\frac 12,2)$ be a nonnegative function such that
\begin{equation}\label{knight}\sum_{\ell=-\infty}^\infty\theta(2^{-\ell}t)=1
\quad  \mbox{for all}\,\, t\in \Bbb R^1.
\end{equation}

Set $\mbox{\bf v}_\ell(x)=\mbox{\bf v}(x)\kappa_\ell(x_0),\mbox{\bf g}_\ell(x)=\mbox{\bf g}(x)\kappa_\ell(x_0)$ where
\begin{equation}\label{knight0}
\kappa_\ell(x_0)=
\theta \left (2^{-\ell}2^\frac{1}{\theta(x_0)^\frac 14}\right ),
\end{equation} where
\begin{eqnarray}\label{Aziza}\theta(x_0)\in C^\infty[-T,T],\quad \theta\vert_{[-T,T/2]}=T+x_0,\quad \theta\vert_{[T/2,T]}=T-x_0,\nonumber\\ \partial_{x_0}\theta<0 \,\,\mbox{on}\,\, (-T,0), \quad \partial_{x_0}\theta>0 \,\,\mbox{on}\,\, (0,T),\quad \partial^2_{x_0}\theta(0)<0 .
\end{eqnarray}
Observe that it
suffices to prove the Carleman estimate (\ref{3.2'1'}) for the function
$\mbox{\bf v}_\ell$ instead of $\mbox{\bf v}$ provided that the constant
$C_1$  and the function $s_0$ are independent of $\ell.$
Observe that if $G\subset \Bbb R^m$ is a bounded domain and $g\in L^2(G)$,
then there exist an independent constants $C_2$ and $C_3$
(see e.g. \cite{Sog}) such that
\begin{equation}\label{gorokn}
C_2\sum_{\ell=-\infty}^{\infty}\Vert \kappa_\ell g\Vert^2_{L^2(G)}
\le \Vert g\Vert^2_{L^2(G)}\le C_3\sum_{\ell=-\infty}^{\infty}\Vert \kappa_
\ell g\Vert^2_{L^2(G)}.
\end{equation}
Denote the norm on the left-hand side of (\ref{3.2'1'}) as $\Vert\cdot\Vert_*.$ Suppose that the estimate  (\ref{3.2'1'}) is true for any function $\mbox{\bf v}_\ell$ with constants $C_1$ and $s_0$ independent of $\ell.$
By (\ref{gorokn})  for some constant $C_4$ independent of $s$ we have

\begin{equation}\label{zombi}
\Vert \mbox {\bf v} e^{s\varphi}\Vert_*
=\Vert \sum_{\ell=-\infty}^{+\infty} \mbox{\bf v}_\ell e^{s\varphi}\Vert_*\le
\sum_{\ell=-\infty}^{+\infty}\Vert \mbox {\bf v}_\ell e^{s\varphi}\Vert_*\le
C_4\sum_{\ell=-\infty}^{\infty}(\Vert \kappa_\ell \mbox{\bf P}(x,D)
\text{\bf v}
e^{s\varphi}\Vert^2_{L^2(Q)}
\end{equation}
$$+\Vert e^{s\varphi}[\kappa_\ell,\mbox{\bf P}(x,D)]
\text{\bf v}
\Vert^2_{L^2(Q)}
+ \Vert \mbox{\bf v}_\ell e^{s\varphi}\Vert^2_{H^{\frac 34,\frac 32,\widetilde s}(\widetilde\Sigma)}+\Vert \partial_{\vec \nu}\mbox{\bf v}_\ell e^{s\varphi}\Vert^2_{H^{\frac 14,\frac 12,\widetilde s}(\widetilde\Sigma)}+\Vert \kappa_\ell\mbox{\bf g}e^{s\varphi}\Vert^2_{H^{\frac 14,\frac 12,\widetilde s}(\widetilde\Sigma)})^\frac 12.
$$

Assume that near $(0,\dots,0)$,
the boundary $\partial \Omega$ is locally given by an equation
$x_n-\vartheta(x_1,\dots, x_{n-1})=0$ and
if $(x_1, \dots, x_n) \in \Omega$, then $x_n -\vartheta (x_1,\dots, x_{n-1}) > 0$
where $\vartheta\in C^3$
and $\vartheta(0) = 0$.  Since $\vec\nu(0)=-\vec e_n$
\begin{equation}\label{napoleon}
(\partial_{x_1}\vartheta(0),\dots \partial_{x_{n-1}}\vartheta(0))=0.
\end{equation}
%The operator $\widetilde L_{\lambda,\mu}(y,D')$ given by formula
%\begin{equation}\label{gramophone}
%\widetilde L_{\lambda,\mu}(y,D')\mbox{\bf u}=\mu \Delta_\ell \mbox{\bf u}
%- (\lambda+\mu)\nabla\left(\text{div}\,\mbox{\bf u}
%- (\partial_{y_n} \mbox{\bf u},\nabla \ell)\right)
%-(\lambda+\mu)\partial_{y_n}\left(
%\text{div}\,\mbox{\bf u}
%- (\partial_{y_n} \mbox{\bf u},\nabla \ell)\right)\nabla\ell,
%\end{equation}
Denote  \begin{equation}\label{kupol}F(x)=(x_0,\dots, x_{n-1}, x_n-\vartheta(x_1,\dots, x_{n-1})).\end{equation}
By Proposition \ref{Fops5} we obtain from (\ref{zombi}):
\begin{equation}\label{ozon}
\Vert \mbox {\bf v} e^{s\varphi}\Vert_*
=\Vert \sum_{\ell=-\infty}^{+\infty} \mbox{\bf v}_\ell e^{s\varphi}\Vert_*\le
\sum_{\ell=-\infty}^{+\infty}\Vert \mbox {\bf v}_\ell e^{s\varphi}\Vert_*\le
C_5(\Vert \mbox{\bf P}(x,D)
\text{\bf v}
e^{s\varphi}\Vert^2_{L^2(Q)}
\end{equation}
$$+\sum_{\ell=-\infty}^{\infty}\Vert e^{s\varphi}[\kappa_\ell,\mbox{\bf P}(x,D)]
\text{\bf v}
\Vert^2_{L^2(Q)}
+ \Vert \mbox{\bf v} e^{s\varphi}\Vert^2_{H^{\frac 34,\frac 32,\widetilde s}(\widetilde\Sigma)}+\Vert \partial_{\vec \nu}\mbox{\bf v} e^{s\varphi}\Vert^2_{H^{\frac 14,\frac 12,\widetilde s}(\widetilde\Sigma)}+\Vert \mbox{\bf g}e^{s\varphi}\Vert^2_{H^{\frac 14,\frac 12,\widetilde s}(\widetilde\Sigma)})^\frac 12.
$$

Using (\ref{knight0}) and (\ref{Aziza}) we
estimate the norm of the commutator $ [\kappa_\ell,\mbox{\bf P}(x,D)]$ we obtain
\begin{eqnarray}\label{ooo}
\sum_{\ell=-\infty}^{\infty}\Vert e^{s\varphi} [\kappa_\ell,\mbox{\bf P}(x,D)]\text{\bf v}
\Vert^2_{L^2(Q)}\le C_{6}\sum_{\ell=-\infty}
^{\infty}(\Vert \partial_{x_0}\kappa_\ell\nabla\text{\bf v}
e^{s\varphi}\Vert^2_{L^2(Q)}
+\Vert\partial^2_{x_0}\kappa_\ell \text{\bf v}
e^{s\varphi}\Vert^2_{L^2(Q)})\nonumber\\
\le C_{7}\sum_{\ell=-\infty}
^{\infty}(\Vert \widetilde \varphi^\frac {5}{12} \chi_{\mbox{supp}\, \kappa_\ell}\nabla\text{\bf v}
e^{s\varphi}\Vert^2_{L^2(Q)}
+\Vert \widetilde \varphi^\frac 32 \chi_{\mbox{supp}\, \kappa_\ell}\text{\bf v}
e^{s\varphi}\Vert^2_{L^2(Q)})\nonumber\\
\le C_{8}(\Vert \widetilde \varphi^\frac 12 \nabla\text{\bf v}
e^{s\varphi}\Vert^2_{L^2(Q)}
+\Vert \widetilde \varphi^\frac 32 \text{\bf v}
e^{s\varphi}\Vert^2_{L^2(Q)})
.
\end{eqnarray}

From  (\ref{ozon}), and (\ref{ooo}) we obtain (\ref {3.2'1'}).

Now, without loss of generality we assume that
\begin{equation}\label{3.55}
\text{supp}\, \mbox{\bf v}\subset B (y^*,\delta)\cap
\mbox{supp}\, \kappa_\ell,
\end{equation}
where $B(y^*,\delta)$ is the ball of the radius $\delta>0$
centered at some point $y^*=(y_0^*,0,\dots,0).$

We set
$$
\Delta_{\vartheta} u=\sum_{j=1}^{n-1}(\partial^2_{y_j} u
- 2\partial_{x_j}\vartheta\circ F^{-1}(y)\partial^2_{y_jy_n} u )
+ (1+\vert \nabla'\vartheta\vert^2\circ F^{-1}(y))\partial^2_{y_n} u.
$$
Henceforth we set $y = (y_0, y ')=(y_0,y_1,\dots, y_n)$.
After the change of variables, the  equations (\ref{poko}) have the forms
\begin{equation}\label{3.9}
\mbox{\bf P}(y,D) \mbox{\bf v}=(  \rho\partial_{y_0} \mbox{\bf v}_1-\mu \Delta_{\vartheta} \mbox{\bf v}_1,\rho\partial_{y_0} v_2
- (\lambda +2\mu)\Delta_{\vartheta} v_2)
= \mbox{\bf q}, \quad
\mbox{on}\, \mathcal Q\triangleq {\Bbb R}^n\times [0, \gamma] ,    \end{equation}
\begin{equation}\label{!poko1}
\widetilde {\mathcal B}(y,D)\mbox{\bf v}=\mbox{\bf g},
\end{equation}
where $\gamma$ is some positive constant.  Without loss of generality, we may assume $\gamma=1$.
%and
%\begin{eqnarray}\label{3.9'}
%P_{\lambda+2\mu}(y,D) w_2 \equiv \rho\partial^2_{y_0} w_2
%- (\lambda +2\mu)\Delta_\ell w_2  = q_4
%\qquad\text{in}\,\mathcal  Q .
%\end{eqnarray}
Here    for  functions $\rho\circ F^{-1}(y), \mu\circ F^{-1}(y),
\lambda \circ F^{-1}(y)$ we used the notations $\rho,\mu,\lambda$.
Similarly by  $\mbox{\bf q}_{3}, q_4$ denote the functions $\mbox{\bf q}_3, q_4$ after the
change of variables:  $\mbox{\bf q}_3\circ F^{-1}(y) , q_4\circ F^{-1}(y).$

The operator $\widetilde {\mathcal B}(y,D)$ is obtained from ${\mathcal B}(x,D)$ in the following way
$$
\widetilde{\mathcal B}(y,D')= (\widetilde B_1(y,D'), B_2(F^{-1}(y))),\quad
\widetilde B_1(y,D')=(\widetilde b_{1}(y,D'),\dots ,\widetilde b_{n}(y,D'))
$$
and
$$
\widetilde b_{k}(y,D')\mbox{\bf v}= -\sum_{j=1, j\ne k}^n\sum_{m=1}^n\mbox{sign}(k-j)
\partial_{x_m} v_{j,k}\partial_{x_j}F_m^{-1}\circ F^{-1}(y)-\frac{(\lambda+2\mu)}{\mu}\sum_{m=1}^n\partial_{x_m} v_{2}\partial_{x_k}F_m^{-1}\circ F^{-1}(y).
$$
Since $F'(y^*)$ is the unit matrix, from the above equality we have
\begin{equation}\label{legioner}
\widetilde b_{k}(y^*,D')\mbox{\bf v}= -\sum_{j=1, j\ne k}^n\mbox{sign}(k-j)
\partial_{y_j} v_{j,k}-\frac{(\lambda+2\mu)}{\mu}(y^*)\partial_{y_{k}}v_2.
\end{equation}
Now we introduce operators
\begin{equation}\label{3.18}
P_{\rho,\mu}(y,D,\widetilde s) = e^{\vert s\vert\varphi}P_{\rho,\mu}(y,D)
e^{-\vert s\vert\varphi},
\quad
P_{\rho,\lambda+2\mu}(y,D,\widetilde s)
= e^{\vert s\vert\varphi}P_{\rho,\lambda+2\mu} (y,D)
e^{-\vert s\vert\varphi},
\end{equation}
where
$$
\widetilde s=s\widetilde\varphi(y^*_0).
$$

We denote the principal symbols
of the operators $P_{\rho,\mu}(y,D,\widetilde s)$ and $P_{\rho,\lambda+2\mu}(y,D,\widetilde s)$ by \newline
$p_{\rho,\mu}(y,\xi,\widetilde s)$ $=p_{\rho,\mu}(y,\xi+{i}\vert s\vert \nabla\varphi)$
and $p_{\rho,\lambda+2\mu}(y,\xi,\widetilde s)=p_{\rho,\lambda+2\mu}(y,\xi
+{i}\vert s\vert \nabla\varphi)$ respectively.

The principal symbol of the operator $P_{\rho,\beta}(y,D,\widetilde s)$ has the
form
\begin{eqnarray}\label{3.34}
p_{\rho,\beta}(y,\xi,\widetilde s) = i\rho(y)(\xi_0+{i}\vert s\vert\varphi_{y_0})
+ \beta\biggl[
\sum_{j=1}^{n-1}(\xi_j+{i}\vert s\vert\varphi_{y_j})^2
                                               \nonumber\\
-2(\nabla'\vartheta, (\xi'+{i}\vert s\vert\nabla'\varphi))
(\xi_n+{i}\vert s\vert\varphi_{y_n})
+(\xi_n+{i}\vert s\vert\varphi_{y_n})^2 G\biggr],
\end{eqnarray}
where $ G(y_1,\dots,y_{n-1})=1+\vert\nabla\vartheta(y_1,\dots,y_{n-1})\vert^2.$
The zeros of the polynomial $p_{\beta}(y,\xi,\widetilde s)$ with respect to variable $\xi_n$  for
$M(\widetilde \xi, \widetilde s) \ge 1, \widetilde \xi=(\xi_1,\dots,\xi_{n-1})$ and $y\in B (y^*,\delta)\cap
\mbox{supp}\, \kappa_\ell$ are
\begin{equation}\label{3.35}
\Gamma^\pm_\beta(y,\widetilde \xi,\widetilde s)
= (-{i}\vert \widetilde s\vert\widetilde\mu_\ell\varphi_{n}\kappa(\widetilde\xi,\widetilde s)
+\alpha^\pm_\beta(y,
\widetilde \xi,\widetilde s)),
\end{equation}
where $\vec \varphi=(\varphi_0,\dots \varphi_n),\varphi_j(y)= \frac{\varphi_{y_j}(y)}{\widetilde\varphi(y^*)}
$,
\begin{equation}\label{bobo}\widetilde\mu_\ell(y)=\eta_*(y)\sum_{k=\ell-30}^{\ell+30}\kappa_\ell(y_0),\quad
\eta_*\in C_0^\infty (B(y^*,\frac{98}{50}\delta)),\quad \eta_*\vert_{B(y^*,\delta)}=1,\end{equation}
the function $\kappa_\ell$ is given by (\ref{knight0}),
\begin{equation}\label{3.36}
\alpha^\pm_\beta(y,\widetilde \xi,\widetilde s)=\widetilde \mu_\ell(y)
\left(\frac{-\sum_{j=1}^{n-1}(\xi_j+{i}\vert \widetilde s\vert
\varphi_j)\partial_{y_j}\vartheta(y_1,\dots,y_{n-1})}{\vert G\vert}\kappa(\widetilde\xi,\widetilde s)\pm\root
\of{r_\beta(y,\widetilde\xi,\widetilde s)}\right),                          \end{equation}
\begin{equation}\label{3.37}
r_\beta(y,\widetilde \xi,\widetilde s)=\kappa^2(\widetilde\xi,\widetilde s)
\frac{(-i\rho\xi_0-\beta\sum_{j=1}^{n-1}(\xi_j+{i}
\vert \widetilde s\vert\varphi_j)^2)
G+\beta(\xi+{i}\vert \widetilde s\vert\vec {\varphi},\nabla\vartheta)
^2}
{\beta G^2},
\end{equation}
where

Let
$\chi_\nu$ be a $C^\infty_0(\Bbb M)$ function on
$\Bbb M=\{(\widetilde \xi,\widetilde s); \thinspace
M(\widetilde \xi, \widetilde s)=1\}$ such that
$\chi_\nu$ is identically equal $1$ in some neighborhood of the $(\widetilde \xi^*,s^*)\in
\Bbb M$ and $\mbox{supp}\, \chi_\nu\subset \mathcal O(\zeta^*,\delta).$
Assume that
\begin{equation}\label{book}
\kappa(\widetilde\xi,\widetilde s)\vert_{\mbox{supp}\, \chi_\nu}=1,\quad \mbox{supp}\, \kappa(\widetilde\xi,\widetilde s)
\subset \mathcal O(\zeta^*,2\delta), \,\,\kappa(\widetilde\xi,\widetilde s)\ge 0\quad\mbox{on}\,\, \Bbb M.
\end{equation}

 We extend the function $\chi_\nu$ on $\Bbb R^{n+1}$ as follows :
$\chi_\nu(\xi_0/M^2(\widetilde\xi, \widetilde s), \xi_1/M(\widetilde\xi, \widetilde s),\dots,  \xi_{n-1}/M(\widetilde\xi, \widetilde s))$ for
$M(\widetilde\xi,\widetilde s)>1$ and
$\chi_\nu(\xi_0/M^2(\widetilde\xi, \widetilde s), \xi_1/M(\widetilde\xi, \widetilde s),\dots,  \xi_{n-1}/M(\widetilde\xi, \widetilde s))\kappa^*(
M(\widetilde\xi,s))$  for $M(\widetilde\xi, \widetilde s)
<1$, where $\kappa^*(t)\in C^\infty(\Bbb R^1), \kappa^*(t)\ge 0, \kappa^*(t)=1$ for
$t\ge 1$ and $\kappa^*(t)=0$ for $t\in [0,1/2].$ In the similar way we extend the function $\kappa(\widetilde \xi,s)$ on $\Bbb R^{n+1}$. Denote by $\chi_\nu(y,D',\widetilde s)$ the pseudodifferential operator
with the symbol $\eta_\ell(y)\chi_\nu(\widetilde \xi,\widetilde s)$ and $\eta_\ell(y)=
\sum_{k=-10}^{10}\kappa_{\ell+k}.$
We set $\mbox{\bf w}_\nu
=\chi_\nu(y,\widetilde D,\widetilde s)\mbox{\bf w}$
and $\mbox{\bf w}=\mbox{\bf v}e^{s\varphi}.$

Let ${\mathcal O}$ be a domain in ${\Bbb R}^n.$

{\bf Definition.} {\it We say that the symbol
$a(\widetilde y,\widetilde \xi,\widetilde  s)\in W^{\widetilde k,\infty}(\overline{\mathcal O}
\times
{\Bbb R}^{n+1})$ belongs to the class
$W^{\widetilde k,\infty}_{cl}S^{\kappa,\widetilde  s}({\mathcal O})$ if

{\bf A}) There exists a compact set $K\subset\subset{\mathcal O}$
such that $a(\widetilde y,\widetilde \xi,\widetilde s)\vert_{{\mathcal O}
\setminus K}=0;$

{\bf B}) For any $\beta=(\beta_0,\dots,\beta_{n})$ there exists a
constant $C_\beta$
$$
\left\Vert \partial^{\beta_0}_{\xi_0}\cdots
\partial^{\beta_{n-1}}_{\xi_{n-1}}
\partial^{\beta_{n}}_{\widetilde s}
a(\cdot,\widetilde\xi,\widetilde s)\right\Vert_{W^{\widetilde k,\infty}({\mathcal O})}\le
C_\beta\left(\widetilde s^2+\vert \xi_0\vert+\sum_{i=1}^{n-1}\xi^2_i\right)
^{\frac{\kappa-\vert\beta\vert}{2}}\quad,
$$
where $\vert \beta\vert=\sum_{j=0}^{n}\beta_j$ and
$M(\widetilde\xi,\widetilde s) \ge 1$;

{\bf C}) For any $N\in \Bbb N$ the symbol $a$ can be represented
as
$$
a(\widetilde y,\widetilde \xi,\widetilde s)=\sum_{j=1}^Na_j(\widetilde y,
\widetilde \xi,\widetilde s)
+ R_N(\widetilde y,\widetilde \xi,\widetilde s)
$$
where the functions $a_j$ have the following properties: for any $\lambda>1$ and for all $(\widetilde y,\widetilde \xi,\widetilde  s)\in\{(\widetilde y,
\widetilde \xi,\widetilde  s)
\vert \widetilde y\in K, \,M (\widetilde\xi,\widetilde s) >1\}$
$$
a_j(\widetilde y,{\lambda}^2 \xi_0,\lambda\xi_1, \dots, \lambda \xi_{n-1},{\lambda}
\widetilde s)={\lambda}^{\kappa-j}a_j(\widetilde y,\widetilde \xi,\widetilde  s);
$$
for any multiindex $\beta$ and any and $(\widetilde \xi,\widetilde  s)$ satisfying
$M(\widetilde\xi,\widetilde s)\ge 1$ there exist a constant $C_\beta$ such that
$$
\left\Vert  \partial^{\beta_1}_{\xi_0}\cdots
\partial^{\beta_{n-1}}_{\xi_{n-1}}
\partial^{\beta_{n}}_{\widetilde s}
a_j(\cdot,\widetilde \xi,\widetilde s)\right\Vert_{W^{\widetilde k,\infty}({\mathcal O})}\le
C_\beta\left(\widetilde s^2+\vert \xi_0\vert+\sum_{i=1}^{n-1}\xi^2_i\right)
^{\frac{\kappa-j-\vert\beta\vert}{2}}
$$
where the term $R_N$ satisfies the estimate
$$
\Vert R_N(\cdot,\widetilde\xi,\widetilde s)\Vert_{W^{\widetilde k,\infty}({\mathcal O})}\le
C_N(\widetilde s^2+\vert \xi_0\vert+\sum_{i=1}^{n-1}\xi^2_i)^{\frac{\kappa-N}{2}}\quad
\forall (\widetilde\xi,\widetilde  s)\,\, \mbox{satisfying}\,\, M
(\widetilde\xi,\widetilde s)
 \ge 1.
$$}

Let $X^{\widetilde k}(\mathcal O)=W^{1,\infty}(\mathcal O)$ or  $X^{\widetilde k}(\mathcal O)=C^{\widetilde k}(\overline{\mathcal O}).$
For the symbol $a$, we introduce the semi-norm
\begin{eqnarray}
\pi_{X^{\widetilde k}}(a)=\sum_{j=1}^{\widehat N}\sup_{\vert \beta\vert\le \widehat N}
\sup_{\vert(\widetilde \xi,s)\vert\ge 1}\left\Vert \frac{\partial^{\beta_0}}
{\partial\xi_0^{\beta_0}}\cdots
\frac{\partial^{\beta_{n-1}}}{\partial\xi_{n-1}^{\beta_{n-1}}}
\frac{\partial^{\beta_{n}}}{\partial s^{\beta_{n}}}
a_j(\cdot,\widetilde\xi,s)\right\Vert_{X^{\widetilde k}(\overline{\mathcal
O})}/(1+\vert(\widetilde\xi,s)\vert)^{\kappa-j-\vert\beta\vert} \nonumber \\
+\sup_{\vert (\widetilde\xi,s)\vert\le 1}\Vert
a(\cdot,\widetilde\xi,s)\Vert_{X^{\widetilde k}(\overline{\mathcal O})}. \nonumber
\end{eqnarray}
Obviously for any  $k\in \{0,1\}$
\begin{equation}
\pi_{W^{k,\infty}(B(0,\delta(y^*)))}(\chi_\nu)\le C_{9}\widetilde\varphi^\frac 23(y^*).
\end{equation}
Obviously the  pseudodifferential operators with the symbols
$\Gamma_\beta^\pm$ belongs \\
to the class $W^{k,\infty}_{cl}S^{1,s}(B(0,\delta(y^*)))$  for any $k\in \{0,1\}$
and
\begin{equation}
\pi_{W^{k,\infty}(B(0,\delta(y^*)))}(\Gamma_\beta^\pm)\le C_{10}\widetilde\varphi^\frac k3(y^*).
\end{equation}

Denote $a_\beta(y,\xi)=\beta(y)\sum_{k=1}^{n-1}\xi_k^2-2\xi_n\sum_{k=1}^{n-1} \partial_{y_k}\vartheta \xi_k +\xi_n^2(1+G),$ $a_\beta(y,\xi,\eta)=\beta(y)\sum_{k=1}^{n-1}\xi_k\eta_k-(\xi_n\sum_{k=1}^{n-1} \partial_{y_k}\vartheta \eta_k+\eta_n\sum_{k=1}^{n-1} \partial_{y_k}\vartheta \xi_k) +\xi_n\eta_n(1+G).$ We set $A(y,D)=\beta\Delta_{\vartheta}.$ Then $P_{\rho,\beta}(y,D)=\rho\partial_{y_0}-A(y,D).$
We have

\begin{proposition}\label{opana} Let $w\in H^{1,2}(\mathcal Q),$ $\mbox{supp}\, w\subset
B(y^*,\delta)\cap \mbox{supp}\,\eta_\ell$  and $P_\beta(y,D,\widetilde s)
\chi_\nu w\in L^2(\mathcal Q).$ Then
there exist positive constants $\delta(y^*), C_{11}, C_{12}$
independent of $s$
 such that for all
and $s\ge s_0$ we have
\begin{eqnarray}\label{klop}
C_{11}\int_{\mathcal Q}(\vert s\vert\widetilde\varphi\sum_{k=1}^n\vert \partial_{y_k}\chi_\nu w\vert^2
+ \vert s\vert^3\widetilde\varphi^3\vert \chi_\nu w\vert^2)dy
+ \varXi_\beta(\chi_\nu w)\\
\le  \Vert P_\beta(y,D,\widetilde s)\chi_\nu w\Vert^2_{L^2(\mathcal Q)}
+ C_{12} \epsilon(
\delta)
\Vert(\partial_{y_n}\chi_\nu w,\chi_\nu w)(\cdot,0) \Vert^2_{H^{0,\frac 12,\widetilde s}(\Bbb R^n)\times H^{\frac 34,\frac 32,\widetilde s}(\Bbb R^n)},\nonumber
\end{eqnarray}
where $\epsilon(\delta)\rightarrow +0$ as $\delta\rightarrow
+0$ and
$$
\varXi_\beta(w)=\sum_{j=1}^3 \frak I_j(\beta,w),\quad \frak I_1(\beta,w)
=\int_{\Bbb R^n}
(\vert\widetilde  s\vert\beta^2(y^*)\varphi_{n}(y^*)\vert\partial_{y_n}
w\vert^2 +\vert\widetilde
s\vert^3\beta^2(y^*)\varphi^3_{n}(y^*)\vert w\vert^2)\vert_{y_n=0}d\widetilde y,
$$
\begin{equation}\label{01}\frak I_2(\beta,w)=- \frac 12 Re
\int_{\Bbb R^n} 2\vert \widetilde s\vert\beta(y^*)\partial_{y_n} w
\overline{(\rho\partial_{y_0}w+(\nabla_{\widetilde \xi}a_\beta(y^*,\widetilde \nabla w,0),
 \vec \varphi(y^*))}\vert_{y_n=0} d\widetilde y,\end{equation}
\begin{equation}\label{02} \frak I_3(\beta,w)=\int_{\Bbb R^n}\vert
\widetilde s\vert\beta(y^*)\varphi_{n}(y^*)(a_\beta(y^*,\widetilde \nabla w,0)
-\widetilde s^2a_\beta(y^*,\varphi_1(y^*),\dots,\varphi_{n-1}(y^*),0)\vert
 w\vert^2)\vert_{y_n=0} d\widetilde y.
\end{equation}
\end{proposition}

{\bf Proof.} It suffices to prove the statement of the lemma separately for
$\mbox{Re}\, w_\nu, \mbox{Im}\, w_\nu.$  Let  $v_\nu=\mbox{Re}\, w_\nu$ or
$v_\nu=\mbox{Im}\, w_\nu.$ We write $A(y,D)=\sum_{j,k=1}^na_{jk}(y)\partial^2_{y_ky_j} .$  For simplicity instead of the notation $a_\beta (y,\xi)$
we use the notation
$a(y,\xi)=\sum_{k,j=1}^na_{kj}(y)\xi_k\xi_j,$ where $a_{kj}=a_{jk}$ for all
$k,j\in\{1,\dots, n\}$.
We set
$$
a_j(y,\xi)=\partial_{y_j}a(y,\xi),\quad
a^{(j)}(y,\xi)=\partial_{\xi_j}a(y,\xi),
$$
$$
a^{(j,m)}(y,\xi)=\partial^2_{\xi_j \xi_m}a(y,\xi), \quad a(y,\eta,\xi)=\sum_{k,j=1}^n
a_{kj}(y)
\eta_k\xi_j.
$$

We introduce the operators
\begin{equation}\label{poker}
L_1(y,D,\widetilde s) v_{\nu}=\sum_{k=1}^n
s\partial_{y_k}\varphi a^{(k)}(y,\nabla v_{\nu})+\rho\partial_{y_0} v_{\nu}, \quad L_2(y,D,\widetilde s)v_{\nu}
= -A(y,D)v_{\nu}
-s^2a(y,\nabla \varphi) v_{\nu}.\nonumber
\end{equation}
We set $P(y,D,\widetilde s)=L_1(y,D,\widetilde s)+L_2(y,D,\widetilde s).$
Then
\begin{equation}\label{kp1}
P_\beta(y,D,\widetilde s)v_\nu= P(y,D,\widetilde s)v_\nu+ f_\nu\quad \mbox{in}\quad\mathcal Q.
\end{equation}
Function  $f_\nu$ satisfies the estimate
\begin{equation}\label{kp2}
\Vert f_\nu\Vert_{L^2(\mathcal Q)}\le C_{13}\Vert v_\nu\Vert_{H^{0,1,\widetilde s}(\mathcal Q)}.
\end{equation}
Taking the $L^2$- norm of (\ref{kp1}) we have
\begin{eqnarray}\label{okop}\Vert P_\beta(y,D,\widetilde s)v_{\nu}-f_\nu\Vert^2
_{L^2(\mathcal Q)}
=\Vert L_1(y,D,\widetilde s)v_{\nu}\Vert^2_{L^2(\mathcal Q)}\nonumber\\
+ \Vert L_2(y,D,\widetilde s)v_
\nu\Vert^2_{L^2(\mathcal Q)}+2\mbox{Re}\,(L_1(y,D,\widetilde s)v_{\nu},
L_2(y,D,\widetilde s)v_{\nu})_{L^2(\mathcal Q)}.
\end{eqnarray}
The following equality is proved in \cite{Im}:
\begin{eqnarray}\label{LK}
\mbox{Re}\,(L_1(y,D,\widetilde s)v_{\nu},L_2(y,D,\widetilde s)v_{\nu})
_{L^2(\mathcal Q)}
                                       \nonumber\\
= -\mbox{Re}\,\int_{\partial\mathcal Q}a(y,\vec e_n,\nabla v_{\nu})
{L_1(y,D,\widetilde s)v_{\nu}} d\Sigma-\vert s\vert\int_{\partial\mathcal Q}
a(y,\vec e_n,\nabla
\varphi)a(y,\nabla v_{\nu},\nabla v_{\nu}) d\Sigma\nonumber\\
+\vert s\vert^3\int
_{\partial\mathcal Q}a(y,\nabla\varphi,\nabla \varphi)a(y,\vec e_n,\nabla \varphi)
\vert v_{\nu}\vert^2 d\Sigma
+ \vert s\vert\int_{\mathcal Q}\mathcal G(y,\widetilde s,v_{\nu})dy\nonumber\\
+ \int_{\mathcal Q}\frac{\vert s\vert}{2}\left (\sum_{k=1}^n a^{(k)}_k(y,\nabla v_{\nu})
L_1(y,D,\widetilde s)v_{\nu}-\widetilde\theta(y)(a(y,\nabla v_{\nu},\nabla v_{\nu})
-s^2a(y,\nabla\varphi,\nabla\varphi)\vert v_{\nu}\vert^2)\right )dy\nonumber\\
+ \int_{\mathcal Q}(\sum_{k,m=1}^n\partial_{y_0}(\rho a_{km})\partial_{y_k}v_\nu\partial_{y_m}v_\nu-s^2\partial_{y_0}(\rho a(y,\nabla \varphi)) \vert v_\nu\vert^2)dy,
\end{eqnarray}
where
$$
\mathcal G(y,\widetilde s,w)=\{a,\{a,\varphi\}\}(y,\nabla w)
+ s^2\sum_{k,j=1}^na_{j}
(y,\nabla \varphi)a^{(k)}
(y,\nabla \varphi) w^2
$$
$$+s^2\sum_{k,j=1}^n\partial^2_{y_ky_j}\varphi a^{(k)}
(y,\nabla\varphi)a^{(j)}(y,\nabla\varphi) w^2
$$
and
$$
\widetilde\theta(y)=\sum_{k,m=1}^n(\partial^2_{y_k y_m}\varphi a^{(k,m)}(y,\nabla \varphi)
+\partial_{y_k}\varphi a_m^{(k,m)}(y,\nabla v_{\nu})).
$$
Observe that the function $\widetilde\theta(y)$ is independent of $w$ and by  (\ref{lomka1})
\begin{equation}\label{ok}
\sup_{y\in \mbox{supp}\, \eta_\ell}\vert \widetilde\theta(y)\vert
\le C_{14}\varphi(y^*).
\end{equation}
If $y\in  B(y^*,\delta)\cap \mbox{supp}\, \eta_\ell$ then either
$\frac 12\le 2^{-\ell}2^\frac{1}{\theta(y_0)^\frac 14}\le 2$ or $\frac 12\le 2^{-\ell-1}
2^\frac{1}{\theta(y_0)^\frac 14}\le 2$ or $\frac 12\le 2^{-\ell+1}2^\frac{1}{\theta(y_0)^\frac 14}\le 2.$
This is equivalent
\begin{equation}\nonumber
2^{\ell-2}\le 2^\frac{1}{\theta(y_0)^\frac 14} \le 2^{\ell+2}\quad \forall y\in  B(y^*,\delta)\cap \mbox{supp}\, \eta_\ell.
\end{equation}
So, we have
\begin{equation}\nonumber
(\ell+2)^{-4}\le \theta(y_0) \le (\ell-2)^{-4}\quad \forall y\in  B(y^*,\delta)\cap \mbox{supp}\, \eta_\ell.
\end{equation}
Hence  by definition (\ref{Aziza}) of the function $\theta$  there exist $\ell_0$ such that for all $\ell\ge \ell_0$ and for all $y$ from $ B(y^*,\delta)\cap \mbox{supp}\, \eta_\ell$
\begin{equation}\label{capitol} (\ell+2)^{-4}\le T+y_0\le (\ell-2)^{-4}
\quad \mbox{or}\,\,(\ell+2)^{-4}\le T-y_0\le (\ell-2)^{-4}.
\end{equation}
Then for any $y_1=(y_{1,0},\dots, y_{1,n}),y_2=(y_{2,0},\dots, y_{2,n})\in  B(y^*,\delta)\cap \mbox{supp}\, \eta_\ell$ we have
\begin{equation}\label{znarok2}
\vert y_{1,0}-y_{2,0}\vert =\vert y_{1,0}+T-T-y_{2,0}\vert\le  (\ell-2)^{-4}-(\ell+2)^{-4}\le C_{15}\ell^{- 5}.
\end{equation}

Let  $\widetilde y^*=(\widetilde y_0^*, \dots ,\widetilde y_n^*)$ be some point from $ B(y^*,\delta)\cap \mbox{supp}\, \eta_\ell.$
We claim that  for any positive $\epsilon$ there exists positive  $\delta(\epsilon)$  independent of $\ell$ such that
\begin{equation}\label{kokom}
\vert\partial^{\alpha}\varphi(y)-\partial^{\alpha}\varphi(\widetilde y^*)\vert \le \epsilon \vert \varphi(\widetilde y^*)\vert\quad \forall y
\in  B(y^*,\delta)\cap \mbox{supp}\, \eta_\ell \quad \forall \vert\alpha\vert\le 2,\,\,\alpha_0=0.
\end{equation}
%Indeed, since $\varphi(y)<C_{20}<0$ on $Q$ for any $\delta_1>0$ we have (\ref{kokom}) on $[-T-\delta_1, T-\delta_1]\times \Omega.$
%Hence it suffices to consider the case $y_0^*\in\{\pm T\}.$ Let $y_0^*=T.$ Proof of the case $y_0^*=-T$ is the same.

Estimate (\ref{kokom}) follows  from the following  inequality:
% . Let $y \in  B(y^*,\delta)\cap \mbox{supp}\, \eta_\ell$ then
\begin{eqnarray}\label{nonsense}
\vert\partial^{\alpha}\varphi(y)-\partial^{\alpha}\varphi(\widetilde y^*)\vert\le \vert\partial^{\alpha}\varphi(\widetilde y^*_0,y')-\partial^{\alpha}\varphi(y)\vert+\vert\partial^{\alpha}\varphi(\widetilde y^*_0,y')-\partial^{\alpha}\varphi(y^*)\vert\nonumber\\\le\sup_{\widetilde y\in  B(y^*,\delta)\cap \mbox{supp}\, \eta_\ell}\vert \nabla' \partial^{\alpha}\varphi (\widetilde y)\vert \vert y-\widetilde y^*\vert +\sup_{\widetilde y\in  B(y^*,\delta)\cap \mbox{supp}\, \eta_\ell}\vert \partial_{y_0}\partial^{\alpha}\varphi (\widetilde y)\vert \vert y_0-\widetilde y_0^*\vert\nonumber\\\le C_{16} (
 \sup_{\widetilde y\in  B(y^*,\delta)\cap \mbox{supp}\, \eta_\ell}\vert  \varphi (\widetilde y)\vert \vert y-\widetilde y^*\vert+\sup_{\widetilde y\in  B(y^*,\delta)\cap \mbox{supp}\, \eta_\ell}\vert  \varphi^\frac 43 (\widetilde y)\vert \vert y_0-\widetilde y_0^*\vert)\nonumber \\
 \le C_{17} (
 \sup_{\widetilde y\in  B(y^*,\delta)\cap \mbox{supp}\, \eta_\ell}\vert  \varphi (\widetilde y)\vert \vert y-\widetilde y^*\vert+\sup_{\widetilde y\in  B(y^*,\delta)\cap \mbox{supp}\, \eta_\ell}\vert  \varphi (\widetilde y) \ell^4\vert \vert y_0-\widetilde y_0^*\vert)\nonumber \\
  \le C_{18} (
 \sup_{\widetilde y\in  B(y^*,\delta)\cap \mbox{supp}\, \eta_\ell}\vert  \varphi (\widetilde y)\vert \vert y-\widetilde y^*\vert+\sup_{\widetilde y\in  B(y^*,\delta)\cap \mbox{supp}\, \eta_\ell}\vert  \varphi (\widetilde y)\vert \ell^{-1})\nonumber \\\le \epsilon \vert  \varphi(\widetilde y^*)\vert .
\end{eqnarray}
In order to get the last two inequalities in (\ref{nonsense}) we used   (\ref{znarok2}), (\ref{lomka1}) and (\ref{lomka}).

We introduce the form $\frak G(\widetilde s,\nabla v_{\nu})$ in the following
way:
In the function $\mathcal G(\widetilde y^*,\widetilde s,\nabla v_{\nu})$ we replace function
$\partial_{y_n}
v_{\nu}$ by $ m_*\left (
\sum_{j,k=1}^{n-1}\partial_{y_j}\varphi(\widetilde y^*)a_{jk}(\widetilde y^*)\partial_{y_k}v_{\nu}+\rho(\widetilde y^*)\partial_{y_0} v_{\nu}/s\right)$\newline where $m_*=-\frac{1}{\sum_{j=1}^na_{jn}(\widetilde y^*)
\partial_{y_j}\varphi(\widetilde y^*)}$ and  set $$\widetilde P(D)v_\nu=m_*\left (
\sum_{j,k=1}^{n-1}\partial_{y_j}\varphi(y^*)a_{jk}(\widetilde y^*)\partial_{y_k}v_{\nu}+\rho(\widetilde y^*)\partial_{y_0} v_{\nu}/s\right).$$  Then
\begin{equation}\label{molot}
\partial_{y_n}
v_{\nu}=\widetilde P(D)v_\nu+\frac{m_*}{s}L_1(\widetilde y^*,s,D)v_\nu.\end{equation}
Since $a_{kj}(\widetilde y^*)=0$ for $k\ne j$ by (\ref{kokom}) and  (\ref{pinok220}) there exist a positive constant $C_{19}$ such that
\begin{equation}\label{prizrak}
\vert \sum_{j=1}^na_{jn}(y)\partial_{y_j}\varphi(y)\vert\ge C_{19}\widetilde \varphi(y)
>0\quad \forall y\in B(y^*,\delta)\cap \mbox{supp}\, \eta_\ell.
\end{equation}
By (\ref{kokom}) for any positive $\epsilon$ one can take a positive $\delta(\epsilon)$ such that
$$
\left\vert \int_{\mathcal Q} \mathcal G(y,\widetilde s,\nabla v_{\nu})
- \frak G(\widetilde s,\nabla v_{\nu}) dy\right\vert\le \left\vert \int_{\mathcal Q} \mathcal G(y,\widetilde s,\nabla v_{\nu})
-  \mathcal G(\widetilde y^*,\widetilde s,\nabla v_{\nu}) dy\right\vert
$$
$$+\left\vert \int_{\mathcal Q}\frak G (\widetilde y^*,\widetilde s,\nabla v_{\nu})
- \mathcal G(\widetilde y^*,\widetilde s,\partial_{y_1} v_{\nu},\dots , \partial_{y_{n-1}} v_\nu, P(D)v_\nu +\frac{m_*}{s}L_1(\widetilde y^*,s,D)v_\nu) dy\right\vert
$$
$$
\le \frac{C_{20}}{\widetilde s^2}\Vert L_1(\widetilde y^*,D,\widetilde s)v_{\nu}\Vert^2
_{L^2(\mathcal Q)}
+\epsilon\Vert v_{\nu}\Vert^2_{H^{0,1,\widetilde s}(\mathcal Q)}
$$
$$
\le \frac{C_{21}}{\widetilde s^2}\Vert L_1(y,D,\widetilde s)v_{\nu}\Vert^2
_{L^2(\mathcal Q)}
+\epsilon\Vert v_{\nu}\Vert^2_{H^{0,1,\widetilde s}(\mathcal Q)}.
$$

Let
\begin{equation}\label{propoganda}
\xi = \left(\widetilde\xi,m_*(\sum_{j,k=1}^{n-1}\partial_{y_j}\varphi(\widetilde y^*)a_{jk}(\widetilde y^*)\xi_k+\rho(\widetilde y^*)\xi_0/s)
\right),
\end{equation}
\begin{equation}\label{propoganda1}
\quad \zeta =(\xi_0,\xi_1+{i}\vert s\vert\partial_{y_1}\varphi(\widetilde y^*),\dots, \xi_n+{i} \vert s\vert\partial_{y_n}\varphi(\widetilde y^*)) .
\end{equation}
We set
$$
q(\widetilde \xi,s)=\sum_{k,j=1}^n\partial^2_{y_ky_j} \varphi(\widetilde y^*)a^{(k)}(\widetilde y^*,\zeta)
\overline{a^{(j)}(\widetilde y^*,\zeta)}+\frac {1}{\vert s\vert} \mbox{Im}\sum_{k=1}^na_k(\widetilde y^*,\zeta)
\overline{a^{(k)}(\widetilde y^*,\zeta)},
$$ where in the right hand side above formula  $\xi$ and $\zeta$ are given  (\ref{propoganda}), (\ref{propoganda1}).
Observe that
\begin{equation}
\int_{\Bbb R^{n+1}_+}q(\widetilde \xi,s)\vert \mbox{F}_{\widetilde y
\rightarrow \widetilde \xi}v_{\nu}\vert^2d\widetilde \xi dy_n=\int_{\Bbb R^{n+1}_+}
\frak G(\widetilde s,\nabla v_{\nu}) dy=\int_{\mathcal Q}
\frak G(\widetilde s,\nabla v_{\nu}) dy.
\end{equation}
Denote $\widetilde w(\widetilde \xi,s, y_n)=(2\pi)^{-\frac{n}{2}}\int_{\Bbb R^n}sign(\mbox{Re}\,
a(\widetilde y^*,\xi+{i}\vert s\vert\nabla \varphi(y^*))\mbox {F}_{\widetilde y\rightarrow
\widetilde \xi}
v_{\nu} e^{{i}<\widetilde \xi,\widetilde y>}d\widetilde \xi$,
where $\mbox {F}_{\widetilde y\rightarrow \widetilde \xi}$
is the Fourier transform given by
$$
F_{\widetilde y\rightarrow
\widetilde \xi}u=\frac{1}{(2\pi)^\frac{n}{2}}\int_{{\Bbb R}^{n}}
e^{-{i}\sum_{j=0}^{n-1}y_j\xi_j}u(y_0,\dots,y_{n-1})d\widetilde y.
$$

Taking the scalar product of the function $L_2(y,D,\widetilde s)v_{\nu}$ and
$\overline{\widetilde w}$  in ${L^2(\mathcal Q)}$ we have
\begin{eqnarray}\label{begemot1!}
\int_{\mathcal Q}(a(y,\nabla v_{\nu},\overline{\nabla\widetilde  w})
- s^2a(y,\nabla \varphi,\nabla \varphi)v_{\nu}\overline{\widetilde w})dy
=\int_{\Bbb R^n}\partial_{\vec\nu_a} v_{\nu}(\widetilde y,0)
\overline{\widetilde w(\widetilde y,0)}d\widetilde y\nonumber\\
-\sum_{k,j=1}^n\int_{\mathcal Q}\partial_{y_k} a_{kj}
\partial_{y_j} v_{\nu}\overline{\widetilde w}dy
+ (L_2(y,D,\widetilde s)v_{\nu},\overline{\widetilde w})_{L^2(\mathcal Q)},
\end{eqnarray}
where $\partial_{\vec\nu_a} w
=\sum_{j=1}^na_{nj}\partial_{y_j}w.$

By (\ref{kokom}) for any  $\epsilon \in (0,1)$ there exists $\delta(\epsilon)$
such that
\begin{eqnarray}\label{begemot2}
\int_{\mathcal Q}\vert\widetilde s\vert s^2\vert a(\widetilde y^*,\nabla\varphi (\widetilde y^*),\nabla\varphi
(y^*))-a(y,\nabla\varphi ,\nabla\varphi)\vert \vert v_{\nu}\vert^2dy
                                           \nonumber\\
+\int_{\mathcal Q}\vert\widetilde s\vert\vert a(y,\nabla v_{\nu},\overline{\nabla \widetilde w})
- a(\widetilde y^*,\nabla v_{\nu},\overline{\nabla \widetilde w})\vert dy
\le \epsilon \vert\widetilde s\vert\Vert v_{\nu}\Vert^2_{H^{0,1,\widetilde s}(\mathcal Q)}.
\end{eqnarray}
The inequalities  (\ref{begemot1!})  and (\ref{begemot2}) imply that for any
positive $\epsilon$ there exists  a constant $C_{22}$ such that
\begin{eqnarray}
\vert\widetilde s\vert\left\vert \int_{\mathcal Q}(-a(\widetilde y^*,\nabla v_{\nu},
\overline{\nabla\widetilde  w})+s^2a(\widetilde y^*,\nabla \varphi(\widetilde y^*),
\nabla \varphi(\widetilde y^*))v_{\nu}\overline{\widetilde w}dy\right\vert\\
\le C_{23}\left\vert\widetilde s \int_{\Bbb R^n}\partial_{\vec\nu_a} v_{\nu}(\widetilde y,0)\overline{\widetilde w}(\widetilde y,0)d
\widetilde y\right\vert
+ \epsilon\Vert L_2(y,D,\widetilde s)v_{\nu}\Vert^2_{L^2(\mathcal Q)}
+ C_{22}\epsilon \vert \widetilde s\vert\Vert v_{\nu}\Vert^2_{H^{0,1,\widetilde s}(\mathcal Q)}.
                                    \nonumber
\end{eqnarray}
We set $\widehat \nabla w=(\partial_{y_0}w, \dots, \partial_{y_{n-1}}w, m_*\left (\sum_{j,k=1}^{n-1}
\partial_{y_j}\varphi(\widetilde y^*)a_{jk}(\widetilde y^*)\partial_{y_k}w+\rho(\widetilde y^*)\partial_{y_0}w/s\right ))$ and \newline$\widehat
\nabla \widetilde w=(\partial_{y_0}\widetilde w, \dots, \partial_{y_{n-1}}
\widetilde w,m_*\left (
\sum_{j,k=1}^{n-1}\partial_{y_j}\varphi(\widetilde y^*)a_{jk}(\widetilde y^*)\partial_{y_k}
\widetilde w+\rho(\widetilde y^*)\partial_{y_0}\widetilde w/s\right )).$
Hence, if $\xi$ is given by (\ref{propoganda}), then we have
\begin{eqnarray}\label{LL0}
\int_{\Bbb R^n_+}\vert\widetilde s\vert \vert i\rho(\widetilde y^*)\xi_0+ a(\widetilde  y^*,\xi+{i} \vert s\vert\nabla \varphi(\widetilde y^*))
\vert \vert\mbox{F}_{\widetilde y\rightarrow \widetilde \xi} v_{\nu}\vert^2
d\widetilde \xi dy_n\\
= \left\vert \widetilde s\int_{\mathcal Q}(-a(\widetilde y^*,\widehat\nabla v_{\nu},
\overline{\widehat \nabla\widetilde  w})+s^2a(\widetilde y^*,\nabla \varphi(y^*),\nabla
\varphi(\widetilde y^*))v_{\nu}\overline{\widetilde w}dy\right\vert
                                                  \nonumber\\
\le C_{24}\vert\widetilde s\vert \int_{\Bbb R^n}\left\vert\partial_{\vec\nu_a} v_{\nu}(\widetilde y,0)\overline{\widetilde w(\widetilde y,0)}
\right\vert
d\widetilde y+\epsilon\Vert L_2(y,D,\widetilde s)v_{\nu}\Vert^2
_{L^2(\mathcal Q)}
+C_{25}\epsilon\vert\widetilde s\vert \Vert v_{\nu}\Vert^2_{H^{0,1,\widetilde s}(\mathcal Q)}.
                                                      \nonumber
\end{eqnarray}

Observe that pseudoconvexity Condition \ref{A1} implies that there exists a
positive constant $C_{26}$ such that
\begin{equation}\label{LL1}
\vert  s\vert q(\widetilde \xi,s)+\frac{1}{\vert\widetilde s\vert}\vert a(\widetilde y^*, \xi)-a(y^*,  s\nabla
\varphi(\widetilde y^*))\vert\ge C_{26}\vert\widetilde s\vert M^2
(\widetilde\xi,\widetilde s)
\quad \forall (\xi,\widetilde s)\in \Bbb M.
\end{equation}
Therefore, from (\ref{LL1}) and (\ref{LL0}),
for some positive  constant $C_{27}$ we have the inequality
\begin{eqnarray}
C_{27}\vert \widetilde s\vert\int_{\mathcal Q} (\sum_{j=1}^{n-1}\vert \partial_{y_j} v_{\nu}
\vert^2
+\widetilde s^2\vert v_{\nu}\vert^2)dy\le C_{28}(\vert\widetilde s\vert \int
_{\Bbb R^n}
\vert \partial_{\vec\nu_a} v_{\nu}(\widetilde y,0)
\overline{\widetilde w(\widetilde y,0)}\vert d\widetilde y\nonumber\\
+\epsilon\Vert L_2(y,D,\widetilde s)v_{\nu}\Vert^2_{L^2(\mathcal Q)})
+ \int_{\mathcal Q}
\vert s\vert\frak G(\widetilde s,\nabla v_{\nu}) dy+C_{29}\epsilon\vert\widetilde s\vert\Vert v_{\nu}\Vert^2
_{H^{0,1,\widetilde s}(\mathcal Q)}.
\end{eqnarray}

Thanks to (\ref{LK})
\begin{eqnarray}\label{mina1}
C_{30}\vert\widetilde s\vert\int_{\mathcal Q}
(\sum_{j=1}^{n-1}\vert \partial_{y_j} v_{\nu}
\vert^2
+\widetilde s^2\vert v_{\nu}\vert^2)dy
\le C_{31}\left(\vert\widetilde s\vert \int_{\Bbb R^n}
\vert\partial_{\vec\nu_a} v_{\nu}(\widetilde y,0)\overline{
\widetilde w(\widetilde y,0)}\vert d\widetilde y
+ \epsilon\Vert L_2(y,D,\widetilde s)w_
\nu\Vert^2_{L^2(\mathcal Q)}\right)              \nonumber\\
+2(L_1(y,D,\widetilde s)v_{\nu}, L_2(y,D,\widetilde s)v_{\nu})
_{L^2(\mathcal Q)}+\mbox{Re}\,
\int_{\Bbb R^n}a(y,\vec e_n,\nabla v_{\nu}){L_1(y,D,\widetilde s)
v_{\nu} }\vert_{y_n=0}
d\Sigma\nonumber\\-\vert s\vert\int_{\Bbb R^n}a(y,\vec e_n,\nabla \varphi)a(y,\nabla
v_{\nu},{\nabla v_{\nu}})\vert_{y_n=0} d\Sigma+\vert s\vert^3\int_{\Bbb R^n}
a(y,\nabla\varphi,\nabla \varphi)a(y,\vec e_n,\nabla \varphi)
\vert v_{\nu}\vert^2\vert_{y_n=0} d
\Sigma                                                 \nonumber\\
-\int_{\mathcal Q}\frac{\vert s\vert}{2}(\sum_{k,m=1}^n a^{(k)}_k(y,\nabla v_{\nu})
L_1(y,D,s)v_\nu-\widetilde\theta(y)(a(y,\nabla v_{\nu})
-s^2a(y,\nabla\varphi)
\vert v_{\nu}\vert^2))dy\nonumber\\
+ \int_{\mathcal Q}(\sum_{k,m=1}^n\partial_{y_0}(\rho a_{km})\partial_{y_k}v_\nu\partial_{y_m}v_\nu-s^2\partial_{y_0}(\rho a(y,\nabla \varphi)) \vert v_\nu\vert^2)dy\nonumber\\
+\frac{C_{32}}{s^2}\Vert L_1(y,D,\widetilde s)v_{\nu}\Vert^2
_{L^2(\mathcal Q)}
+C_{33}\epsilon\vert\widetilde s\vert\Vert v_{\nu}\Vert^2_{H^{0,1,\widetilde s}(\mathcal Q)}.
\end{eqnarray}

 Now we estimate the derivative of the function $v_\nu$ respect to variable $y_n.$
Taking the scalar product of the function $L_2(y,D,\widetilde s)v_{\nu}$ and
$v_\nu$ we have
\begin{eqnarray}\label{begemot1l}
\int_{\mathcal Q}(a(y,\nabla v_{\nu},\nabla v_\nu)
-s^2a(y,\nabla \varphi,\nabla \varphi)v^2_{\nu})dy
=\int_{\Bbb R^n}\partial_{\vec\nu_a} v_{\nu}(\widetilde y,0)
v_\nu (\widetilde y,0)d\widetilde y\nonumber\\
-\int_{\mathcal Q}\sum_{k,j=1}^n\partial_{y_k} a_{kj}
\partial_{y_j} v_{\nu} v_\nu dy
+ (L_2(y,D,\widetilde s)v_{\nu},v_\nu)_{L^2(\mathcal Q)}.
\end{eqnarray}

Thanks to (\ref{begemot1l}) for any positive $\epsilon$  there exists a constant $C_{34}$ independent of $\widetilde s$
such that
\begin{eqnarray}\label{mina}
\vert\widetilde s\vert\int_{\mathcal Q} \vert \partial_{y_n} v_{\nu}\vert^2 dy\le C_{34}(\int_{\mathcal Q} (\sum_{j=1}^{n-1}\vert \partial_{y_j} v_{\nu}
\vert^2
+\widetilde s^2
\vert v_{\nu}\vert^2)dy\nonumber\\
+\vert\widetilde s\int_{\Bbb R^n}\partial_{\vec\nu_a} v_{\nu}(\widetilde y,0)
v_\nu (\widetilde y,0)d\widetilde y\vert)+\epsilon\Vert L_2(y,D,\widetilde s)v_{\nu}\Vert^2
_{L^2(\mathcal Q)}.
\end{eqnarray}

Using (\ref{okop}), from (\ref{mina}) and (\ref{mina1}) we obtain
\begin{eqnarray}\label{logotip8}
C_{35}\vert\widetilde s\vert\int_{\mathcal Q} (\sum_{j=1}^{n}\vert \partial_{y_j} v_{\nu}
\vert^2
+\widetilde s^2
\vert v_{\nu}\vert^2)dy
\le C_{36}\vert\widetilde s\vert\left\vert \int_{\Bbb R^n}
\partial_{\vec\nu_a} v_{\nu}(\widetilde y,0)
\overline{\widetilde w(\widetilde y,0)}d\widetilde y\right\vert
                                         \\
+\nonumber
\vert\vert\widetilde s\vert \int_{\Bbb R^n}\partial_{\vec\nu_a} v_{\nu}(\widetilde y,0)
v_\nu (\widetilde y,0)d\widetilde y\vert
+\epsilon\Vert L_2(y,D,\widetilde s)v_{\nu}\Vert^2_{L^2(\mathcal Q)}+
\int_{\Bbb R^n}a(y,\vec e_n,\nabla v_{\nu}){L_1(y,D,\widetilde s)
v_{\nu}}\vert_{y_n=0}
d\Sigma\\-\vert s\vert\int_{\Bbb R^n}a(y,\vec e_n,\nabla \varphi)
a(y,\nabla v_{\nu},{\nabla v_{\nu}}) \vert_{y_n=0}d\Sigma+s^3\int_{\partial\mathcal Q}
a(y,\nabla\varphi,\nabla \varphi)a(y,\vec e_n,\nabla \varphi)\vert v_{\nu}
\vert^2
d\Sigma                       \nonumber\\
-\int_{\mathcal Q}\frac{\vert s\vert}{2}\left(\sum_{k,m=1}^n a^{(k)}_
k(y,\nabla v_{\nu})\partial_{y_m}\varphi
a^{(m)}(y,\nabla v_{\nu})-\widetilde\theta(y)(a(y,\nabla v_{\nu},\nabla  v_{\nu})
-s^2a(y,\nabla\varphi,\nabla\varphi)\vert v_{\nu}\vert^2)\right)dy
                                            \nonumber\\
                                             +\int_{\mathcal Q}(\sum_{k,m=1}^n\partial_{y_0}(\rho a_{km})\partial_{y_k}v_\nu\partial_{y_m}v_\nu-s^2\partial_{y_0}(\rho a(y,\nabla \varphi)) \vert v_\nu\vert^2)dy\nonumber\\
+\Vert P(y,D,\widetilde s)v_{\nu}\Vert^2_{L^2(\mathcal Q)}-\frac 12\Vert
L_1(y,D,\widetilde s)
v_{\nu}\Vert^2_{L^2(\mathcal Q)}-\Vert L_2(y,D,\widetilde s)v_{\nu}\Vert^2
_{L^2(\mathcal Q)}.\nonumber
\end{eqnarray}
Now we estimate some integrals on the right hand side of (\ref{logotip8}):
\begin{equation}\label{puk1}
\frac{1}{2}\vert \int_{\mathcal Q}\sum_{k=1}^n a^{(k)}_k(y,\nabla v_{\nu})
L_1(y,D,\widetilde s)v_{\nu}dy\vert\le C_{37}\int_{\mathcal Q}\sum_{j=1}^{n}\vert \partial_{y_j} v_{\nu}
\vert^2dy+\epsilon \Vert L_1(y,D,\widetilde s)v_{\nu}\Vert^2
_{L^2(\mathcal Q)}
\end{equation}
and
\begin{equation}\label{puk2}
\vert \int_{\mathcal Q}(\sum_{k,m=1}^n\partial_{y_0}(\rho a_{km})\partial_{y_k}v_\nu\partial_{y_m}v_\nu-s^2\partial_{y_0}(\rho a(y,\nabla \varphi)) \vert v_\nu\vert^2)dy\vert\le C_{38}\widetilde \varphi^\frac 13(\widetilde y^*)\Vert v_\nu\Vert^2_{H^{0,1,\widetilde s}(\mathcal Q)}.
 \end{equation}
Integrating by parts we have
\begin{eqnarray}
\int_{\mathcal Q}\widetilde\theta(a(y,\nabla v_{\nu},\nabla v_{\nu})-s^2a(
y,\nabla\varphi,\nabla\varphi)\vert v_{\nu}\vert^2)dy=\int _{\mathcal Q}
\widetilde\theta L_2(y,D,\widetilde s)v_{\nu}{v_{\nu}}dy\\ +\int_{\Bbb R^n}
\widetilde\theta\partial_{\vec\nu_a} v_{\nu}(\widetilde y,0){v_{\nu}(
\widetilde y,0)}d\widetilde y+\sum_{j,k=1}^n\int_{\mathcal Q}(\widetilde\theta a_k^{(j)}
(y,\nabla v_{\nu}){v_{\nu}}+a^{(j)}(y,
\nabla v_{\nu})\partial_{y_k}\widetilde\theta{v_{\nu}})dy.\nonumber
\end{eqnarray}
Therefore (\ref{ok}) yields
\begin{eqnarray}\label{puk}
\vert s\vert\left\vert \int_{\mathcal Q}\widetilde\theta(a(y,\nabla v_{\nu},\nabla  v_{\nu})
-s^2a(y,\nabla\varphi)\vert v_{\nu}\vert^2)dy\right\vert
\le \epsilon\Vert L_2(y,D,\widetilde s)
v_{\nu}\Vert^2_{L^2(\mathcal Q)}\nonumber\\
+ C_{39}(\Vert v_{\nu}\Vert^2_{H^{0,1,\widetilde s}(\mathcal Q)}
+ \Vert (\partial_{y_n} v_{\nu}(\cdot,0),v_
\nu(\cdot,0))\Vert^2
_{L^2(\Bbb R^n)\times H^{0,1,\widetilde s}(\Bbb R^n)}).
\end{eqnarray}
Using (\ref{puk}),  (\ref{puk1}) and (\ref{puk2}), from (\ref{logotip8}) we obtain
\begin{eqnarray}\label{logotip1}
C_{40}\vert\widetilde s\vert\int_{\mathcal Q} (\sum_{j=1}^{n}\vert \partial_{y_j} v_{\nu}
\vert^2
+\widetilde s^2\vert v_{\nu}\vert^2)dy
\le C_{41}\vert\widetilde s\vert\Vert (
\partial_{y_n} v_{\nu}(\cdot,0),v_{\nu}(\cdot,0)
)\Vert^2
_{L^2(\Bbb R^n)\times H^{0,1,\widetilde s}(\Bbb R^n)}\nonumber\\+
\int_{\Bbb R^n}a(y,\vec e_n,\nabla v_{\nu}){L_1(y,D,\widetilde s)
v_{\nu}}\vert_{y_n=0}
d\Sigma\\-\vert s\vert\int_{\Bbb R^n}a(y,\vec e_n,\nabla \varphi)
a(y,\nabla v_{\nu},{\nabla v_{\nu}})\vert_{y_n=0} d\Sigma+s^3\int_{\Bbb R^n}
a(y,\nabla\varphi,\nabla \varphi)a(y,\vec e_n,\nabla \varphi)\vert v_{\nu}
\vert^2\vert_{y_n=0} d\Sigma
\nonumber\\
+\Vert P(y,D,\widetilde s)v_{\nu}\Vert^2_{L^2(\mathcal Q)}
-\frac 12\Vert L_1(y,D,\widetilde s)v_{\nu}
\Vert^2_{L^2(\mathcal Q)}-\frac 12\Vert L_2(y,D,\widetilde s)v_{\nu}\Vert^2
_{L^2(\mathcal Q)}.\nonumber
\end{eqnarray}
Next we estimate the difference between the boundary integrals in (\ref{logotip1}) and
$\sum_{j=1}^3 \frak I_j(\beta,v_{\nu})$. Using (\ref{kokom}) we have
\begin{eqnarray}\label{logotip2}
\biggl\vert -\int_{\partial\mathcal Q}a(y,\vec e_n,\nabla v_{\nu})
{L_1(y,D,\widetilde s)v_{\nu}}\vert_{y_n=0} d\Sigma-\vert s\vert\int_{\Bbb R^n}
a(y,\vec e_n,\nabla
\varphi)a(y,\nabla v_{\nu},{\nabla v_{\nu}})\vert_{y_n=0} d\Sigma\nonumber\\
+\vert s\vert^3\int_{\Bbb R^n}a(y,\nabla\varphi,\nabla \varphi)a(y,\vec e_n,\nabla
\varphi)\vert v_{\nu}\vert^2\vert_{y_n=0} d\Sigma-\sum_{j=1}^3 \frak I_j(\beta,v_{\nu})
\biggr\vert                                             \nonumber\\
\le C_{42} \epsilon
\int_{\Bbb R^n}(\vert \widetilde s\vert \sum_{j=1}^{n}\vert \partial_{y_j} v_{\nu}
\vert^2
+\vert \widetilde s\vert^3 \vert v_{\nu}\vert^2 )(\widetilde y,0)
d\widetilde y.
\end{eqnarray}
From (\ref{logotip1}) and (\ref{logotip2}), we obtain (\ref{klop}).
$\blacksquare$

In some cases, we can represent the operator $
P_{\rho,\beta}(y,D,\widetilde s)$ as a product of two first order
pseudodifferential operators.
\begin{proposition}\label{gorokx1}  Let $w\in H^{1,2}(\mathcal Q),$ $\mbox{supp}\, w\subset
B(y^*,\delta)\cap \mbox{supp}\,\eta_\ell, y^*\in \mbox{supp}\,\eta_\ell$  and $P_\beta(y,D,\widetilde s)
\chi_\nu w\in L^2(\mathcal Q).$ assume that  $r_\beta(y^*,\zeta^*)\ne 0$
and $ supp\, {\chi_\nu}\subset \mathcal O(\zeta^*,\delta_1).$
Then we can
factorize the operator $P_{\rho,\beta}(y,D,\widetilde s)$ into the product of two
first order pseudodifferential operators:
\begin{equation} \label{min}
P_{\rho,\beta}(y,D,\widetilde s)w_{\nu}=\beta G(\frac 1i\partial_{y_n}-\Gamma^-_
\beta(y,{\widetilde D},\widetilde s))
(\frac 1i\partial_{y_n}-\Gamma^+_\beta(y,{\widetilde D},\widetilde s))
w_{\nu}+T_\beta w_\nu,
\end{equation}
where $T_\beta : H^{\frac 12 ,1,\widetilde s}({\Bbb R}^{n+1}_+)
\rightarrow L^2(0,\gamma;L^2({\Bbb R}^n))$ operator satisfies the estimate
\begin{equation}\label{lokom}
\Vert T_\beta w_{\nu}\Vert_{L^2(0,\gamma;L^2({\Bbb
R}^n))}\le C_{43}\widetilde\varphi^\frac {5}{12}(\widetilde y^*)\Vert w\Vert_{ H^{\frac 12 ,1,\widetilde s}({\Bbb R}^{n+1}_+)}.
\end{equation}
\end{proposition}
{\bf Proof.} Let
$$
\widetilde R(y,\widetilde D,\widetilde s)=\left [\frac{i\rho(\frac 1 i\partial_{y_0}+{i}
\vert
\widetilde s
\vert\varphi_0)}{\beta G}+\frac{\sum_{j=1}^{n-1}(\frac 1i\partial_{y_j}+{i}\vert
\widetilde s\vert\varphi_j)^2}{G}\right]
$$
and  $\Gamma(y,{\widetilde D},\widetilde s)$ is the operator with symbol
$\Gamma^-_\beta(y,{\widetilde \xi},\widetilde  s)\Gamma^+_\beta
(y,{\widetilde \xi},\widetilde s):$
$$
\Gamma(y,\widetilde \xi,\widetilde s)=(-\vert \widetilde s\vert^2(
\widetilde \mu_\ell\varphi_n)^2+\alpha^+_\beta {i}\vert \widetilde s\vert
\varphi_n+\alpha_\beta^-{i}\vert \widetilde s\vert\varphi_n+\alpha_\beta^+
\alpha_\beta^-)
$$
$$
= \kappa^2(\widetilde \xi,\widetilde s)\left [-\vert \widetilde s\vert^2(
\widetilde \mu_\ell\varphi_n)^2+\frac{i\rho(\xi_0
+{i}
\vert\widetilde s\vert\varphi_0)}{\beta G}+\frac{\sum_{j=1}^{n-1}
(\xi_j+{i}\vert\widetilde s\vert\varphi_j)^2}{G}\right].
$$
We set \begin{equation}\label{zontik}\Upsilon_\ell=B(y^*,2\delta)\cap\text{supp}\,\widetilde \mu_\ell.\end{equation}
Then
\begin{eqnarray}\label{lesopoval1}
\Gamma(y,{\widetilde D},\widetilde s)w_{\nu}=[\Gamma,\eta_\ell]
\chi_\nu(\widetilde D,\widetilde s)w
+\eta_\ell\Gamma(y,\widetilde D,\widetilde s)\chi_\nu(\widetilde D,\widetilde s)w               \\
= [\Gamma,\eta_\ell]\chi_\nu(\widetilde D,\widetilde s)w + \eta_\ell \widetilde
R(y,\widetilde D,\widetilde  s)\chi_\nu(\widetilde D,\widetilde  s)w
                                                \nonumber\\
=[\Gamma,\eta_\ell]\chi_\nu(\widetilde D,\widetilde s)w+\widetilde R(y,\widetilde D,\widetilde s)
w_{\nu}+[\eta_\ell,
\widetilde R(y,\widetilde D,\widetilde s)]\chi_\nu(\widetilde D,\widetilde s)w.\nonumber
\end{eqnarray}
In order to obtain the second equality in
(\ref{lesopoval1}) we used (\ref{book}).  The short computations imply
\begin{equation}\label{lecopoval}
(\frac 1i\partial_{y_n}-\Gamma^-_\beta(y,{\widetilde D},\widetilde s))
(\frac 1i\partial_{y_n}-\Gamma^+_
\beta(y,{\widetilde D},\widetilde s))
\end{equation}
\begin{eqnarray*}
= &&-\partial_{y_n}^2-\frac 1i[\partial_{y_n},\Gamma^+_
\beta(y,{\widetilde D},\widetilde s)]+\Gamma^-_\beta(y,{\widetilde D},
\widetilde s)\Gamma^+_\beta(
y,{\widetilde D},\widetilde s)  \\
+&& {i}\Gamma^-_\beta(y,{\widetilde D},\widetilde s)\partial_{y_n}
+{i}\Gamma^+_\beta(y,{\widetilde D},\widetilde s)\partial_{y_n}.
\end{eqnarray*}
By Lemma \ref{Fops3}, we have
$$
\Gamma^-_\beta(y,{\widetilde D},\widetilde s)\Gamma^+_\beta(y,{\widetilde D},
\widetilde s)
= \Gamma(y,{\widetilde D},\widetilde s)+R_0,
$$
where
$$\Vert R_{0}\Vert_{{\mathcal L}(H_0^{\frac 12 ,1,\widetilde s}({
\Upsilon_\ell}),L^2({\Upsilon_\ell}))}
$$
$$\le C_{44}(\pi_{W^{0,\infty}( \Upsilon_\ell)}
(\Gamma^+_\beta)\pi_{W^{1,\infty}( \Upsilon_\ell)}(\Gamma^-_\beta)+\pi_{W^{1,\infty}( \Upsilon_\ell)}
(\Gamma^+_\beta)\pi_{W^{0,\infty}( \Upsilon_\ell)}(\Gamma^-_\beta))
\le C_{45}\widetilde\varphi^\frac {5}{12}( y^*_0).
$$

The commutator
$
[\partial_{y_n},\Gamma^+_\beta(y,{\widetilde D},\widetilde s)]
$ is the  pseudodifferential operator with the symbol
$\partial_{y_n}\Gamma^+_\beta(y,{\widetilde \xi},\widetilde s)$.
By Lemma \ref{Fops2} we have
$$
\Vert [\partial_{y_n},\Gamma^+_\beta(y,{\widetilde D},\widetilde s)]\Vert
_{\mathcal L (H^{\frac 12 ,1,\widetilde s}_0(\Upsilon_\ell),L^2(\Upsilon_\ell))}
\le C_{46}\widetilde\varphi^\frac {5}{12}(y^*_0).
$$
Denote
$$
R(y,\widetilde D,\widetilde s) = \left(2\vert \widetilde s\vert \varphi_n+
\frac{\sum_{j=1}^{n-1}\partial_{y_j}\vartheta (y_1,\dots,y_{n-1})(\partial_{y_j}
-\vert \widetilde s\vert \varphi_j)}{ G}\right).
$$
By (\ref{3.34}) - (\ref{3.36}), (\ref{book}) and the fact that
$\widetilde\mu_\ell(y)\eta_\ell(y)=\eta_\ell(y)$ the following is true:
\begin{eqnarray}\label{lecopoval3}
({i}\Gamma^-_\beta(y,{\widetilde D},\widetilde s)\partial_{y_n}
+{i}\Gamma^+_\beta(y,{\widetilde D},\widetilde s)\partial_{y_n}) w_{\nu}
=\widetilde\mu_\ell R(y,\widetilde D,\widetilde s)\partial_{y_n}\kappa(\widetilde D,\widetilde s)
w_{\nu}\\
=\widetilde\mu_\ell  [R(y,\widetilde D,\widetilde s)\partial_{y_n}\kappa,\eta_\ell]
\chi_\nu(\widetilde D,\widetilde s)w
+\widetilde\mu_\ell\eta_\ell
R(y,\widetilde D,\widetilde s)\partial_{y_n}\chi_\nu(\widetilde D,\widetilde s)w      \nonumber\\
=\widetilde\mu_\ell  [R(y,\widetilde D,\widetilde  s)\partial_{y_n}\kappa,\eta_\ell]
\chi_\nu(\widetilde D,\widetilde s)w
+[\eta_\ell ,R(y,\widetilde D,\widetilde s)\partial_{y_n}]\chi_\nu(\widetilde D,\widetilde s)w
+R(y,\widetilde D,\widetilde s)\partial_{y_n}w_{\nu}.
                             \nonumber
\end{eqnarray}
Since $-\partial^2_{y_n}w_{\nu}+R(y,\widetilde D,\widetilde s)\partial_{y_n}w_{\nu}
+ \widetilde R(y,\widetilde D,\widetilde s)
w_{\nu}=\frac{1}{\beta G}P_{\rho,\beta}(y,D,\widetilde s)w_{\nu}$,
setting
$$
T_\beta =-R_0+[\partial_{y_n},\Gamma^+_\beta(y,{\widetilde D},\widetilde s)]
\chi_\nu(y,\widetilde D,\widetilde s)
-[\Gamma,\eta_\ell]\chi_\nu(\widetilde D,\widetilde s)-[\eta_\ell, \widetilde R(y,\widetilde D,\widetilde s)]\chi_\nu(\widetilde D,\widetilde s)
$$
$$
-\widetilde\mu_\ell  [R(y,\widetilde D,\widetilde  s)\partial_{y_n}\kappa,\eta_\ell]\chi_\nu
(\widetilde D,\widetilde s)
-[\eta_\ell ,R(y,\widetilde D,\widetilde  s)\partial_{y_n}]\chi_\nu(\widetilde D,\widetilde s)
$$
and using (\ref{lesopoval1}) - (\ref{lecopoval3}), we obtain (\ref{min}).
Now we prove estimate (\ref{lokom}). Lemma \ref{Fops2} yields
\begin{eqnarray}\label{AQ1}
\Vert [\Gamma,\eta_\ell]\Vert
_{\mathcal L (H^{\frac 12 ,1,\widetilde  s}_0( \Upsilon_\ell),L^2( \Upsilon_\ell))}
\le C_{47}(\pi
_{C^0(\Upsilon_\ell)}(\Gamma)\pi_{C^0(\Upsilon_\ell)}(\eta_\ell)+\pi
_{C^0(\Upsilon_\ell)}(\Gamma)\pi_{W^{1,\infty}( \Upsilon_\ell)}(\eta_\ell)\nonumber\\+\pi
_{W^{1,\infty}(\Upsilon_\ell)}(\Gamma)\pi_{C^0(\Upsilon_\ell)}(\eta_\ell))\le C_{48}\widetilde\varphi^\frac {5}{12}( y^*_0).
\end{eqnarray}
For differential operators $R$ and $\widetilde R$,  we obtain the
estimates
\begin{equation}\label{LQ2}
 \Vert [\mu_\ell, \widetilde R(y,\widetilde D,\widetilde s)]\Vert
_{\mathcal L (H^{\frac 12 ,1,\widetilde s}
_0(\Upsilon_\ell),L^2(\Upsilon_\ell))}\le C_{50}\widetilde\varphi^\frac {5}{12}( y^*_0),
\end{equation}
and
\begin{eqnarray}\label{Q3}
\Vert\mu_\ell  [R(y,\widetilde D,\widetilde s)\partial_{y_n}\kappa,\eta_\ell]\Vert
_{\mathcal L (H^{\frac 12 ,1,\widetilde s}_0(\Upsilon_\ell),L^2(\Upsilon_\ell))}
                                     \nonumber\\
\le \Vert
  [R(y,\widetilde D,\widetilde s)\partial_{y_n}\kappa,\eta_\ell]\Vert
_{\mathcal L (H^{\frac 12 ,1,\widetilde s}_0( \Upsilon_\ell),L^2(\Upsilon_\ell))}
\le C_{51}\widetilde\varphi^\frac {5}{12}(y^*_0),
\end{eqnarray}
and
\begin{equation}\label{LQ1}
\Vert [\eta_\ell ,R(y,\widetilde D,\widetilde s)\partial_{y_n}]\Vert _{\mathcal L
 (H^{\frac 12 ,1,\widetilde  s}_0(\Upsilon_\ell),L^2(\Upsilon_\ell))}
\le C_{52} \Vert \eta_\ell\Vert
_{C^1(\Upsilon_\ell)}\le C_{53}\widetilde\varphi^\frac {5}{12}(y^*_0).
\end{equation}
Form (\ref{AQ1})-(\ref{LQ1}) we obtain (\ref{lokom}).
The proof of the proposition is complete. $\blacksquare$

Let $$ y^*\in\mbox{supp}\, \eta_\ell$$ and $\ell$ is fixed.

We have
\begin{proposition}\label{gopnik} Let $-\infty<\alpha<a<b<\beta<+\infty,$ $p\in \Bbb N_+$ and  $\mbox{supp} \,\mbox{\bf v}\subset I_1= [a, b]\times \Bbb R^{n-1}.$
Then there exists an independent constant $C_{54}$ such that
$$\Vert M^p(\widetilde D, s)\mbox{\bf v}
\Vert_{L^2 (B(0,R)\setminus [\alpha,\beta]\times \Bbb R^{n-1})}\le \frac{C_{54}}{(\min \{a-\alpha, \beta-b\})^p}\Vert \mbox{\bf v}\Vert
_{L^2(\Bbb R^n)}.
$$
\end{proposition}

{\bf Proof.}  Let $\widetilde y \in B(0,R)\setminus ((\alpha,\beta)\times {\Bbb R}^{n-1}).$   Consider two cases.
{\bf Case 1.}   Let  $\vert y_0 \vert  > \mbox{max}\,\{ \vert a\vert, \vert b\vert\}=\tilde c.$ Integrating by parts we have

$$
 M(\widetilde D, s)\mbox{\bf v}=
\frac{1}{(2\pi)^n}\int_{\Bbb R^n\times \{supp\mbox{\bf v}\}} \frac{\partial^p_{\xi_0}M^p(\widetilde\xi,s)}{(i(x_0-y_0))^p} e^{i< \widetilde y- \widetilde x,\widetilde\xi>}\mbox{\bf v}(\widetilde x)d \tilde\xi d\widetilde x=
$$
%$$
%\sum_{k=1}^3\frac{1}{(2\pi)^n}\left \vert\int_{\Bbb R^n\times \mathcal Z_k} \frac{\partial_{\xi_0}M(\widetilde\xi,s)}{i(x_0-y_0)} e^{i<\widetilde y-\widetilde %x,\widetilde \xi>}\mbox{\bf u}(x)d\widetilde\xi d\widetilde x\right\vert,
%$$
$$
-\frac{1}{(2\pi)^ny_0i^p}\int_{\Bbb R^n\times \{supp\mbox{\bf v}\}} \frac{\partial^p_{\xi_0}M^p(\widetilde\xi,s)}{(1-\frac{x_0}{y_0})^p} e^{i<\widetilde y-\widetilde x,\widetilde \xi>}\mbox{\bf v}(\widetilde x)d\widetilde\xi d\widetilde x=
$$
$$
\frac{i^p}{y_0^p}\frac{(-1)^{p+1}}{(2\pi)^n}\sum_{k=1}^\infty \int_{\Bbb R^n\times \{supp\mbox{\bf v}\} } c_k(p)(\frac{x_0}{y_0})^k\partial^p_{\xi_0}M^p(\widetilde\xi,s) e^{i<\widetilde y-\widetilde x,\widetilde \xi>}\chi_1\mbox{\bf v}(\widetilde x)d\widetilde\xi d\widetilde x=
$$
$$
\frac{i^p}{y_0^p}\frac{(-1)^{p+1}}{(2\pi)^\frac n2 }\sum_{k=1}^\infty c_k(p) K(\tilde y,\widetilde D,s)(\left (\frac{x_0}{y_0}\right )^k\mbox{\bf v}).
$$

Therefore
$$
\Vert  M^p(\widetilde D,s)\mbox{\bf u}\Vert_{L^2(B(0,R)\setminus  [\alpha,\beta]\times \Bbb R^{n-1})}\le \sup_{\widetilde y\in B(0,R)\setminus [\alpha,\beta]\times \Bbb R^{n-1})}\left\{\frac{C_{55}}{\vert y_0\vert ^p}\frac{1}{(2\pi)^\frac n 2} \sum_{k=1}^\infty c_k\Vert  (\frac{x_0}{y_0} )^k\chi_1\mbox{\bf v}\Vert_{L^2(\Bbb R^n)}\right\}
$$
$$\le  \frac{C_{56}}{\vert y_0\vert } \sum_{k=1}^\infty c_k(p) (\frac{\tilde c}{y_0} )^k\Vert  \mbox{\bf v}\Vert_{L^2(\Bbb R^n)}
\le C_{57} \frac{\Vert  \mbox{\bf v}\Vert_{L^2(\Bbb R^n)}}{ (\min \{a-\alpha, \beta-b\})^p }.
$$

{\bf Case 2.}
Let  $\vert y_0 \vert  \le  \mbox{max}\,\{ \vert a\vert, \vert b\vert\} .$ If $\vert a\vert >\vert b\vert$ then
 $y_0\in (b, \vert a\vert] $ and we set $\tilde c=b.$ On the other hand if  $\vert a\vert <\vert b\vert$ then
 $y_0\in (\vert b\vert ,  a] $ and we set $\tilde c=a.$ In both cases $\vert x_0-y_0+1\vert>1$ on $[a,b].$
The short computations imply
$$
 M^p(\widetilde D, s)\mbox{\bf v}=
\frac{1}{(2\pi)^n}\int_{\Bbb R^n\times  \{supp\mbox{\bf v}\}} \frac{\partial^p_{\xi_0}M^p(\widetilde\xi,s)}{i^p(x_0-y_0)^p} e^{i< \widetilde y- \widetilde x,\widetilde\xi>}\mbox{\bf v}(\widetilde x)d \tilde\xi d\widetilde x=
$$
%$$
%\sum_{k=1}^3\frac{1}{(2\pi)^n}\left \vert\int_{\Bbb R^n\times \mathcal Z_k} \frac{\partial_{\xi_0}M(\widetilde\xi,s)}{i(x_0-y_0)} e^{i<\widetilde y-\widetilde %x,\widetilde \xi>}\mbox{\bf u}(x)d\widetilde\xi d\widetilde x\right\vert,
%$$
$$
-\frac{1}{(2\pi)^ni^p}\int_{\Bbb R^n\times  \{supp\mbox{\bf v}\} } \frac{\partial^p_{\xi_0}M^p(\widetilde\xi,s)}{(x_0-y_0+1)^p(1-\frac{1}{x_0-y_0+1})^p} e^{i<\widetilde y-\widetilde x,\widetilde \xi>}\mbox{\bf v}(\widetilde x)d\widetilde\xi d\widetilde x=
$$
$$
\frac{(-1)^{p+1}i^p}{(2\pi)^n}\sum_{k=1}^\infty \int_{\Bbb R^n\times \{supp\mbox{\bf v}\} } \frac{1}{\vert x_0-y_0+1\vert^p} c_k(p)(\frac{1}{x_0-y_0+1})^{k+1}\partial^p_{\xi_0}M^p(\widetilde\xi,s) e^{i<\widetilde y-\widetilde x,\widetilde \xi>}\mbox{\bf v}(\widetilde x)d\widetilde\xi d\widetilde x=
$$
$$
\frac{i}{(2\pi)^\frac n2 }\sum_{k=1}^\infty c_k(p) K(\tilde y,\widetilde D, s)(\frac{1}{\vert x_0-y_0+1\vert^p}(\frac{1}{x_0-y_0+1})^{k+1}\mbox{\bf v}).
$$

Therefore
$$
\Vert  M^p(\widetilde D,s)\mbox{\bf u}\Vert_{L^2(B(0,R)\setminus  [\alpha,\beta]\times \Bbb R^{n-1})}
$$
$$\le \sup_{\widetilde y\in B(0,R)\setminus  [\alpha,\beta]\times \Bbb R^{n-1})}\left\{\frac{C_{58}}{(2\pi)^\frac n 2} \sum_{k=1}^\infty c_k\Vert \frac{1}{\vert x_0-y_0+1\vert^p} (\frac{1}{x_0-y_0+1})^{k+1} \mbox{\bf v}\Vert_{L^2(\Bbb R^n)}\right\}
$$
$$\le  \frac{C_{59}}{\vert \tilde c-y_0+1\vert^p } \sum_{k=1}^\infty c_k(p) (\frac{1}{\tilde c-y_0+1})^{k+1}\Vert  \mbox{\bf v}\Vert_{L^2(\Bbb R^n)}
\le C_{60} \frac{\Vert  \mbox{\bf v}\Vert_{L^2(\Bbb R^n)}}{(\min \{a- y_0, y_0-b\})^p  }.
$$

Proof of the proposition is complete.
$\blacksquare$

We apply the  Proposition \ref{gopnik} in order to estimate  the  $H^{\frac 12,1,\widetilde s}$ norm of the function $(1-\eta_\ell)\chi_\nu(\widetilde D, \widetilde s) \mbox{\bf w}.$ Observe that by (\ref{capitol}) for all sufficiently large  $\ell$
$$
\mbox{supp}\, \mbox{\bf w}\subset [-T+(\ell+2)^{-4}, -T+(\ell-2)^{-4}]\cup [T-(\ell-2)^{-4}, T-(\ell+2)^{-4}]\times \Bbb R^n
$$
and
$$
\mbox{supp}\,(1- \eta_\ell)\subset [-T,T]\setminus[-T+(\ell+11)^{-4}, -T+(\ell-11)^{-4}]\cup [T-(\ell-11)^{-4}, T-(\ell+11)^{-4}].
$$
Therefore
\begin{equation}\label{cosmos}
\Vert (1-\eta_\ell)\chi_\nu(\widetilde D, \widetilde s) \mbox{\bf w}
\Vert_{H^{\frac 12,1,\widetilde s} (B(0,R))}\le C_{61}\ell^{5}\Vert \mbox{\bf w}\Vert
_{L^2(\Bbb R^n)} \le C_{62}\widetilde\varphi^\frac{5}{12}(y^*)\Vert \mbox{\bf w}\Vert
_{L^2(\Bbb R^n)}.
\end{equation}
%By (\ref{lomka1}), (\ref{pauk}), (\ref{knight0}) and  Proposition 9.3 p. 95 of \cite{Tay1} we have
%By (\ref{gabon})
%\begin{equation}\label{po}
%\vert \widetilde s\vert\Vert (1-\eta_\ell)\chi_\nu(\widetilde D, \widetilde s) \mbox{\bf w}
%\Vert_{H^{\frac 12,1,\widetilde s} (B(0,R))}\le C_{57}\widetilde \varphi(y^*_0)\Vert \mbox{\bf w}\Vert
%_{L^2(\Bbb R^n)}.
%\end{equation}
By arguments, same as in Proposition \ref{Fops4} we obtain
\begin{equation}\label{po1}
\vert \widetilde s\vert\Vert (1-\eta_\ell)\chi_\nu(\widetilde D, \widetilde s)
\mbox{\bf w}\Vert_{H^{\frac 12,1,\widetilde s} (\Bbb R^n\setminus B(0,R))}\le{ C_{63}}
\Vert \mbox{\bf w}\Vert_{L^2(\Bbb R^n)}.
\end{equation}
By (\ref{cosmos}) and (\ref{po1})
\begin{equation}\label{nikolas}\Vert (1-\eta_\ell)\chi_\nu(\widetilde D, \widetilde s) \mbox{\bf w}
\Vert_{H^{\frac 12 ,1,\widetilde s} (\Bbb R^n)}\le C_{64}\widetilde \varphi(y^*_0)\Vert
\mbox{\bf w}\Vert_{L^2(\Bbb R^n)}.
\end{equation}

Denote
$$
V^\pm_{\mu}(k,j)=(\frac 1i \partial_{y_n}-\Gamma^\pm_{\mu}(y,{\widetilde D},
\widetilde s))
w_{k,j,\nu}, \quad
V^\pm_{\lambda+2\mu}=(\frac 1i \partial_{y_n}-\Gamma^\pm_{\lambda+2\mu}
(y,{\widetilde D},\widetilde s))w_{2,\nu}.
$$

Let us consider the equation
\begin{equation}\label{barmalei0}
(\frac 1i \partial_{y_n}-\Gamma^-_\beta(y,{\widetilde D},\widetilde s))V=\mbox{\bf p},
\quad V\vert_{y_n=\gamma}=0.
\end{equation}
For solutions of this problem, we can prove an a priori estimate.
\begin{proposition}\label{gorokx} Let $r_\beta(y^*,\zeta^*)\ne 0$
and $V=V^+_{\mu}(k,j)$ if $\beta =\mu$ and $V=V^+_{\lambda+2\mu}$ if $\beta
=\lambda+2\mu.$ There exists a
constant $C_{61}>0$ such that
\begin{equation} \label{min1}
\Vert V(\cdot, 0)\Vert_{H^{\frac 14,\frac 12,\widetilde s}({\Bbb R}^n)}\le
C_{61}(\tilde\varphi^\frac{5}{12}(y^*)\Vert \mbox{\bf w}\Vert_{H^{\frac 12,1,\widetilde s}(\mathcal Q)}+\Vert \mbox{\bf p}\Vert_{L^2(\mathcal Q)}).
\end{equation}
\end{proposition}

${\bf Proof.}$  Let $M( \widetilde D,\widetilde s)$ be the pseudodifferential operator with the symbol $M(\widetilde \xi, \widetilde s).$ We taking the scalar product of the equation (\ref{barmalei0})
and the function $-i\overline{M( \widetilde D,\widetilde s) V}$ in $L^2(\Bbb R^n)$ and integrating by parts
we obtain
\begin{eqnarray}\label{fosset}
\frac 12\Vert M^\frac 12( \widetilde D,\widetilde s) V(\cdot, 0)\Vert^2_{L^2(\Bbb R^n)}\quad\quad\quad\quad\quad\\+i\int_0^\gamma(\Gamma^-_\beta
(y,{\widetilde D},\widetilde s)V,\overline{M( \widetilde D,\widetilde s)V})_{L^2(\Bbb R^n)}dy_n=-i\int_0^\gamma(\mbox{\bf p},\overline{M( \widetilde D,\widetilde s)V})
_{L^2(\Bbb R^n)}dy_n.\nonumber
\end{eqnarray}
By (\ref{3.35})- (\ref{3.37}) and assumption $r_\beta(y^*,\zeta^*)\ne 0$
for some positive constant $C_{62}$
\begin{equation}\label{3.38}
-\text{Im}\,  \Gamma^-_\beta(y,\widetilde \xi,\widetilde s)
\ge  C_{62}M(\widetilde \xi, \widetilde s)
\quad \forall (y,\widetilde \xi,\widetilde s)\in \Upsilon_\ell\times
\mathcal O(y^*,\delta_1(y^*)),
\end{equation} where $\Upsilon_\ell$ is determined by (\ref{zontik}).

We set $\tilde  V(y)=\eta^*_1(y)\sum_{k=\ell-17}^{\ell+17} \kappa_k(y_0) V(y), \eta^*_1\in C^\infty_0(B(y^*,\frac{93}{50}\delta)), \eta_1\vert_{B(y^*,\frac{91}{50}\delta)}=1,$  \newline $ \tilde\mu_* =\eta^*_2(y)\sum_{k=\ell-24}^{\ell+34} \kappa_k(y_0),\eta^*_2\in C^\infty_0(B(y^*,2\delta)), \eta_2\vert_{\mbox{supp}\, (B(y^*,\frac{98}{50}\delta)) }=1,$ where function $\eta_*$ is defined in (\ref{bobo}). The short computations imply

$$
-\int_0^\gamma\mbox{Im}(\Gamma^-_\beta(y,{\widetilde D},\widetilde s)V,\overline{M( \widetilde D,\widetilde s)V})
_{L^2(\Bbb R^n)}dy_n=
-\int_0^\gamma\mbox{Im}(\tilde \mu_* \Gamma^-_\beta(y,{\widetilde D},\widetilde s)V,\overline{M( \widetilde D,\widetilde s)V})
_{L^2(\Bbb R^n)}dy_n=
$$
$$-\int_0^\gamma\mbox{Im}(\Gamma^-_\beta(y,{\widetilde D},\widetilde s)\tilde V,\overline{\tilde \mu_* M( \widetilde D,\widetilde s)\tilde V})
_{L^2(\Bbb R^n)}dy_n
$$
$$
-\int_0^\gamma\mbox{Im}(\Gamma^-_\beta(y,{\widetilde D},\widetilde s)(V-\tilde V),\overline{\tilde \mu_* M( \widetilde D,\widetilde s)\tilde V})
_{L^2(\Bbb R^n)}dy_n
$$
$$
-\int_0^\gamma\mbox{Im}(\Gamma^-_\beta(y,{\widetilde D},\widetilde s)\tilde V,\overline{\tilde \mu_* M( \widetilde D,\widetilde s)(V-\tilde V)})
_{L^2(\Bbb R^n)}dy_n=\sum_{j=1}^3 \mathcal I_j.
$$

We estimate each term in the above inequality separately.
Observe that $\tilde\mu_*\vert_{\mbox{supp}\, \tilde \mu}=1.$
By G\aa rding's inequality (\ref{goblin}) there exists a positive constant $C_{63}$
such that
\begin{eqnarray}\label{qupol11}
-\int_0^\gamma\mbox{Im}(\Gamma^-_\beta(y,{\widetilde D},\widetilde s)\tilde V,\overline{\tilde \mu_* M( \widetilde D,\widetilde s)\tilde V})
_{L^2(\Bbb R^n)}dy_n=-\int_0^\gamma\mbox{Im}(M( \widetilde D,\widetilde s)\tilde \mu\Gamma^-_\beta(y,{\widetilde D},\widetilde s)\tilde V,\overline{\tilde V})
_{L^2(\Bbb R^n)}dy_n=\nonumber\\
-\int_0^\gamma\mbox{Im}(\tilde \mu_* M( \widetilde D,\widetilde s)\tilde \mu_* \Gamma^-_\beta(y,{\widetilde D},\widetilde s)\tilde V,\overline{\tilde V})
_{L^2(\Bbb R^n)}dy_n\ge
\nonumber\\
\int_0^\gamma\mbox{Re}(p(y,{\widetilde D},\widetilde s)\tilde V,\overline{\tilde V})
_{L^2(\Bbb R^n)}dy_n=
\nonumber\\ C_{63}
\int_0^\gamma\Vert \tilde V(\cdot,y_n)\Vert^2_{H^{\frac 12,1,\widetilde s}(\Bbb R^n)}
dy_n
-C_{64}\ell^{10}  \Vert \tilde V\Vert^2_{L^2(\mathcal Q)}\ge
\nonumber\\ C_{63}
\int_0^\gamma\Vert \tilde V(\cdot,y_n)\Vert^2_{H^{\frac 12,1,\widetilde s}(\Bbb R^n)}
dy_n
-C_{65}\tilde \varphi^{\frac 5 6}(y^*) \Vert \tilde V\Vert^2_{L^2(\mathcal Q)}.
\end{eqnarray}
Here $p(y,\widetilde D,\widetilde s)$ is the operator  with symbol
$
 p(y,\widetilde\xi,\widetilde s)=i\tilde\mu_* M( \widetilde \xi,\widetilde s)\tilde \mu_* \Gamma^-_\beta(y,{\widetilde \xi},\widetilde s).
$
where
$$\mbox{ Re}p(y,\widetilde\xi,\widetilde s)=-\tilde\mu^2_*M(( \widetilde \xi,\widetilde s)\mbox{Im}\, \Gamma^-_\beta(y,{\widetilde \xi},\widetilde s)\ge C_{66}\tilde\mu M^2( \widetilde \xi,\widetilde s)\quad\forall y\in B(y^*,2\delta(y^*)).
$$
$$
(\pi_{W^{0,\infty}(\mathcal O_3)}(p)+1)\pi_{C^{0}(\mathcal O_3)}
(\widetilde \gamma)+\sum_{k=0}^1(\pi_{W^{k,\infty}(\mathcal O_3)}(p)+1)\pi_{C^{1-k}(\mathcal O_3)}(\widetilde\gamma)+\pi_{C^{1}(\mathcal O_3)}(\widetilde\gamma)]^2\le C_{67}\tilde\varphi^\frac 56(y^*),
$$ and
$\mathcal O=B(y^*,2\delta(y^*)),\mathcal O_1= (B(y^*,\frac{93}{50}\delta)) \cap \mbox{supp}\,\sum_{k=\ell-17}^{\ell+17} \kappa_k(y_0), \mathcal O_2=B(y^*,\frac{97}{50}\delta)) \cap \mbox{supp}\,\sum_{k=\ell-19}^{\ell+19} \kappa_k(y_0), \eta_2\vert_{B(y^*,\frac{99}{50}\delta)}=1, \mathcal O_3,\tilde\gamma=\sum_{k=\ell-23}^{\ell+23} \kappa_k(y_0)\eta_2, \eta_2\vert_{B(y^*,\frac{98}{50}\delta)}=1, \eta_2\in C^\infty_0(B(y^*,\frac{99}{50}\delta)). $
Using Lemma \ref{Fops1} and Lemma \ref{Fops0}  and  Proposition \ref{gopnik} we obtain
\begin{equation}\label{Aqupol1}
\vert \mathcal I_3\vert\le \frac{C_{68}}{4}
\int_0^\gamma\Vert \tilde V(\cdot,y_n)\Vert^2_{H^{\frac 12,1,\widetilde s}(\Bbb R^n)}dy_n-C_{69}\Vert\tilde \mu  M(V- \tilde V)\Vert^2_{L^2(\mathcal Q)}.
\end{equation}
On the other hand by Proposition \ref{Fops2} and Proposition \ref{gopnik}
\begin{eqnarray}\label{sun}
\Vert\tilde\mu M(V- \tilde V)\Vert_{L^2(\mathcal Q)}\le \Vert \tilde\mu M(\gamma (\frac{1}{i}\partial_{y_n} \mbox{\bf w}-\Gamma^+_\beta \mbox{\bf w}))\Vert_{L^2(\mathcal Q)}+  \Vert M[\Gamma^+_\beta,\gamma] \mbox{\bf w})\Vert_{L^2(\mathcal Q)}
\\
\le\Vert\gamma M( (\frac{1}{i}\partial_{y_n} \mbox{\bf w}-\Gamma^+_\beta \mbox{\bf w}))\Vert_{L^2(\mathcal Q)}+\Vert[M,\gamma]( (\frac{1}{i}\partial_{y_n} \mbox{\bf w}-\Gamma^+_\beta \mbox{\bf w})\Vert_{L^2(\mathcal Q)}+ C_{70}\Vert [\Gamma^+_\beta,\gamma] \mbox{\bf w})\Vert_{H^{\frac 12, 1,\tilde s}(\mathcal Q)}
\nonumber\\
%\le\Vert\gamma M( (\frac{1}{i}\partial_{y_n} \mbox{\bf w}-\Gamma^+_\beta \mbox{\bf w}))\Vert_{L^2(\mathcal Q)}+\Vert[M,\gamma]( %(\frac{1}{i}\partial_{y_n} \mbox{\bf w}-\Gamma^+_\beta \mbox{\bf w})\Vert_{L^2(\mathcal Q)}+ \Vert M[\Gamma^+_\beta,\gamma] \mbox{\bf %w})\Vert_{L^2(\mathcal Q)}.
%\nonumber\\
\le C_{71}( \tilde\varphi^{\frac {5}{12}}(y^*)\Vert(\frac{1}{i}\partial_{y_n} \mbox{\bf w}-\Gamma^+_\beta \mbox{\bf w}))\Vert_{L^2(\mathcal Q)}+ \tilde\varphi^{\frac {5}{12}}(y^*)\Vert \mbox{\bf w}\Vert_{H^{\frac 12, 1, \tilde s}(\mathcal Q)}).\nonumber
\end{eqnarray}
Therefore
\begin{equation}\label{line}
\vert \mathcal I_3\vert\le \frac{C_{65}}{4}
\int_0^\gamma\Vert \tilde V(\cdot,y_n)\Vert^2_{H^{\frac 12,1,\widetilde s}(\Bbb R^n)}dy_n-C_{72}\tilde\varphi^{\frac {5}{6}}(y^*)\Vert \mbox{\bf w}\Vert^2_{H^{\frac 12, 1, \tilde s}(\mathcal Q)}.
\end{equation}

The short computations imply and Lemma \ref{Fops1} imply
$$
\vert \mathcal I_2\vert=\vert \int_0^\gamma\mbox{Im}((V-\tilde V),\overline{{\Gamma^{-,*}_\beta}( \widetilde D,\widetilde s)\tilde\mu M({\widetilde D},\widetilde s)\tilde V})
_{L^2(\Bbb R^n)}dy_n
$$
$$+\int_0^\gamma\mbox{Im}((V-\tilde V),\overline{R\tilde\mu M\tilde V})
_{L^2(\Bbb R^n)}dy_n\vert\le
$$
$$
\vert \int_0^\gamma\mbox{Im}(\tilde\mu M({\widetilde D},\widetilde s)(V-\tilde V),\overline{{\Gamma^{-,*}_\beta}( \widetilde D,\widetilde s)\tilde V})
_{L^2(\Bbb R^n)}dy_n\vert
$$
$$
+\vert \int_0^\gamma\mbox{Im}((V-\tilde V),\overline{[{\Gamma^{-,*}_\beta},\tilde\mu_* M]\tilde V})\vert
_{L^2(\Bbb R^n)}dy_n
$$
$$+\vert\int_0^\gamma\mbox{Im}((V-\tilde V),\overline{R \tilde\mu_* M\tilde V})
_{L^2(\Bbb R^n)}dy_n\vert
=\sum_{k=1}^3 \mathcal J_k.
$$
Using Lemma \ref{Fops1} and Lemma \ref{Fops0}  and  (\ref{sun}) we obtain
\begin{equation}\label{zopa}
\mathcal J_1\le \frac{C_{65}}{20}
\int_0^\gamma\Vert \tilde V(\cdot,y_n)\Vert^2_{H^{\frac 12,1,\widetilde s}(\Bbb R^n)}dy_n-C_{73}\tilde \varphi^{\frac 5 6}(y^*) \Vert \mbox{\bf w} \Vert^2_{L^2(\mathcal Q)}
\end{equation}

Proposition \ref{Fops2} implies
\begin{equation}\label{zopa1}
\mathcal J_2\le \frac{C_{65}}{20}
\Vert \tilde V\Vert^2_{L^2(\mathcal Q)}-C_{74}\tilde \varphi^{\frac 5 6}(y^*) \Vert  \mbox{\bf w}\Vert^2_{H^{\frac 12,1,\tilde s}(\mathcal Q)}.
\end{equation}

 By Proposition \ref{Fops1} we have
\begin{equation}\label{zopa2}
\mathcal J_3\le \frac{C_{65}}{2}
\Vert \tilde V\Vert^2_{L^2(\mathcal Q)}-C_{75}\Vert  \mbox{\bf w}\Vert^2_{H^{\frac 12,1,\tilde s}(\mathcal Q)}.
\end{equation}
By  (\ref{qupol11}),  (\ref{Aqupol1}), (\ref{line}), (\ref{zopa})- (\ref{zopa2})
$$
-\int_0^\gamma\mbox{Im}(\Gamma^-_\beta(y,{\widetilde D},\widetilde s)V,\overline{M( \widetilde D,\widetilde s)V})
_{L^2(\Bbb R^n)}dy_n\ge -C_{76}\Vert  \mbox{\bf w}\Vert^2_{H^{\frac 12,1,\tilde s}(\mathcal Q)}.
$$
This inequality, (\ref{fosset}) and Proposition \ref{gorokx1} imply
(\ref{min1}).
$\blacksquare$

%Now we apply the Lemma \ref{Fops4} to estimate $(1-\eta_\ell)\chi_\nu(\widetilde D,
%\widetilde s) \mbox{\bf w}$.
% Set $\mathcal O_1
%=\mbox{supp}\, \kappa_\ell\cap B(0,R)$, $\mathcal O=B(0,2R)$ there radius $R$
%is sufficiently large, $\mathcal O_2=supp(1-\eta_\ell)\cap B(0,R).$
%We claim that for all sufficiently large $\ell$ one can find a constant $C_{54}$ independent of $\ell$ such that
%\begin{equation}\label{POPO}
%\mbox{dist}\, (\mathcal O_1,\mathcal O_2)\ge C_{54}\ell^{-\frac{9}{5}}.
%\end{equation}
%Indeed, if $y\in \mathcal O_1$ then $\frac 12\le 2^{-\ell} 2^\frac{1}{\theta^\frac 54(y_0)}\le 2.$ This inequality implies $(\ell+1)^{-\frac %54}\le \theta(y_0)\le (\ell-1)^{-\frac 54}.$ Provided that $\ell\ge \ell_0$ where $\ell_0$ taken sufficiently large we have
%\begin{equation}\label{Iznarok} (\ell+2)^{-\frac 45}\le T+y_0\le (\ell-2)^{-\frac 45}
%\quad \mbox{or}\,\,(\ell+2)^{-\frac 45}\le T-y_0\le (\ell-2)^{-\frac 45}\quad y\in  B(y^*,\delta)\cap \mbox{supp}\, \eta_\ell.
%\end{equation}

%On the other hand the function $\eta_\ell$ is identically equal to $1$ on the set $\{y\in \Bbb R^{n+1}\vert \frac{1}{(\ell+4)^\frac 45}\le %\theta(y_0)\le \frac{1}{(\ell-4)^\frac 45}\}.$
%Hence, provided that $\ell\ge \ell_0$
%\begin{eqnarray}\label{IIznarok}
%\mathcal O_2\subset\mathcal O_3=\\\{y\in \Bbb R^{n+1}\vert T-y_0\ge (\ell-4)^{-\frac 45}\,\, \mbox{or}\,  T-y_0\le (\ell+4)^{-\frac 45}\mbox{or}\,  %T+y_0\le (\ell+4)^{-\frac 45}\,\,\mbox{or}\,\, T+y_0\ge (\ell-4)^{-\frac 45}\}.\nonumber
%\end{eqnarray}
%From (\ref{Iznarok}) and  (\ref{IIznarok}) we obtain (\ref{POPO}).

We will separately consider the two cases
$r_\mu(y^*,\zeta^*)=0$ in Section \ref{gromA} and $r_{\lambda+2\mu}(y^*,\zeta^*)=0$
in Section \ref{grom2}.

\bigskip
\section{ Case $r_{\mu}(y^*,\zeta^*)=0.$}\label{gromA}
\bigskip

In this section, we treat the case when
$\text{supp}\,\chi_\nu\subset \mathcal O(y^*,\delta_1(y^*)),$ and
$(y^*,\zeta^*)$ be a point on $\Bbb R^{n+1}\times \Bbb M$ such that
$r_\mu(y^*,\zeta^*)=0.$ By (\ref{normal}) and (\ref{3.37}) this equality implies
\begin{equation}\label{medved} \sum_{j=1}^{n-1}({\xi_j^*})^2=(\widetilde s^*)^2\sum_{j=1}^{n-1} \vert\partial_{y_j}{\varphi}(y^*)\vert^2\quad
\mbox{and}\quad \frac{\rho(y^*)\xi_0^*}{\mu(y^*)}+\widetilde s^*\sum_{j=1}^{n-1}\xi^*_j\partial_{y_j}{\varphi}(y^*)=0.\end{equation}
By (\ref{medved}) and (\ref{pinok220})
\begin{equation}\label{medved2}
\widetilde s^*\ne 0.
\end{equation}

By (\ref{medved}) there exists $C_1>0$ such that for all $(y^*,\zeta)$ from $\mathcal O(y^*,\delta_1(y^*))$ we have
\begin{eqnarray}\label{2.7} \left \vert \sum_{j=1}^{n-1}{\xi_j^2}- \widetilde s^2\sum_{j=1}^{n-1} \vert\partial_{y_j}{\varphi}(y^*)\vert^2\right\vert+
\left\vert \frac{\rho(y^*)\xi_0}{\mu(y^*)}+\widetilde s\sum_{j=1}^{n-1}\xi_j\partial_{y_j}{\varphi}(y^*)\right\vert \le \delta_1C_1M(\widetilde \xi,\widetilde s).
\end{eqnarray}
Hence, by (\ref{2.7}) and (\ref {nikolas}) for some independent constants
$C_2,C_3$
\begin{eqnarray}\label{2.99}
\vert \frak J_3(\mu, w_{k,j,\nu})\vert\le C_2\delta_1\Vert
(\partial_{y_n} \mbox{\bf w}_{\nu}(\cdot,0),\mbox{\bf  w}
_{\nu}(\cdot,0))\Vert^2
_{L^2(\Bbb R^n)\times H^{0,1,\widetilde s}(\Bbb R^n)}\nonumber\\+C_3\widetilde \varphi( y_0^*)
\Vert
(\partial_{y_n} \mbox{\bf w}(\cdot,0),\mbox{\bf  w}(\cdot,0))\Vert^2
_{L^2(\Bbb R^n)\times H^{0,1,\widetilde s}(\Bbb R^n)}.
\end{eqnarray}

We recall that by (\ref{klop}) there exist $C_4>0$ and $C_5>0$ such
that
\begin{eqnarray}\label{2.1}C_4(\vert \widetilde s\vert \Vert  w_{k,j,\nu}
\Vert^2_{H^{0,1}(\mathcal Q)}+\vert \widetilde s\vert^3\Vert
w_{k,j,\nu}\Vert^2_{L^2(\mathcal Q)})+\varXi_\mu(w_{k,j,\nu})\le C_5\Vert
P_\mu(y,D,\widetilde s) w_{k,j,\nu} \Vert^2_{L^2(\mathcal Q)}\nonumber\\
+\epsilon(\delta)\Vert
(\partial_{y_n} \mbox{\bf w}_{\nu}(\cdot,0),\mbox{\bf  w}
_{\nu}(\cdot,0))\Vert^2
_{H^{\frac{1}{4},\frac 12,\widetilde s}(\Bbb R^n)\times H^{\frac{3}{4},\frac 32,\widetilde s}(\Bbb R^n)} ,
\end{eqnarray}
where $\epsilon(\delta)\rightarrow 0$ as $\delta\rightarrow +0.$

Since by (\ref{medved2})  and (\ref {nikolas}) $s^*\ne 0$
inequality (\ref{2.7}) yields
\begin{eqnarray}\label{2.9}
\vert \frak J_2(\mu, w_{j,n,\nu})\vert
\le \frac{C_6\delta_1\mu(y^*)\vert\widetilde s\vert}{\vert\widetilde s^*\vert}
\Vert
(\partial_{y_n} \mbox{\bf w}_{\nu}(\cdot,0),\mbox{\bf  w}
_{\nu}(\cdot,0))\Vert^2
_{L^2(\Bbb R^n)\times H^{0,1,\widetilde s}(\Bbb R^n)}\nonumber\\+C_7\widetilde \varphi( y_0^*)
\Vert
(\partial_{y_n} \mbox{\bf w}(\cdot,0),\mbox{\bf  w}(\cdot,0))\Vert^2
_{L^2(\Bbb R^n)\times H^{0,1,\widetilde s}(\Bbb R^n)}.
\end{eqnarray}
Hence from  (\ref{2.9}), (\ref{2.99}) we obtain
\begin{equation} \label{2.4}
\varXi_{\mu}(\mbox{\bf w}_{1,\nu}) \ge C_{8}\int_{{\Bbb R^n}} (\vert \widetilde
s\vert\sum_{j=1}^n
\vert\partial_{y_j}\mbox{\bf  w}_{1,\nu}\vert^2+\vert \widetilde s\vert^3
\vert \mbox{\bf  w}_{1,\nu}\vert^2)(\widetilde y,0) d\widetilde y
\end{equation}
$$
- \epsilon(\delta)\vert\widetilde s\vert \Vert
(\partial_{y_n} \mbox{\bf w}_{\nu}(\cdot,0),\mbox{\bf  w}
_{\nu}(\cdot,0))\Vert^2
_{L^2(\Bbb R^n)\times H^{0,1,\widetilde s}(\Bbb R^n)}-C_{9}\widetilde \varphi( y_0^*)\Vert
(\partial_{y_n} \mbox{\bf w}(\cdot,0),\mbox{\bf  w}(\cdot,0))\Vert^2
_{L^2(\Bbb R^n)\times H^{0,1,\widetilde s}(\Bbb R^n)},
$$
where $\epsilon(\delta)\rightarrow 0$ as $\delta\rightarrow +0.$
Now we consider two subcases:

{\bf Subcase A.} Let $r_{\lambda+2\mu}(y^*,\zeta^*)=0.$

Similarly to (\ref{2.4}), we obtain
\begin{equation} \label{2.444}
\varXi_{\lambda+2\mu}(w_{2,\nu}) \ge C_{10}\int_{{\Bbb R^n}} (\vert \widetilde
s \vert \sum_{j=1}^n
\vert\partial_{y_j} w_{2,\nu}\vert^2+\vert \widetilde s\vert^3
\vert w_{2,\nu}\vert^2)(\widetilde y,0) d\widetilde y
\end{equation}
$$
- \epsilon(\delta)\vert\widetilde s\vert \Vert
(\partial_{y_n} \mbox{\bf w}_{\nu}(\cdot,0),\mbox{\bf  w}
_{\nu}(\cdot,0))\Vert^2
_{L^2(\Bbb R^n)\times H^{0,1,\widetilde s}(\Bbb R^n)}-C_{11}\widetilde \varphi( y_0^*)\Vert
(\partial_{y_n} \mbox{\bf w}(\cdot,0),\mbox{\bf  w}(\cdot,0))\Vert^2
_{L^2(\Bbb R^n)\times H^{0,1,\widetilde s}(\Bbb R^n)},
$$
where $\epsilon(\delta)\rightarrow 0$ as $\delta\rightarrow +0.$
Combining estimates  (\ref{2.4}) and (\ref{2.444}), we have
\begin{eqnarray} \label{pulemetQ}
\root\of{\vert\widetilde s\vert }\Vert
(\partial_{y_n} \mbox{\bf w}_{\nu}(\cdot,0),\mbox{\bf  w}
_{\nu}(\cdot,0))\Vert
_{L^2(\Bbb R^n)\times H^{0,1,\widetilde s}(\Bbb R^n)} + \root\of{\vert\widetilde s\vert}
\Vert
\mbox{\bf w}_{\nu}\Vert_{H^{0,1,\widetilde s}(\mathcal Q)}\le
C_{12}(\widetilde\varphi^\frac {5}{12}( y^*_0)\Vert \text{\bf w}
\Vert_{H^{\frac 12,1,\widetilde s}(\mathcal Q)}                   \nonumber\\
+\root\of{ \vert\widetilde s\vert}\Vert
\mbox{\bf g}e^{\vert s\vert\varphi}\Vert_{L^2({\Bbb R^n})}
+\Vert \mbox{\bf P}(y,D,\widetilde s)\mbox{\bf
w}_{\nu}\Vert_{L^2(\mathcal Q)}
\nonumber\\+\root\of{\widetilde \varphi( y_0^*)}\Vert
(\partial_{y_n} \mbox{\bf w}(\cdot,0),\mbox{\bf  w}(\cdot,0))\Vert
_{L^2(\Bbb R^n)\times H^{0,1,\widetilde s}(\Bbb R^n)}).
\end{eqnarray}

By (\ref{pulemetQ}), (\ref{medved2}) and (\ref{nikolas})
\begin{eqnarray} \label{pulemet}
\Vert
(\partial_{y_n} \mbox{\bf w}_{\nu}(\cdot,0),\mbox{\bf  w}
_{\nu}(\cdot,0))\Vert
_{H^{\frac{1}{4},\frac 12,\widetilde s}(\Bbb R^n)\times H^{\frac{3}{4},\frac 32,\widetilde s}(\Bbb R^n)} + \root\of{\vert\widetilde s\vert}
\Vert
\mbox{\bf w}_{\nu}\Vert_{H^{0,1,\widetilde s}(\mathcal Q)}\le
C_{13}(\widetilde\varphi^\frac {5}{12}( y^*_0)\Vert \text{\bf w}
\Vert_{H^{\frac 12,1,\widetilde s}(\mathcal Q)}                   \nonumber\\
+ \Vert
\mbox{\bf g}e^{\vert s\vert\varphi}\Vert_{H^{\frac{1}{4},\frac 12,\widetilde s}(\Bbb R^n)}
+\Vert \mbox{\bf P}(y,D,\widetilde s)\mbox{\bf
w}_{\nu}\Vert_{L^2(\mathcal Q)}\nonumber\\
+\root\of{\widetilde \varphi( y_0^*)}\Vert
(\partial_{y_n} \mbox{\bf w}(\cdot,0),\mbox{\bf  w}(\cdot,0))\Vert
_{H^{\frac{1}{4},\frac 12,\widetilde s}(\Bbb R^n)\times H^{\frac{3}{4},\frac 32,\widetilde s}(\Bbb R^n)}).
\end{eqnarray}

{\bf Subcase B.} Let  $r_{\lambda+2\mu}(y^*,\zeta^*)\ne 0.$
Then by Proposition \ref{gorokx} and Proposition \ref{gorokx1}
there exists a constant $C_{14}$ independent of $s$  such that
\begin{equation}\label{ix}
\Vert (\frac 1i \partial_{y_n} w_{2,\nu}
-\Gamma^+
_{\lambda+2\mu}(y,{\widetilde D},\widetilde s))w_{2,\nu}\vert_{y_n=0}\Vert
_{H^{\frac 14,\frac 12, \widetilde s}(\Bbb R^n)}
\end{equation}
$$
\le C_{14}(\Vert P_{\rho,\lambda+2\mu}(y,D,\widetilde s)w_{2,\nu}\Vert
_{L^2(\mathcal Q)}
+ \widetilde\varphi^\frac {5}{12}( y^*_0)\Vert w_{2}\Vert_{H^{\frac 12,1,\widetilde s}(\mathcal Q)}).
$$

On the other hand, on $\Bbb R^n$ from the boundary condition (\ref{!poko1}) and (\ref{legioner}) we have
\begin{equation}\label{ix1}
((\lambda +2\mu)(y^*)\left(\partial_{y_n} w_{2,\nu}
-\vert \widetilde s\vert{\varphi}_{n}(y^*)w_{2,\nu}\right)
- \mu(y^*)\sum_{j=1}^{n-1} \left(\partial_{y_j}
w_{j,n,\nu}-\vert\widetilde s\vert{\varphi}_{j}(y^*)w_{j,n,\nu}
\right))(\cdot,0)
= \mbox{\bf r},
\end{equation}
where the function $\mbox{\bf r}$ satisfies
\begin{eqnarray}\label{garmodon}
\Vert \mbox{\bf r}\Vert^2_{L^2(\Bbb R^n)}
\le \epsilon(\delta)\Vert
(\partial_{y_n} \mbox{\bf w}_{\nu}(\cdot,0),\mbox{\bf  w}
_{\nu}(\cdot,0))\Vert^2
_{L^2(\Bbb R^n)\times H^{0,1,\widetilde s}(\Bbb R^n)} \nonumber\\+ \frac{C_{15}
}{1+\vert \widetilde s\vert}\Vert
(\partial_{y_n} \mbox{\bf w}(\cdot,0),\mbox{\bf  w}(\cdot,0))\Vert^2
_{L^2(\Bbb R^n)\times H^{0,1,\widetilde s}(\Bbb R^n)}
+ C_{16}\Vert \mbox{\bf g}e^{\vert s\vert\varphi}\Vert^2_{L^2(\Bbb R^n)}
\end{eqnarray}
with some constants $C_{15}, C_{16}$ independent of $\widetilde s.$

From (\ref{ix1}), (\ref{garmodon}) and (\ref{2.4}), it follows that
\begin{eqnarray}\label{02.4}
\vert\widetilde s\vert\left\Vert(\lambda +2\mu)(y^*)
\left(\partial_{y_n} w_{2,\nu}-\vert\widetilde s\vert
{\varphi}_{n}(y^*)w_{2,\nu}\right)(\cdot,0)\right\Vert^2_{L^2(\Bbb R^n)} \\
\le C_{17}\sum_{j=1}^{n-1}\varXi_\mu(w_{j,n,\nu}) + \epsilon(\delta)
\vert\widetilde s\vert \Vert
(\partial_{y_n} \mbox{\bf w}_{\nu}(\cdot,0),\mbox{\bf  w}
_{\nu}(\cdot,0))\Vert^2
_{L^2(\Bbb R^n)\times H^{0,1,\widetilde s}(\Bbb R^n)}\nonumber\\
+C_{18}(\widetilde \varphi( y_0^*)\Vert
(\partial_{y_n} \mbox{\bf w}(\cdot,0),\mbox{\bf  w}
(\cdot,0))\Vert^2
_{L^2(\Bbb R^n)\times H^{0,1,\widetilde s}(\Bbb R^n)}
+ \vert\widetilde s\vert\Vert \mbox{\bf g}
e^{\vert s\vert\varphi}
\Vert_{L^2(\Bbb R^n)}^2).\nonumber
\end{eqnarray}
Then this estimate, the G\aa rding inequality (\ref{goblin}) and (\ref{ix})
imply
\begin{eqnarray}\label{robin}
\vert\widetilde s\vert^3\Vert w_{2,\nu}(\cdot,0)\Vert^2
_{L^2(\Bbb R^n)}
+ \vert\widetilde s\vert\sum_{j=1}^n\Vert \partial_{y_j} w_{2,\nu}
(\cdot,0)\Vert_{L^2(\Bbb R^n)}^2 \\
\le C_{19}\sum_{j=1}^{n-1}\varXi_\mu(w_{j,n,\nu})+\epsilon(\delta)\vert
\widetilde s\vert
\Vert (\partial_{y_n} \mbox{\bf w}_{\nu}(\cdot,0),\mbox{\bf  w}
_{\nu}(\cdot,0))\Vert^2
_{L^2(\Bbb R^n)\times H^{0,1,\widetilde s}(\Bbb R^n)}\nonumber\\
+ C_{20}(\vert\widetilde s\vert\Vert \mbox{\bf g}
e^{\vert s\vert\varphi}
\Vert_{L^2(\Bbb R^n)}^2+\widetilde \varphi( y_0^*)\Vert
(\partial_{y_n} \mbox{\bf w}(\cdot,0),\mbox{\bf  w}
(\cdot,0))\Vert^2
_{L^2(\Bbb R^n)\times H^{0,1,\widetilde s}(\Bbb R^n)}\nonumber\\
+\Vert P_{\rho,\lambda+2\mu}(y,D,\widetilde s)w_{2,\nu}\Vert^2
_{L^2(\mathcal Q)}
+\widetilde\varphi^\frac {5}{6}( y^*)\Vert w_{2}\Vert^2_{H^{\frac 12,1,\widetilde s}(\mathcal Q)}). \nonumber
\end{eqnarray}
Inequalities (\ref{2.4}), (\ref{robin}) and (\ref{02.4}) imply (\ref{pulemet}).

\bigskip
\section{ Case $r_{\lambda+2\mu}(y^*,\zeta^*)=0.$}\label{grom2}
\bigskip
Let $(y^*,\zeta^*)$ be a point on ${\Bbb R}^{n+1}\times \Bbb M$ such
that $r_{\lambda+2\mu}(y^*,\zeta^*)=0$ and $\text{supp}\,
\chi_\nu\subset \mathcal O(y^*,\delta_1(y^*)).$  Since the case $r_\mu(y^*,\zeta^*)=r_{\lambda+2\mu}(y^*,\zeta^*)=0$ was treated in the previous section one can  assume that
\begin{equation}
r_\mu(y^*,\zeta^*)\ne 0.
\end{equation}
By (\ref{normal}) and (\ref{3.37})
the  equality $r_{\lambda+2\mu}(y^*,\zeta^*)=0$ implies that
\begin{equation}\label{medved4}\sum_{j=1}^{n-1}({\xi_j^*})^2=(\widetilde s^*)^2\sum_{j=1}^{n-1} \vert\partial_{y_j}{\varphi}(y^*)\vert^2\quad
\mbox{and}\quad \frac{\rho(y^*)\xi_0^*}{(\lambda+2\mu)(y^*)}-\widetilde s^*\sum_{j=1}^{n-1}\xi^*_j\partial_{y_j}{\varphi}(y^*)=0.\end{equation}

Form (\ref{medved4}) and (\ref{pinok220}) we immediately obtain $s^*\ne 0.$

By (\ref{medved4}) there exists $C_1>0$ such that for all $(y^*,\zeta)$ from $\mathcal O(y^*,\delta_1(y^*))$ we have
\begin{equation}\label{2.7y}\left  \vert \sum_{j=1}^{n-1}{\xi^2_j}- \widetilde s^2\sum_{j=1}^{n-1} \vert\partial_{y_j}{\varphi}(y^*)\vert^2\right \vert+
\left \vert \frac{\rho(y^*)\xi_0}{(\lambda+2\mu)(y^*)}+\widetilde s\sum_{j=1}^{n-1}\xi_j\partial_{y_j}{\varphi}(y^*)\right\vert \le \delta_1C_1M(\widetilde \xi,\widetilde s).
\end{equation}

By (\ref{klop}) there exists $C_{2}>0$ independent of $s$ such that
\begin{equation}\label{3.19}
\varXi_{\lambda+2\mu}(w_{2,\nu})+C_{2}({\vert \widetilde s\vert}\Vert
w_{2,\nu}\Vert^2_{H^{0,1}(\mathcal Q)}+ {\vert \widetilde s\vert}^3\Vert
w_{2,\nu}\Vert^2_{L^2(\mathcal Q)})
\end{equation}
\begin{eqnarray*}
\le &&C_{3}\Vert P_{\rho,\lambda+2\mu}(y,D,\widetilde s) w_{2,\nu} \Vert^2
_{L^2(\mathcal Q)}
+ \epsilon  \Vert
(\partial_{y_n} \mbox{\bf w}_{\nu}(\cdot,0),\mbox{\bf  w}
_{\nu}(\cdot,0))\Vert^2
_{H^{\frac{1}{4},\frac 12,\widetilde s}(\Bbb R^n)\times H^{\frac{3}{4},\frac 32,\widetilde s}(\Bbb R^n)}
\end{eqnarray*}
where $\epsilon(\delta)\rightarrow 0$ as $\delta\rightarrow +0.$

By (\ref{01}), (\ref{02}), (\ref{2.7y}) and (\ref{3.18}), we have
\begin{eqnarray}\label{3.21}
\vert \frak J_2(\lambda+2\mu,w_{2,\nu})+\frak J_3(\lambda+2\mu,w_{2,\nu})\vert
\le C_{4}\delta_1 \Vert (\partial_{y_n}
w_{2,\nu}(\cdot,0),w_{2,\nu}(\cdot,0))\Vert^2
_{H^{\frac{1}{4},\frac 12,\widetilde s}(\Bbb R^n)\times H^{\frac{3}{4},\frac 32,\widetilde s}(\Bbb R^n)}\nonumber\\
+C_{5}\widetilde \varphi( y_0^*)\Vert (\partial_{y_n}
w_{2}(\cdot,0),w_{2}(\cdot,0))\Vert^2
_{L^2(\Bbb R^n)\times H^{0,1,\widetilde s}(\Bbb R^n)}.
\end{eqnarray}

By (\ref{3.21}), there exists a constant $C_{6}>0$ such that
\begin{eqnarray}\label{3.22}
\varXi_{\lambda+2\mu}(w_{2,\nu})\ge C_{6}\int_{{\Bbb R^n}} \left(
\vert\widetilde s\vert
\vert\partial_{y_n} w_{2,\nu}\vert^2
+ {\vert \widetilde s\vert}^3 \vert
w_{2,\nu}\vert^2\right)(\widetilde y,0)d \widetilde y        \nonumber\\
-\epsilon
\Vert (\partial_{y_n}
w_{2,\nu}(\cdot,0),w_{2,\nu}(\cdot,0))\Vert^2
_{H^{\frac{1}{4},\frac 12,\widetilde s}(\Bbb R^n)\times H^{\frac{3}{4},\frac 32,\widetilde s}(\Bbb R^n)}
-C_{7}\widetilde \varphi( y_0^*)\Vert (\partial_{y_n}
w_{2}(\cdot,0),w_{2}(\cdot,0))\Vert^2
_{L^2(\Bbb R^n)\times H^{0,1,\widetilde s}(\Bbb R^n)},\nonumber
\end{eqnarray} where $\epsilon(\delta)\rightarrow 0$ as $\delta\rightarrow +0.$
Since $\widetilde s^*\ne 0$, we have
\begin{eqnarray}\label{3.22}
\varXi_{\lambda+2\mu}(w_{2,\nu})\ge C_{8}\int_{{\Bbb R^n}} (
{\vert \widetilde s\vert}\sum_{j=1}^{n-1} \vert\partial_{y_j}
w_{2,\nu}\vert^2+{\vert \widetilde s\vert}^3
\vert w_{2,\nu}\vert^2)(\widetilde y,0)d \widetilde y\\
- \epsilon\Vert (\partial_{y_n}
w_{2,\nu}(\cdot,0),w_{2,\nu}(\cdot,0))\Vert^2
_{H^{\frac{1}{4},\frac 12,\widetilde s}(\Bbb R^n)\times H^{\frac{3}{4},\frac 32,\widetilde s}(\Bbb R^n)}\nonumber\\-C_{9}\widetilde \varphi(y_0^*)\Vert
(\partial_{y_n}
w_{2}(\cdot,0),w_{2}(\cdot,0))\Vert^2
_{L^2(\Bbb R^n)\times H^{0,1,\widetilde s}(\Bbb R^n)}.\nonumber
\end{eqnarray}
Since $r_{\mu}(y^*,\zeta^*)\ne 0$, then by Proposition \ref{gorokx} and Proposition \ref{gorokx1}
there exists a constant $C_{10}$ independent of $s$  such that
\begin{equation}\label{ix0}
\Vert (\frac 1i \partial_{y_n} \mbox{\bf w}_{1,\nu}
-\Gamma^+
_{\mu}(y,{\widetilde D},\widetilde s)\mbox{\bf w}_{1,\nu})\vert_{y_n=0}\Vert
_{H^{\frac 14,\frac 12, \widetilde s}(\Bbb R^n)}
\end{equation}
$$
\le C_{10}(\Vert P_{\rho,\mu}(y,D,\widetilde s)\mbox{\bf w}_{1,\nu}\Vert
_{L^2(\mathcal Q)}
+ \widetilde\varphi^\frac {5}{12}( y^*)\Vert \mbox{\bf w}_{1}\Vert_{H^{\frac 12,1,\widetilde s}(\mathcal Q)}).
$$

By (\ref{!poko1}), (\ref{legioner}) for any $k\in\{1,\dots, n-1\}$, we obtain
\begin{equation}\label{barmaleil}
\mu(y^*)\left(\partial_{y_n} w_{ k,n,\nu}-{\vert \widetilde s
\vert} \varphi_{n}(y^*) w_{k,n,\nu}\right)(\widetilde y,0)
\end{equation}
$$
= (\lambda+\mu)(y^*)\left(
\partial_{y_k} w_{2,\nu}-{\vert \widetilde s\vert}
\varphi_{{k}}(y^*) w_{2,\nu}\right)(\widetilde y,0)
+ \mbox{\bf r}(\widetilde y)\quad \mbox{in}\,\,\Bbb R^n,
$$
where the function $\mbox{\bf r}$ satisfies estimate (\ref{garmodon}).
From  (\ref{barmaleil}), (\ref{garmodon}) and (\ref{ix0})  we obtain the
estimate
\begin{eqnarray}\label{3.25}
{\vert \widetilde s\vert}\Vert \partial_{y_n} \mbox{\bf w}_{1,\nu}(\cdot,0)-\vert \widetilde s\vert{\varphi}_{n}(y^*)\mbox{\bf w}_{1,\nu}(\cdot,0)\Vert^2
_{L^2({\Bbb R^n})}\nonumber\\
\le C_{11}\biggl(\int_{{\Bbb R^n}} (
{\vert \widetilde s\vert}\sum_{j=1}^{n-1}\vert\partial_{y_j} w_{2,\nu}(\widetilde y,0)\vert^2
+{\vert \widetilde s\vert}^3\vert w_{2,\nu}(\widetilde y,0)\vert^2)d
\widetilde y\nonumber\\
+\Vert P_{\rho,\mu}(y,D,\widetilde s)
\mbox{\bf w}_{1,\nu} \Vert^2_{L^2(\mathcal Q)}+\widetilde\varphi^\frac {5}{6}( y^*_0)\Vert
\text{\bf w}\Vert^2_{H^{\frac 12,1,\widetilde s}(\mathcal Q)}\\
+ \epsilon(\delta)
\vert\widetilde s\vert \Vert
(\partial_{y_n} \mbox{\bf w}_{\nu}(\cdot,0),\mbox{\bf  w}
_{\nu}(\cdot,0))\Vert^2
_{L^2(\Bbb R^n)\times H^{0,1,\widetilde s}(\Bbb R^n)}\nonumber\\
+\widetilde \varphi( y_0^*)\Vert
(\partial_{y_n} \mbox{\bf w}(\cdot,0),\mbox{\bf  w}
(\cdot,0))\Vert^2
_{L^2(\Bbb R^n)\times H^{0,1,\widetilde s}(\Bbb R^n)}
+ \vert\widetilde s\vert\Vert \mbox{\bf g}
e^{\vert s\vert\varphi}
\Vert_{L^2(\Bbb R^n)}^2\biggr).\nonumber
\end{eqnarray}

From the G\aa rding inequality and (\ref{3.25}), (\ref{ix0})
\begin{eqnarray}\label{lika}
\vert \widetilde s\vert \Vert \mbox{\bf w}_{1,\nu}(\cdot,0)\Vert^2_{H^{1,\widetilde s}(\Bbb R^n)}\le C_{12}\biggl(\Vert P_{\rho,\mu}(y,D,\widetilde s)\mbox{\bf w}_{1,\nu}\Vert^2
_{L^2(\mathcal Q)}
+\widetilde\varphi^\frac {5}{6}( y^*_0) \Vert \mbox{\bf w}_{1,\nu}\Vert^2_{H^{\frac 12,1,\widetilde s}(\mathcal Q)})\nonumber\\
+\int_{{\Bbb R^n}} (
{\vert \widetilde s\vert}\sum_{j=1}^{n-1}\vert\partial_{y_j} w_{2,\nu}(\widetilde y,0)\vert^2
+{\vert \widetilde s\vert}^3\vert w_{2,\nu}(\widetilde y,0)\vert^2)d
\widetilde y\nonumber\\
+ \epsilon(\delta)
\vert\widetilde s\vert \Vert
(\partial_{y_n} \mbox{\bf w}_{\nu}(\cdot,0),\mbox{\bf  w}
_{\nu}(\cdot,0))\Vert^2
_{L^2(\Bbb R^n)\times H^{0,1,\widetilde s}(\Bbb R^n)}\nonumber\\
+\widetilde \varphi( y_0^*)\Vert
(\partial_{y_n} \mbox{\bf w}(\cdot,0),\mbox{\bf  w}
(\cdot,0))\Vert^2
_{L^2(\Bbb R^n)\times H^{0,1,\widetilde s}(\Bbb R^n)}
+ \vert\widetilde s\vert\Vert \mbox{\bf g}
e^{\vert s\vert\varphi}
\Vert_{L^2(\Bbb R^n)}^2\biggl ).\nonumber
\end{eqnarray}

Inequalities (\ref{3.22}), (\ref{3.25}), (\ref{lika}) imply
\begin{eqnarray}\label{3.26}
\varXi_{\lambda+2\mu}(w_{2,\nu})\ge C_{13}{\vert \widetilde s\vert}
\Vert
(\partial_{y_n} \mbox{\bf w}_{\nu}(\cdot,0),\mbox{\bf  w}
_{\nu}(\cdot,0))\Vert^2
_{L^2(\Bbb R^n)\times H^{0,1,\widetilde s}(\Bbb R^n)}\nonumber\\
- C_{14}(\Vert
\mbox{\bf  P}(y,D,\widetilde s) \mbox{\bf w}_{\nu}\Vert^2
_{L^2(\mathcal Q)}+{\vert \widetilde s\vert}\Vert  \mbox{\bf g}
e^{\vert s\vert\varphi}\Vert^2_{L^2({\Bbb R^n})}+\widetilde\varphi^\frac {5}{6}( y^*_0)\Vert \text{\bf w}
\Vert^2_{H^{\frac 12,1,\widetilde s}(\mathcal Q)})\nonumber\\
-C_{15}\widetilde \varphi( y_0^*)\Vert (\partial_{y_n}
\text{\bf w}(\cdot,0),\text{\bf w}(\cdot,0))\Vert^2
_{L^2(\Bbb R^n)\times H^{0,1,\widetilde s}(\Bbb R^n)},
\end{eqnarray}
where $C_{13}>0.$  From (\ref{3.19}) and (\ref{3.26}), we obtain
(\ref{pulemet}).

\bigskip
\section{ Case $r_\mu(y^*,\zeta^*)\ne 0$ and
$r_{\lambda+2\mu}(y^*,\zeta^*)\ne 0.$}\label{QQ5}
\bigskip

In this section we consider the conic neighborhood $ \mathcal O(y^*,
\delta_1(y^*))$ of the point $(y^*,\zeta^*)$ such that
\begin{equation}\label{4.1}
r_\mu(y^*,\zeta^*)\ne 0\quad\mbox{
and }
r_{\lambda+2\mu}(y^*,\zeta^*)\ne 0.\end{equation}

In that case, thanks to (\ref{4.1}) and Proposition \ref{gorokx1},
factorization
(\ref{min}) holds true for $\beta=\mu$ and $\beta=\lambda+2\mu.$
Then Proposition \ref{gorokx} yields the a priori estimate
\begin{eqnarray}\label{4.4}
\sum_{k,j=1, k<j}^n\Vert
 V^+_{\mu}(k,j)(\cdot,0)\Vert_{H^{\frac 14,\frac 12,\widetilde s}({\Bbb R^n})}+
\Vert
V^+_{\lambda+2\mu}(\cdot,0)\Vert_{H^{\frac 14,\frac 12,\widetilde s}({\Bbb R^n})}   \\
\le C_1(\Vert \mbox{\bf P}(y,D,\widetilde s)\mbox{\bf w}_{\nu}\Vert
_{L^2(\mathcal Q)} +\widetilde\varphi^\frac {5}{12}( y^*)\Vert
\text{\bf w}\Vert_{H^{0,1,\widetilde s}(\mathcal Q)}).   \nonumber
\end{eqnarray}

Using (\ref{ix}), (\ref{ix0}),  we rewrite (\ref{ix1}) and  (\ref{barmaleil})  as
\begin{equation}\label{4.5}
(\frac{\lambda+2\mu}{\mu}(y^*)\left(\partial_{y_j}  w_{2,\nu}
-\vert\widetilde  s\vert{\varphi}_{j}
w_{2,\nu}\right)
- {i}\alpha^+_\mu(\widetilde y,0,\widetilde D,\widetilde s)w_{j,n,\nu})(\cdot,0)
=V^+_{\mu}(i,n)(\cdot,0)-r_{j,n,\nu}\quad \mbox{in}\,\,\Bbb R^n,
\end{equation}
where $k\in\{1,\dots, n-1\}$ and
\begin{eqnarray}\label{4.6}
(\sum_{k=1}^{n-1}\frac{\mu}{\lambda+2\mu}(y^*)\left(-\partial_{y_k}
w_{k,n,\nu}+{\vert \widetilde s\vert}{\varphi}_{k}
w_{k,n,\nu}\right)\\
- {i}\alpha^+_{\lambda+2\mu}(\widetilde y,0,\widetilde D,\widetilde s)
w_{2,\nu})(\cdot,0)=V^+_{\lambda+2\mu}(\cdot,0)-r_{2,\nu}, \nonumber
\end{eqnarray}
where the function $\mbox{\bf r}=(r_{1,n,\nu},\dots, r_{n-1,n,\nu},r_{2,\nu})$
satisfies estimate (\ref{garmodon}).
 Let
$\text{\bf B}(\widetilde y,{\widetilde D},\widetilde s)$ be the matrix
pseudodifferential operator with the symbol
\begin{equation}\label{logo}
\mbox{\bf B}(\widetilde y,\widetilde\xi,\widetilde s)=
\left( \begin{matrix}
-i\alpha^+_\mu(\widetilde y,0,\widetilde\xi,\widetilde s)&0 &\dots&\frac{\lambda+2\mu}{\mu}
(i\xi_1-\vert\widetilde
s\vert{\varphi}_{1}) \\
0&\dots&\dots &\dots\\
0&
-i\alpha^+_\mu(\widetilde y,0,\widetilde\xi,\widetilde s)&\dots &\frac{\lambda+2\mu}{\mu}
(i\xi_j-\vert\widetilde
s\vert{\varphi}_{j}) \\
\dots&\dots&\dots &\dots\\
 \frac{\mu}{\lambda+2\mu}(-i\xi_1+\vert\widetilde s\vert{\varphi}_{1}) &\dots
&\frac{\mu}{\lambda+2\mu}(-i\xi_{n-1}+\vert\widetilde s\vert{\varphi}_{{n-1}})
 &-i\alpha^+_{\lambda+2\mu}(\widetilde y,0,\widetilde\xi,\widetilde s)
\end{matrix}\right).
\end{equation}
%
%\begin{equation}\label{logo}
%\mbox{\bf B}(\widetilde y,\zeta)=
%\left( \begin{matrix}
%-{i}\alpha^+_\mu(\widetilde y,0,\zeta)&0 &\dots&
%\frac{\lambda+2\mu}{\mu}({i}\xi_1-\vert\widetilde
%s\vert{\varphi}_{y_1}) \\
%0&\dots&\dots &\dots\\
%0& -{i}\alpha^+_\mu(\widetilde y,0,\zeta)&\dots &\frac{\lambda+2\mu}
%{\mu}({i}\xi_i-\vert\widetilde s\vert{\varphi}_{y_i}) \\
%\dots&\dots&\dots &\dots\\
%\frac{\mu}{\lambda+2\mu}(-{i}\xi_1+\vert\widetilde s\vert{\varphi}_{y_1})
%&\dots
%&\frac{\mu}{\lambda+2\mu}(-{i}\xi_{n-1}+\vert\widetilde s\vert
%{\varphi}_{y_{n-1}})
%&-{i}\alpha^+_{\lambda+2\mu}(\widetilde y,0,\zeta)
%\end{matrix}\right).
%\end{equation}

We have
\begin{proposition} Let $\zeta=(\xi_0,\xi_1+i\vert\widetilde s\vert{\varphi}_{1},\dots, \xi_{n-1}+i\vert\widetilde s\vert{\varphi}_{n-1}).$
The following formula is true:
\begin{eqnarray}\label{pardon}
det\,\text{\bf B}(y^*,\widetilde \xi,\widetilde s)=(-{i})^n (\alpha^+_\mu(y^*,\widetilde\xi,\widetilde s))^{n-1}
\alpha^+_{\lambda+2\mu}(y^*,\widetilde\xi,\widetilde s)                           \nonumber\\
+ (-1)^{n-1}(-{i})^{n-2} (\alpha^+_\mu(y^*,\widetilde\xi,\widetilde s))^{n-2}\sum_{j=1}
^{n-1}
(-{i}\xi_j+\vert\widetilde s\vert{\varphi}_{j}(y^*))^2.
\end{eqnarray}
\end{proposition}
{\bf Proof.} By $\text{\bf B}_{n}(y,\widetilde\xi,\widetilde s)$ we denote the matrix determined by
(\ref{logo}) of the size $n\times n$ and $\text{\bf B}_{k,j,n}(y,\widetilde\xi,\widetilde s)$ be the minor
obtained
from the matrix  $\text{\bf B}_{n}(y,\widetilde\xi,\widetilde s)$ by crossing out the $k$-th row and
the $j$-th column.   Our proof is based on the  induction method.
Except the formula (\ref{pardon}), we claim
\begin{equation}\label{paka}
\vert \text{\bf B}_{1,n-1,n}(y^*,\widetilde\xi,\widetilde s)\vert
= (-{i})^{n} (\alpha^+_\mu(y^*,\widetilde\xi,\widetilde s))^{n-2}(-{i}\xi_1
+\vert \widetilde s\vert{\varphi}_{1}(y^*))^2.
\end{equation}
For $n=2,3$, we can easily verify the formulae by direct computations.
Suppose that (\ref{pardon}) and (\ref{paka}) are true for $n-1$. Then
\begin{eqnarray}\label{pardon3}
det\,\text{\bf B}_{n-1}(y^*,\widetilde\xi,\widetilde s )=(-{i})^{n-1}
(\alpha^+_\mu(y^*,\widetilde\xi,\widetilde s))^{n-2}\alpha^+_{\lambda+2\mu}(y^*,\widetilde\xi,\widetilde s)
                                         \nonumber\\
+ (-1)^{n-1}(-{i})^{n-3} (\alpha^+_\mu(y^*,\widetilde\xi,\widetilde s))^{n-3}
\sum_{j=1}^{n-2}(-{i}\xi_j+\vert \widetilde s\vert{\varphi}_{j}(y^*))^2
\end{eqnarray}
and
\begin{equation}\label{pardon1}
\vert \text{\bf B}_{1,n-1,n}(y^*,\widetilde\xi,\widetilde s)\vert
= (-{i})^{n-1}\frac{\mu}{\lambda+2\mu}(y^*)(\alpha^+_
\mu(y^*,\widetilde\xi,\widetilde s))^{n-2}
(-{i}\xi_1+\vert \widetilde s\vert{\varphi}_{1}(y^*))^2.
\end{equation}

Since $det\,\text{\bf B}_{n}(y^*,\widetilde\xi,\widetilde s) = -{i}\alpha^+_
\mu(y^*,\widetilde\xi,\widetilde s)
\vert \text{\bf B}_{1,1,n}(y^*,\widetilde\xi,\widetilde s)\vert$ \newline$+ (-1)^{1+n}\frac{\lambda+2\mu}{\mu}(y^*)
({i}\xi_1-\vert\widetilde s\vert{\varphi}_{1}(y^*))\vert \text{\bf B}
_{1,n,n}(y^*,\widetilde\xi,\widetilde s)\vert$,
by (\ref{pardon}) and (\ref{pardon1}) we have
\begin{eqnarray}
\mbox{det}\,\text{\bf B}_{n}(y^*,\widetilde\xi,\widetilde s)
= -{i}\alpha^+_\mu( y^*,\widetilde\xi,\widetilde s)((-{i})^{n-1}
(\alpha^+_\mu(y^*,\widetilde\xi,\widetilde s))^{n-2}\alpha^+_{\lambda+2\mu}(y^*,\widetilde\xi,\widetilde s)
                                                       \nonumber\\
+ (-1)^{n-1}(-{i})^{n-3} (\alpha^+_\mu(y^*,\widetilde\xi,\widetilde s))^{n-3}
\sum_{j=2}^{n-1}
(-{i}\xi_j+\vert \widetilde s\vert{\varphi}_{j}(y^*))^2            \nonumber\\
+(-1)^{1+n}\frac{\lambda+2\mu}{\mu}(y^*)(-{i}\alpha^+_
\mu(y^*,\widetilde\xi,\widetilde s))^{n-2}
({i}\xi_1-\vert\widetilde
s\vert{\varphi}_{1}(y^*))\frac{\mu}{\lambda+2\mu}(y^*)
(-{i}\xi_1+\vert\widetilde s\vert{\varphi}_{1}(y^*))       \nonumber\\
= (-{i})^{n} (\alpha^+_\mu(y^*,\widetilde\xi,\widetilde s))^{n-1}\alpha^+_{\lambda+2\mu}
(y^*,\widetilde\xi,\widetilde s)                           \nonumber\\
+ (-{i})^{n-2}(-1)^{n-1} (\alpha^+_\mu(y^*,\widetilde\xi,\widetilde s))^{n-2}\sum_{j=2}^{n-1}
(-{i}\xi_j+\vert \widetilde s\vert{\varphi}_{j}(y^*))^2)\nonumber\\
+ (-1)^{1+n}(-{i}\alpha^+_\mu(y^*,\widetilde\xi,\widetilde s))^{n-2}(-{i}\xi_1
+\vert\widetilde s\vert{\varphi}_{1}(y^*))^2.\nonumber
\end{eqnarray}
The proof of the proposition is complete. $\blacksquare$
\\

If $(\widetilde \xi^*,\widetilde s^*)\ne \{  (\widetilde \xi,\widetilde s)\in \Bbb M\vert det\,\text{\bf B}(y^*,\widetilde\xi,\widetilde s)=0\}$
by  (\ref{4.5}), (\ref{4.6}) we have
\begin{eqnarray} \label{ziza}
\Vert
(\partial_{y_n} \mbox{\bf w}_{\nu}(\cdot,0),\mbox{\bf  w}
_{\nu}(\cdot,0))\Vert
_{H^{\frac{1}{4},\frac 12,\widetilde s}(\Bbb R^n)\times H^{\frac{3}{4},\frac 32,\widetilde s}(\Bbb R^n)} \le
C_{2}(\widetilde\varphi^\frac {5}{12}( y^*_0)\Vert \text{\bf w}
\Vert_{H^{0,1,\widetilde s}(\mathcal Q)}
+ \Vert
\mbox{\bf g}e^{\vert s\vert\varphi}\Vert_{H^{\frac{1}{4},\frac 12,\widetilde s}(\Bbb R^n)}\nonumber\\
+\Vert \mbox{\bf P}(y,D,\widetilde s)\mbox{\bf
w}_{\nu}\Vert_{L^2(\mathcal Q)}
+\frac{C_{3}}{\root\of{1+\vert\widetilde s\vert}}\Vert
(\partial_{y_n} \mbox{\bf w}(\cdot,0),\mbox{\bf  w}(\cdot,0))\Vert
_{H^{\frac{1}{4},\frac 12,\widetilde s}(\Bbb R^n)\times H^{\frac{3}{4},\frac 32,\widetilde s}(\Bbb R^n)}).
\end{eqnarray}
From this inequality and Proposition \ref{opana}  we obtain (\ref{pulemet}).

By (\ref{3.35}) -(\ref{3.37}) if $det\,\text{\bf B}(y^*,\widetilde \xi^*,\widetilde s^*)=0$ and
(\ref{4.1}) holds true, then
\begin{equation}\label{4.7}
(\widetilde \xi^*,\widetilde s^*)\in \mathcal U=\left \{(\widetilde \xi,\widetilde s)\in {\Bbb R}^{n}; \thinspace
\sum_{j=1}^{n-1} (\xi_j+{i}\vert\widetilde  s\vert{\varphi}_{j}
(y^*))^2=-\frac{\rho(y^*)i\xi_0}
{(\lambda+3\mu)(y^*)}\right\}.\nonumber
\end{equation}
If $(\widetilde \xi,\widetilde s)\in \mathcal U$ then
\begin{equation}\label{ono}
\sum_{j=1}^{n-1}\xi_j^2=\widetilde s^2\sum_{j=1}^{n-1}\varphi_j^2(y^*).
\end{equation}
By (\ref{3.35}), (\ref{3.36}) and (\ref{napoleon})
$$
\Gamma^\pm_\beta(y^*,\widetilde \xi^*,\widetilde s^*)=-i\vert \widetilde s^*\vert{\varphi}_{n}(y^*)\pm e^{i\frac {\pi }{4}}\root\of{\frac{(\lambda+3\mu-\beta)(y^*)}{\beta(y^*)}}\root\of{\vert\widetilde s^*\vert\vert\sum_{j=1}^{n-1}\xi_j^*{\varphi}_j(y^*)\vert}
$$
So
\begin{eqnarray}\label{oko}
\left \vert\mbox{Im}\left\{\,e^{i\frac {\pi }{4}}\root\of{\frac{(\lambda+3\mu-\beta)(y^*)}{\beta(y^*)}}\root\of{\vert\widetilde s^*\vert\vert\sum_{j=1}^{n-1}\xi_j^*{\varphi}_j(y^*)\vert}\right\}\right\vert=\frac{1}{\root\of{2}}\root\of{\frac{(\lambda+3\mu-\beta)(y^*)}{\beta(y^*)}}\root\of{\vert\widetilde s^*\vert\vert\sum_{j=1}^{n-1}\xi_j^*{\varphi}_j(y^*)\vert}\nonumber\\
\le\root\of{ \frac{\vert \widetilde s^*\vert}{2}}\root\of{\frac{\mu(y^*)}{(\lambda+2\mu)(y^*)}}  \root\of{\sum_{j=1}^{n-1}\varphi_j^2(y^*)}.
\end{eqnarray}
Here in order to get the  last inequality in (\ref{oko}) we used (\ref{ono}).
By (\ref{pinok22})
\begin{equation}
-\mbox{Im}\,\Gamma^\pm_{\lambda+2\mu}(y^*,\widetilde \xi^*,\widetilde s^*)>0.
\end{equation}
Then by Proposition \ref{gorokx}
there exists a constant $C_{4}$ independent of $\widetilde s$  such that
\begin{equation}\label{ixq}
\Vert (\frac 1i \partial_{y_n} w_{2,\nu}
-\Gamma^\pm
_{\lambda+2\mu}(y,{\widetilde D},\widetilde s))w_{2,\nu}\vert_{y_n=0}\Vert
_{H^{\frac 14,\frac 12, \widetilde s}(\Bbb R^n)}
\end{equation}
$$
\le C_{4}(\Vert P_{\rho,\lambda+2\mu}(y,D,\widetilde s)w_{2,\nu}\Vert
_{L^2(\mathcal Q)}
+ \widetilde\varphi^\frac {5}{12}( y^*_0)\Vert w_{2}\Vert_{H^{\frac 12,1,\widetilde s}(\mathcal Q)}).
$$
Inequalities (\ref{ixq})  imply
\begin{equation}\label{zov}
\Vert \alpha
_{\lambda+2\mu}(y,{\widetilde D},\widetilde s)w_{2,\nu}\vert_{y_n=0}\Vert
_{H^{\frac 14,\frac 12, \widetilde s}(\Bbb R^n)}
\le C_{5}(\Vert P_{\rho,\lambda+2\mu}(y,D,\widetilde s)w_{2,\nu}\Vert
_{L^2(\mathcal Q)}+ \widetilde\varphi^\frac {5}{12}( y^*_0)\Vert w_{2}\Vert_{H^{\frac 12,1,\widetilde s}(\mathcal Q)}
).
\end{equation}
By  G\aa rding's  inequality we obtain from (\ref{zov})
\begin{equation}\label{ixqq}
\Vert w_{2,\nu}(\cdot,0)\Vert
_{H^{\frac 34,\frac 32, \widetilde s}(\Bbb R^n)}
\le C_{6}(\Vert P_{\rho,\lambda+2\mu}(y,D,\widetilde s)w_{2,\nu}\Vert
_{L^2(\mathcal Q)}
+ \widetilde\varphi^\frac {5}{12}( y^*_0)\Vert w_{2}\Vert_{H^{0,1,\widetilde s}(\mathcal Q)}).
\end{equation}
From (\ref{ixq}) and (\ref{ixqq}) we obtain
\begin{equation}\label{ixqq1}
\Vert (\partial_{y_n}w_{2,\nu},w_{2,\nu})(\cdot,0)\Vert
_{H^{\frac 14,\frac 12, \widetilde s}(\Bbb R^n)\times H^{\frac 34,\frac 32, \widetilde s}(\Bbb R^n)}\le C_{7}(\Vert P_{\rho,\lambda+2\mu}(y,D,\widetilde s)w_{2,\nu}\Vert
_{L^2(\mathcal Q)}
\end{equation}
$$
+\widetilde\varphi^\frac {5}{12}( y^*_0)\Vert w_{2}\Vert_{H^{\frac 12,1,\widetilde s}(\mathcal Q)}
).
$$
On the other hand the (\ref{ix}) and (\ref{barmaleil}) holds true. By (\ref{ixqq1}) we obtain from (\ref{barmaleil})
\begin{equation}\label{barmaleilk}
\mu(y^*)\left(\partial_{y_n} w_{ k,n,\nu}-{\vert \widetilde s
\vert} \varphi_{n}(y^*) w_{k,n,\nu}\right)(\cdot,0)=\widetilde{\mbox{\bf r}}\quad \mbox{in}\,\,\Bbb R^n,
\end{equation}
where function $\widetilde{\mbox{\bf r}}$ satisfies the estimate
\begin{eqnarray}\label{garmodon1}
\Vert\widetilde {\mbox{\bf r}}\Vert^2_{L^2(\Bbb R^n)}
\le \epsilon(\delta)\Vert
(\partial_{y_n} \mbox{\bf w}_{\nu}(\cdot,0),\mbox{\bf  w}
_{\nu}(\cdot,0))\Vert^2
_{L^2(\Bbb R^n)\times H^{0,1,\widetilde s}(\Bbb R^n)} \nonumber\\+ \frac{C_{8}
}{1+\vert \widetilde s\vert}(\Vert
(\partial_{y_n} \mbox{\bf w}(\cdot,0),\mbox{\bf  w}(\cdot,0))\Vert^2
_{L^2(\Bbb R^n)\times H^{0,1,\widetilde s}(\Bbb R^n)}
+ \nonumber\\
\Vert\mbox{\bf  P}(y,D,\widetilde s)\mbox{\bf  w}_{\nu}\Vert^2
_{L^2(\mathcal
Q)})+C_{9}\Vert \mbox{\bf g}e^{\vert s\vert\varphi}\Vert^2_{L^2(\Bbb R^n)}.
\end{eqnarray}
By (\ref{garmodon1}), (\ref{ixqq}) and (\ref{4.4})
\begin{equation}\label{lllxxx1}
\Vert \alpha_\mu(\widetilde y,0,\widetilde D,\widetilde s)\mbox{\bf w}_{1,\nu}(\cdot,0)\Vert
_{ H^{\frac 34,\frac 32, \widetilde s}(\Bbb R^n)}\le C_{10}(\Vert \mbox{\bf P}(y,D,\widetilde s)\mbox{\bf w}_{\nu}\Vert
_{L^2(\mathcal Q)}
\end{equation}
$$
+ \widetilde\varphi^\frac {5}{12}( y^*_0)\Vert w_{2}\Vert_{H^{\frac 12,1,\widetilde s}(\mathcal Q)}
+ \Vert
(\partial_{y_n} \mbox{\bf w}(\cdot,0),\mbox{\bf  w}(\cdot,0))\Vert
_{L^2(\Bbb R^n)\times H^{0,1,\widetilde s}(\Bbb R^n)}+\Vert \mbox{\bf g}e^{\vert s\vert\varphi}\Vert_{H^{\frac 14,\frac 12,\widetilde s}(\Bbb R^n)}).
$$
By  G\aa rding's inequality from (\ref{lllxxx1}) we have
\begin{equation}\label{lllxxx}
\Vert \mbox{\bf w}_{1,\nu}(\cdot,0)\Vert
_{H^{\frac 34,\frac 32, \widetilde s}(\Bbb R^n)}\le C_{11}(\Vert\mbox{\bf  P}(y,D,\widetilde s)\mbox{\bf  w}_{\nu}\Vert
_{L^2(\mathcal
Q)}
\end{equation}
$$
+ \widetilde\varphi^\frac {5}{12}( y^*_0)\Vert w_{2}\Vert_{H^{\frac 12,1,\widetilde s}(\mathcal Q)}
+ \Vert
(\partial_{y_n} \mbox{\bf w}(\cdot,0),\mbox{\bf  w}(\cdot,0))\Vert
_{L^2(\Bbb R^n)\times H^{0,1,\widetilde s}(\Bbb R^n)}+\Vert \mbox{\bf g}e^{\vert s\vert\varphi}\Vert_{H^{\frac 14,\frac 12,\widetilde s}(\Bbb R^n)}).
$$
Form (\ref{lllxxx}) and (\ref{ixqq}) we have (\ref{ziza}). From (\ref{ziza}) and Proposition \ref{opana}  we obtain (\ref{pulemet}).

\section{\bf End of the proof.}\label{Q2}
First we finish the proof of the Proposition \ref{zoopa}.
Now let us take the covering of the surface $\Bbb M$ of the sphere by
conical neighborhoods $\mathcal O(\zeta^*,\delta_1(\zeta^*)).$
From this covering we take the finite subcovering $\cup_{\nu=1}^N\mathcal O
(\zeta_\nu^*,\delta_1(\zeta_\nu^*))$.  Let $\chi_\nu$ be the partition of
unity associated to this subcovering. Hence $\sum_{\nu=1}^N\chi_\nu
(\widetilde \xi,\widetilde s)\equiv 1$ or all $(\widetilde \xi,\widetilde s)$ such that $M(\widetilde \xi,\widetilde s)\ge 1.$ Let $\chi_0(\widetilde \xi,\widetilde s)\in C^\infty_0(\Bbb R^{n+1})$ be a nonegative function  which is identically equal one  if $M(\widetilde \xi,\widetilde s)\le 1.$ Then by (\ref{pulemet})
\begin{eqnarray}\label{golifax}
\Vert(
\partial_{y_n} \mbox{\bf w}(\cdot,0),\mbox{\bf  w}(\cdot,0)
)\Vert_{H^{\frac 14,\frac 12, \widetilde s}(\Bbb R^n)\times H^{\frac 12,\frac 32,\widetilde s}(\Bbb R^n)}
+\root\of{ \vert\widetilde s\vert} \Vert
\mbox{\bf w}\Vert_{H^{\frac 12,1,\widetilde s}(\mathcal Q)}
                         \\
\le C_1\sum_{\nu=0}^N \left(\Vert \eta_\ell(
\partial_{y_n} (\chi_\nu\mbox{\bf w})(\cdot,0),\chi_\nu\mbox{\bf w}
(\cdot,0))\Vert_{H^{\frac 14,\frac 12, \widetilde s}(\Bbb R^n)\times H^{\frac 12,\frac 32,\widetilde s}(\Bbb R^n)}
+ \root\of{\vert \widetilde s\vert}\Vert\eta_\ell
\chi_\nu\mbox{\bf w}\Vert_{H^{\frac 12,1,\widetilde s}(\mathcal Q)}+\right.\nonumber\\
% \sum_{\nu=0}^N
 \left.\Vert (1-\eta_\ell)(
\partial_{y_n} (\chi_\nu\mbox{\bf w})(\cdot,0),\chi_\nu\mbox{\bf w}
(\cdot,0))\Vert_{H^{\frac 14,\frac 12, \widetilde s}(\Bbb R^n)\times H^{\frac 12,\frac 32,\widetilde s}(\Bbb R^n)}
+ \root\of{\vert \widetilde s\vert}\Vert(1-\eta_\ell)
\chi_\nu\mbox{\bf w}\Vert_{H^{\frac 12,1,\widetilde s}(\mathcal Q)}\right)\nonumber\\
\le C_{2}\biggl(\root\of{\widetilde \varphi(y_0^*)}\Vert \text{\bf w}
\Vert_{H^{0,1,\widetilde s}(\mathcal Q)} +\Vert
\mbox{\bf g}e^{\vert s\vert \varphi}\Vert_{H^{\frac 14,\frac 12, \widetilde s}(\Bbb R^n)}
+ \sum_{\nu=0}^N\Vert \mbox{\bf P}(y,D,\widetilde s)\mbox{\bf w}_\nu\Vert
_{L^2(\mathcal Q)}\nonumber\\
+ \Vert (
\partial_{y_n} \mbox{\bf w}(\cdot,0),\mbox{\bf  w}(\cdot,0)
)\Vert_{H^{\frac 14,\frac 12, \widetilde s}(\Bbb R^n)\times H^{\frac 12,\frac 32,\widetilde s}(\Bbb R^n)}\biggr)                     \nonumber\\
+C_3\sum_{\nu=0}^N \left(\Vert (1-\eta_\ell)(
\partial_{y_n} (\chi_\nu\mbox{\bf w})(\cdot,0),\chi_\nu\mbox{\bf w}
(\cdot,0))\Vert_{H^{\frac 14,\frac 12, \widetilde s}(\Bbb R^n)\times H^{\frac 12,\frac 32,\widetilde s}(\Bbb R^n)}
+ \root\of{\vert \widetilde s\vert}\Vert(1-\eta_\ell)
\chi_\nu\mbox{\bf w}\Vert_{H^{0,1,\widetilde s}(\mathcal Q)}\right)\nonumber .
\end{eqnarray}
By (\ref{nikolas}) there exist a constant $C_4$ independent of $\widetilde s,\ell$ and $\nu$ such that
\begin{eqnarray}\label{nikolas1}
\sum_{\nu=0}^N \left(\Vert (1-\eta_\ell)(
\partial_{y_n} (\chi_\nu\mbox{\bf w})(\cdot,0),\chi_\nu\mbox{\bf w}
(\cdot,0))\Vert_{H^{\frac 14,\frac 12, \widetilde s}(\Bbb R^n)\times H^{\frac 12,\frac 32,\widetilde s}(\Bbb R^n)}
+ \root\of{\vert \widetilde s\vert}\Vert(1-\eta_\ell)\chi_\nu
\mbox{\bf w}\Vert_{H^{\frac 12,1,\widetilde s}(\mathcal Q)}\right)\nonumber\\
\le
C_4\left(\root\of{\widetilde \varphi(y_0^*)}\Vert (
\partial_{y_n} \mbox{\bf w}(\cdot,0),\mbox{\bf w}
(\cdot,0))\Vert_{L^2(\Bbb R^n)\times H^{0,1,\widetilde s}(\Bbb R^n)}
+ \root\of{\widetilde \varphi(y_0^*)}\Vert
\mbox{\bf w}\Vert_{H^{\frac 12,1,\widetilde s}(\mathcal Q)}\right)
\end{eqnarray}
Using  inequality (\ref{nikolas1}) in order to estimate the last terms in
(\ref{golifax}) we obtain
\begin{eqnarray}
\Vert (
\partial_{y_n} \mbox{\bf w}(\cdot,0),\mbox{\bf  w}(\cdot,0)
)\Vert_{H^{\frac 14,\frac 12, \widetilde s}(\Bbb R^n)\times H^{\frac 12,\frac 32,\widetilde s}(\Bbb R^n)}
+ \root\of{\vert\widetilde s\vert}\Vert
\mbox{\bf w}\Vert_{H^{\frac 12,1,\widetilde s}(\mathcal Q)}
                         \nonumber\\
\le C_{5}\biggl(\root\of{\widetilde \varphi(y_0^*)}\Vert \text{\bf w}
\Vert_{H^{\frac 12,1,\widetilde s}(\mathcal Q)} +\Vert
\mbox{\bf g}e^{\vert s\vert \varphi}\Vert_{H^{\frac 14,\frac 12,\widetilde s}(\Bbb R^n)}
+\Vert \mbox{\bf P}(y,D,\widetilde s)
\mbox{\bf w}\Vert_{L^2(\mathcal Q)}\nonumber\\
+ \root\of{\widetilde \varphi(y_0^*)}\Vert (
\partial_{y_n} \mbox{\bf w}(\cdot,0),\mbox{\bf w}
(\cdot,0))\Vert_{L^2(\Bbb R^n)\times H^{0,1,\widetilde s}(\Bbb R^n)}
+ \sum_{\nu=0}^N\Vert[\chi_\nu,\mbox{\bf  P}(y,D,\widetilde s)]\mbox{\bf
w}_{1}\Vert_{L^2(\mathcal Q)}
\biggr)  +           \nonumber\\
C_6\left(\root\of{\widetilde \varphi(y_0^*)}\left\Vert \left(
\partial_{y_n}\mbox{\bf w}(\cdot,0),\mbox{\bf w}
(\cdot,0)\right)\right\Vert_{H^{\frac 14,\frac 12, \widetilde s}(\Bbb R^n)\times H^{\frac 12,\frac 32,\widetilde s}(\Bbb R^n)}
+ \root\of{\widetilde \varphi(y_0^*)}\Vert
\mbox{\bf w}\Vert_{H^{\frac 12,1,\widetilde s}(\mathcal Q)}\right)\nonumber\\
\le C_{7}\biggl(\root\of{\widetilde \varphi(y_0^*)}\Vert \text{\bf w}
\Vert_{H^{\frac 12,1,\widetilde s}(\mathcal Q)} +\Vert
\mbox{\bf g}e^{\vert s\vert\varphi}\Vert_{ H^{\frac 14,\frac 12,\widetilde s}(\Bbb R^n)}
+ \Vert\mbox{\bf  P}(y,D,\widetilde s)\mbox{\bf
w}\Vert_{L^2(\mathcal Q)}\nonumber\\
+ \Vert (
\partial_{y_n} \mbox{\bf w}(\cdot,0),\mbox{\bf  w}(\cdot,0)
)\Vert_{H^{\frac 14,\frac 12, \widetilde s}(\Bbb R^n)\times H^{\frac 12,\frac 32,\widetilde s}(\Bbb R^n)}\biggr)+\nonumber\\
C_8\left(\root\of{\widetilde \varphi(y_0^*)}\Vert (
\partial_{y_n} \mbox{\bf w}(\cdot,0),\mbox{\bf w}
(\cdot,0))\Vert_{L^2(\Bbb R^n)\times H^{0,1,\widetilde s}(\Bbb R^n)}
+ \root\of{\widetilde \varphi(y_0^*)}\Vert
\mbox{\bf w}\Vert_{H^{\frac 12,1,\widetilde s}(\mathcal Q)}\right).
\end{eqnarray}
Hence there exists $s_0>1$  such that for all
$s\ge s_0$  we see
\begin{eqnarray}\label{zmpolit}
\Vert(
\partial_{y_n} \mbox{\bf w}(\cdot,0),\mbox{\bf  w}(\cdot,0)
)\Vert_{H^{\frac 14,\frac 12, \widetilde s}(\Bbb R^n)\times H^{\frac 12,\frac 32,\widetilde s}(\Bbb R^n)}
+\root\of{ \vert \widetilde s\vert\Vert}
\mbox{\bf w}\Vert_{H^{\frac 12,1,\widetilde s}(\mathcal Q)}\nonumber\\
\le C_{9}\biggl(\root\of{\widetilde \varphi(y_0^*)}\Vert \text{\bf w}
\Vert_{H^{0,1,\widetilde s}(\mathcal Q)} +\Vert
\mbox{\bf g}e^{ s\varphi}\Vert_{ H^{\frac 14,\frac 12,\widetilde s}(\Sigma)}
+ \Vert \mbox{\bf P}(y,D,\widetilde s)\mbox{\bf
w}\Vert_{L^2(\mathcal Q)}  \biggr).
\end{eqnarray}
Proof of Proposition \ref{zoopa} is complete. $\blacksquare$

The remaining part of this section we obtain from the estimate (\ref{3.2'1'}) the Carleman estimate (\ref{2.9'}).

The short calculations imply that the functions
 the functions $d\omega_{\mbox{\bf u}}, \text{div}\,
\mbox{\bf u}$ satisfy the equations
\begin{eqnarray}\label{3.3}
P_{\rho,\mu}(x,D) d\omega_{\mbox{\bf u}} -\int_0^{x_0}\widetilde
\mu(x,\widetilde x_0) \Delta d\omega_{\mbox{\bf u}}\, d\widetilde x_0
= \mbox{\bf q}_1 \qquad \mbox{in}\,Q,\nonumber\\
P_{\rho,\lambda+2\mu}(x,D)\text{div}\, \mbox{\bf u}-\int_0^{x_0}
(\widetilde \lambda+2\widetilde \mu)(x,\widetilde x_0) \Delta \text{div}\,
\mbox{\bf u}\, d\widetilde x_0= q_2\quad \mbox{in}\,Q,
\end{eqnarray}
where
\begin{eqnarray}\label{3.4}
\mbox{\bf q}_1=K_1 (x,D)d\omega_{\mbox{\bf u}} + K_2 (x,D) \text{div}\,
\mbox{\bf u}+\int_0^{x_0}(\widetilde K_1 (x,\widetilde x_0,D)d\omega
_{\mbox{\bf u}} \nonumber\\+\widetilde K_2 (x,\widetilde x_0,D) \text{div}\,
\mbox{\bf u})d\widetilde x_0+\rho d\omega_{\mbox{\bf F}/\rho},\\
\quad q_2=K_3 (x,D) d\omega_{\mbox{\bf u}} + K_4 (x,D) \text{div}\,
\mbox{\bf u}+\int_0^{x_0}(\widetilde K_3 (x,\widetilde x_0,D)
d\omega_{\mbox{\bf u}}\nonumber\\+\widetilde K_4 (x,\widetilde x_0,D) \text{div} \,
\mbox{\bf u})d\widetilde x_0+\rho\text{div} \,(\mbox{\bf F}/\rho),
\end{eqnarray}
where $K_j (x,D), \widetilde K_j (x,\widetilde x_0,D)$ are first order differential
operators with
$C^1$ coefficients.

Now we introduce new unknown function $\mbox{\bf v}=(\mbox{\bf v}_1,v_2)$  by formulae
\begin{equation}\label{victory}
\mbox{\bf v}_1=d\omega_{\mbox{\bf u}}+ \int_0^{x_0}\frac{\widetilde
\mu(x,\widetilde x_0)}{\mu(x)} d\omega_{\mbox{\bf u}} \,d\widetilde x_0,
\quad v_2=\text{div}\, \mbox{\bf u}+\int_0^{x_0}\frac{(\widetilde \lambda
+2\widetilde \mu)(x,\widetilde x_0)}{(\lambda+2\mu)(x)} \text{div}\,
\mbox{\bf u}\, d\widetilde x_0.
\end{equation}
More specifically
$$
\mbox{\bf v}_1=(v_{1,2}, \dots  v_{n-1,n}),\quad  v_{k,j}=\partial_{x_j} u_k-\partial_{x_k} u_j+\int_0^{x_0}\frac{\widetilde
\mu(x,\widetilde x_0)}{\mu(x)}(\partial_{x_j} u_k-\partial_{x_k} u_j)(\widetilde x_0,x')d\widetilde x_0.
$$
Then from (\ref{3.3}) we have
\begin{equation}\label{pokozad}
\mbox{\bf P}(x,D) \mbox{\bf v}=(
P_{\rho,\mu} (x,D)\mbox{\bf v}_1, P_{\rho,\lambda+2\mu} (x,D) v_2)=\mbox{\bf q} \quad \mbox{in}\,\, Q,
\end{equation}
where $\mbox{\bf q}=(\mbox{\bf q}_3,q_4)=(\mbox{\bf q}_1+\widetilde{\mbox{\bf q}}_1,q_2
+\widetilde q_2):$
\begin{equation}\label{slava1}
\,\, \widetilde{\mbox{\bf q}}_1=-\mu\int_0^{x_0}\left (2(\nabla'\frac{\widetilde
\mu(x,\widetilde x_0)}{\mu(x)},\nabla' d\omega_{\mbox{\bf u}})+ \Delta\left(
\frac{\widetilde \mu(x,\widetilde x_0)}{\mu(x)}\right)d\omega_{\mbox{\bf u}}\right)
\,d\widetilde x_0+\rho\partial_{x_0}\int_0^{x_0}\frac{\widetilde \mu(x,
\widetilde x_0)}{\mu(x)} d\omega_{\mbox{\bf u}} \,d\widetilde x_0,
\end{equation}
\begin{eqnarray}\label{slava2}
 \widetilde q_2=-(\lambda+2\mu)\int_0^{x_0}\left (2(\nabla'\frac{(\widetilde\lambda
+2\widetilde \mu)(x,\widetilde x_0)}{\lambda(x)+2\mu(x)},\nabla'  \text{div} \,
\mbox{\bf u})+ \Delta\left(\frac{\widetilde \mu(x,\widetilde x_0)}{(\lambda+2\mu)(x)}\right)
\text{div} \, \mbox{\bf u}\right)\,d\widetilde x_0\nonumber\\
+\rho\partial_{x_0}\int_0^{x_0}\frac{(\widetilde\lambda+2\widetilde \mu)
(x,\widetilde x_0)}{(\lambda+2\mu)(x)}  \text{div} \, \mbox{\bf u}
\,d\widetilde x_0.
\end{eqnarray}
%The formulae (\ref{slava1}) and (\ref{slava2}) immediately yield the estimate
%\begin{eqnarray}\label{nina}
%\Vert \mbox{\bf q}e^{s\varphi}\Vert_{L^2(Q)}\le C_{11}( \Vert\nabla'  d\omega_{\mbox{\bf u}} e^{s\varphi}\Vert
%_{L^2(Q)}
%+\Vert\nabla'  \mbox{div}\, \mbox{\bf u}e^{s\varphi}\Vert_{L^2(Q)}+\Vert d\omega
%_{\mbox{\bf u}}e^{s\varphi}\Vert_{L^2(Q)}\nonumber\\+\Vert  \mbox{div}\,\mbox{\bf u}e^{s\varphi}\Vert_{L^2(Q)}+\Vert d\omega_{\mbox{\bf %F}}e^{s\varphi}
%\Vert_{L^2(Q)}+\Vert  \mbox{div}\, \mbox{\bf F}e^{s\varphi}\Vert_{L^2(Q)}).
%\end{eqnarray}

On the next step we show that the boundary conditions (\ref{poko1}) holds true for some function $\mbox{\bf g}.$

We start with the following proposition:
\begin{proposition}\label{gonka}
Let $R(x,\widetilde x_0, D')$ be a partial differential
operator with smooth coefficients of order $\widetilde j$ and $\kappa\ge 0.$
Then there exist a constant $C_{10}$ independent of $s$ such that
\begin{equation}\label{qupol}
\left\Vert e^{s\varphi}\widetilde\varphi^\kappa \int_0^{x_0}R(x,\widetilde x_0, D')
\mbox{ \bf u}
(\widetilde x_0,\cdot)d\widetilde x_0\right\Vert_{L^2[-T,T]}
\le \frac{C_{10}}{\root\of{s}}\sum_{\vert\alpha\vert\le \widetilde j,\alpha_0=0}
\Vert  e^{s\varphi}\widetilde\varphi^\kappa\partial^{\alpha}\mbox{ \bf u}
\Vert_{L^2[-T,T]}
\end{equation}
for all $s\ge 1.$
If $r_0\in C^1(\overline\Sigma)$ and $ p\in (0,1)$ there exists a constant $C_{11}$ independent of $s$ such that
\begin{equation}\label{granata}
\Vert e^{s\varphi}\int_0^{x_0}r_0
\mbox{ \bf u}
(\widetilde x_0,\cdot)d\widetilde x_0\Vert_{H^{\frac p 2,p,\widetilde s}(\Sigma)}
\le \frac{C_{11}}{\root\of{s}}
\Vert  e^{s\varphi}\mbox{ \bf u}
\Vert_{H^{\frac p 2,p,\widetilde s}(\Sigma)}.
\end{equation}

\end{proposition}

{\bf Proof.}
Instead of inequality (\ref{qupol}) it suffices to prove
\begin{equation}\label{qupol1}
\left\Vert e^{s\varphi}\widetilde\varphi^\kappa \int_0^{x_0}\vert \mbox{ \bf u}
(\widetilde x_0,\cdot)\vert d\widetilde x_0\right\Vert_{L^2[-T,T]}
\le \frac{C_{12}}{\root\of{s}}
\Vert  e^{s\varphi}\widetilde\varphi^\kappa\mbox{ \bf u}
\Vert_{L^2[-T,T]} \quad \forall x'\in \Omega.
\end{equation}
The inequality  (\ref{qupol}) is equivalent to
\begin{equation}\label{qupol1}
\left\Vert e^{s\varphi}\widetilde\varphi^\kappa \int_0^{x_0}\vert g
(\widetilde x_0)\vert d\widetilde x_0\right\Vert_{L^2[0,T]}
\le \frac{C_{13}}{\root\of{s}}
\Vert  e^{s\varphi}\widetilde\varphi^\kappa g
\Vert_{L^2[0,T]}\quad \forall s\ge 1,
\end{equation} where $g$ in an arbitrary function from $L^2(0,T).$
Applying the Cauchy inequality we have
$$
\left\Vert e^{s\varphi} \widetilde \varphi^\kappa\int_0^{x_0}\vert g
(\widetilde x_0)\vert d\widetilde x_0\right\Vert_{L^2[0,T]}^2\le \int_0^T x_0
\widetilde\varphi^{2\kappa} e^{2s\varphi} \int_0^{x_0} g^2
(\widetilde x_0) d\widetilde x_0 dx_0=
$$
$$\int_0^T\frac{x_0\widetilde\varphi^{2\kappa}}{2s\partial_{x_0}\varphi}
(\partial_{x_0}e^{2s\varphi})\int_0^{x_0} g^2
(\widetilde x_0) d\widetilde x_0 dx_0
\le -\int_0^T\partial_{x_0}\left (\frac{x_0\widetilde\varphi^{2\kappa}}
{2s\partial_{x_0}\varphi}\int_0^{x_0} g^2
(\widetilde x_0) d\widetilde x_0\right) e^{2s\varphi} dx_0=
$$
$$
-\int_0^T\partial_{x_0}\left (\frac{x_0\widetilde\varphi^{2\kappa}}
{2s\partial_{x_0}\varphi}\right)\int_0^{x_0} g^2
(\widetilde x_0) d\widetilde x_0 e^{2s\varphi} dx_0
-\int_0^T\frac{x_0\widetilde\varphi^{2\kappa}}{2s\partial_{x_0}\varphi} g^2
 e^{2s\varphi} dx_0.
$$
By (\ref{pinok22}) we see that
$\frac{x_0\widetilde \varphi^{2\kappa}}{\partial_{x_0}\varphi}\in W^{1}_\infty [0,T/2]$  and $\vert \frac{x_0\widetilde \varphi^{2\kappa}}{\partial_{x_0}\varphi} \vert \le C_{14}\widetilde\varphi^{2\kappa} .$ Therefore
$$
\left\Vert e^{s\varphi} \widetilde\varphi^{\kappa}\int_0^{x_0}\vert g
(\widetilde x_0)\vert d\widetilde x_0\right\Vert_{L^2[0,T]}^2
$$
$$
\le C_{15}(\int_0^T\frac{\widetilde \varphi^{2\kappa-1}}{2s}
e^{2s\varphi}\int_0^{x_0} g^2
(\widetilde x_0) d\widetilde x_0 dx_0
+ \int_0^T\frac{\widetilde\varphi^{2\kappa-1}}{2s} g^2
 e^{2s\varphi} dx_0)
$$
$$
\le \int_0^T\frac{C_{16}\widetilde\varphi^{2\kappa-1}}{2s} g^2 e^{2s\varphi} dx_0.
$$
Here in order to obtain the last inequality,
we used (\ref{pinok22}) and (\ref{giga1}).

In order to prove the estimate (\ref{granata})  we first show that
\begin{eqnarray}\label{OAA}
\left\Vert e^{s\varphi}\int_0^{x_0}r_0
\mbox{ \bf u}
(\widetilde x_0,\cdot)d\widetilde x_0\right\Vert_{H^{0,p,\widetilde s}(\Sigma)}
\le C_{17}\Vert e^{s\varphi}
\mbox{ \bf u}
\Vert_{H^{0,p,\widetilde s}(\Sigma)}.
\end{eqnarray}
Consider the operator $\mathcal K w=e^{s\varphi(x)}\int_0^{x_0}e^{-s\varphi(\widetilde x_0,x')}r_0(\widetilde x_0)w(\widetilde x_0,x')d\widetilde x_0$.
By (\ref{pinok22})
\begin{equation}\label{AAA}
\Vert \mathcal K\Vert_{\mathcal L(L^2(\Sigma);L^2(\Sigma))}\le C_{18}(\Vert r_0\Vert_{L^\infty(\Sigma)})
\end{equation}
and by (\ref{qupol1})
\begin{equation}\label{AAA1}
\Vert \mathcal K\Vert_{\mathcal L(H^{0,1,\widetilde s}(\Sigma);H^{0,1,\widetilde s}(\Sigma))}\le C_{19}(\Vert r_0\Vert_{W^1_\infty(\Sigma)}).
\end{equation} From (\ref{AAA}), (\ref{AAA1}) and interpolation argument  we obtain (\ref{OAA}).
 By (\ref{qupol1})
\begin{equation}\label{AAA2}
\Vert \mathcal K\Vert_{\mathcal L(H^{1,0,\widetilde s}(\Sigma);H^{1,0,\widetilde s}(\Sigma))}\le C_{20}(\Vert r_0\Vert_{W^1_\infty(\Sigma)}).
\end{equation}
From (\ref{AAA}), (\ref{AAA2}) and interpolation argument  we have
\begin{eqnarray}\label{OAA1}
\Vert e^{s\varphi}\int_0^{x_0}r_0
\mbox{ \bf u}
(\widetilde x_0,\cdot)d\widetilde x_0\Vert_{H^{\frac{p}{2},0,\widetilde s}(\Sigma)}
\le C_{21} \Vert e^{s\varphi}
\mbox{ \bf u}\Vert_{H^{\frac{p}{2},0,\widetilde s}(\Sigma)}.
\end{eqnarray} Estimates (\ref{OAA}) and (\ref{OAA1}) imply (\ref{granata}).
$\blacksquare$

We have
\begin{proposition}\label{loskut}
Let function $\mbox{\bf v}$ is given by (\ref{victory}) and $\uu$ satisfy
(\ref{gora1A}). Then $\mbox{\bf v}$ satisfies problem (\ref{poko})
- (\ref{poko1}) with  functions $\mbox{\bf g}$ and $\mbox{\bf q}$ such that
\begin{eqnarray}\label{nina}
\Vert \mbox{\bf q}e^{s\varphi}\Vert_{L^2(Q)}\le C_{22}( \Vert(\nabla'
d\omega_{\mbox{\bf u}}) e^{s\varphi}\Vert
_{L^2(Q)}
+\Vert(\nabla'  \mbox{div}\, \mbox{\bf u})e^{s\varphi}\Vert_{L^2(Q)}+\Vert (d\omega
_{\mbox{\bf u}})e^{s\varphi}\Vert_{L^2(Q)}\nonumber\\+\Vert  (\mbox{div}\,\mbox{\bf u})e^{s\varphi}\Vert_{L^2(Q)}+\Vert (d\omega_{\mbox{\bf F}})e^{s\varphi}
\Vert_{L^2(Q)}+\Vert  (\mbox{div}\, \mbox{\bf F})e^{s\varphi}\Vert_{L^2(Q)}+\Vert \mbox{\bf F}e^{s\varphi}\Vert_{L^2(Q)})
\end{eqnarray}
and
\begin{equation}\label{robot}
\Vert\mbox{\bf  g}e^{s\varphi}\Vert_{H^{\frac 14, \frac 12, \widetilde s}(\Sigma)}
\le C_{23}(\Vert \mbox{\bf F}e^{s\varphi}\Vert_{H^{\frac 14, \frac 12, \widetilde s}(\Sigma)}
+ \Vert(\partial_{\vec \nu}
\mbox{\bf u})e^{s\varphi}\Vert_{H^{\frac 14, \frac 12, \widetilde s}(\Sigma)}).
\end{equation}
%and
%\begin{equation}\label{robot0}
%\Vert\root\of{\varphi}(\partial_{x_0}\mbox{\bf  g})e^{s\varphi}\Vert_{L^2(\Sigma)}
%\le C_{3}(\Vert\root\of{\varphi}(\partial_{x_0} \mbox{\bf F})e^{s\varphi}\Vert_{L^2(\Sigma)}
%+ \Vert\root\of{\varphi}(\partial_{\vec \nu}\partial_{x_0}
%\mbox{\bf u})e^{s\varphi}\Vert_{L^2(\Sigma)}).
%\end{equation}
\end{proposition}

{\bf Proof.} Let function $\mbox{\bf q}=(\mbox{\bf q}_3,q_4)$ given by (\ref{3.4})- (\ref{slava2}).
Applying (\ref{qupol}) to estimate $\mbox{\bf q}$, we obtain (\ref{nina}).
Next we show that the function $\mbox{\bf v}$ satisfies the boundary conditions
(\ref{poko1}).
We set all the components of the  function $\mbox{\bf g}$ starting from $n+1$
be equal to zero:
\begin{equation} \label{makaka} g_k=0\quad\mbox{for all}\quad k\ge n+1.
\end{equation}
By (\ref{gora1A}) and the zero Dirichlet boundary conditions for the function
$\mbox{\bf u}$ we have
\begin{equation}\label{gora1AA}
- L_{\lambda,\mu}(x,D')\mbox{\bf u}
+\int_0^{x_0}L_{\widetilde \lambda, \widetilde\mu}
(x, \widetilde x_0, D')\mbox{\bf u}(\widetilde x_0, x')d\widetilde x_0
= \mbox{\bf F} \quad \mbox{on}\,\,\Sigma .
\end{equation}
Next we move all the terms containing the first
derivatives of the function
$\mbox{\bf u}$ into the right-hand side, divide both sides by $\mu$ and
denote the right-hand side of obtained equality  as $g_k$.
Then the first  $n$ components of the function $\mbox{\bf g}$ are defined by formula:
\begin{equation}\label{gorbun}
(g_1,\dots, g_n)=\frac{1}{\mu}(\mbox{\bf F}-(\mbox{div}\,\mbox{\bf u})\nabla'\lambda
- (\nabla'\mbox{\bf u}
+ (\nabla'\mbox{\bf u})^T) \nabla'\mu
\end{equation}
$$
- \int_0^{x_0}((\mbox{div}\,\mbox{\bf u})\nabla'\widetilde \lambda
+ (\nabla'\mbox{\bf u}
+ (\nabla'\mbox{\bf u})^T) \nabla'\widetilde \mu) d\widetilde x_0)\quad
\mbox{on}\quad \Sigma .
$$
Then by (\ref{1.1}), (\ref{1.2}) we have
\begin{eqnarray}\label{goraAAA}
-\Delta \mbox{\bf u}-\frac{(\lambda+\mu)}{\mu}\nabla'\mbox{div}\,\mbox{\bf u}
-\int_0^{x_0}\left(\frac{\widetilde \mu(x,\widetilde x_0)}{\mu(x)}
\Delta \mbox{\bf u}
(\widetilde x_0,x')+\frac{\widetilde \lambda(x,\widetilde x_0)
+ \widetilde \mu(x,\widetilde x_0)}
{\mu(x)}\nabla'\mbox{div}\,\mbox{\bf u}\right) dx_0\nonumber\\
= (g_1,\dots, g_n)\quad \mbox{on}\quad \Sigma.
\end{eqnarray}
By (\ref{gorbun}) and $\mbox{\bf u}\vert_{\Sigma}=0$ we have
\begin{eqnarray}\label{robot11}
\sum_{j=1}^n\Vert  g_je^{s\varphi}\Vert_{H^{\frac 14, \frac 12, \widetilde s}(\Sigma)}
\le C_{24}(\Vert \mbox{\bf F}e^{s\varphi}\Vert_{H^{\frac 14, \frac 12, \widetilde s}(\Sigma)}
+ \Vert(\partial_{\vec \nu}
\mbox{\bf u})e^{s\varphi}\Vert_{H^{\frac 14, \frac 12, \widetilde s}(\Sigma)}            \\
+ \sum_{k,m=1}^n\Vert\int_0^{x_0}p_{k,m}(\widetilde x_0,x)
\partial_{\vec \nu} u_m(\widetilde x_0,x') d\widetilde x_0e^{s\varphi}\Vert
_{H^{\frac 14, \frac 12, \widetilde s}(\Sigma)}).                \nonumber
\end{eqnarray}
Here $p_{k,m} \in C^1([-T,T]^2 \times\overline \OOO)$ are some functions.
Estimating the last term on the right-hand side  of (\ref{robot11})
using (\ref{granata}) we obtain (\ref{robot}).

For any index $\widehat j\in\{1,\dots, n\}$ the short computations imply
\begin{eqnarray}\label{gavnoed}
-\Delta u_{\widehat j} =\sum_{j=1, j\ne \widehat j}^n
-\partial_{x_j}(\partial_{x_j} u_{\widehat i}-\partial_{x_{\widehat j}} u_j)
-\partial^2_{x_{\widehat i}} u_{\widehat i} -\sum_{j=1,j\ne \widehat j}^n
\partial_{x_j}\partial_{x_{\widehat i} }u_j\nonumber\\
=-\sum^n_{j=1, j\ne \widehat j}\partial_{x_j}(\partial_{x_j} u_{\widehat j}
-\partial_{x_{\widehat j}} u_j)-\partial_{x_{\widehat j} }\text{div}\,
\mbox{\bf u }\quad \text{in}\,\,\Omega.
\end{eqnarray}
Using the equality (\ref{gavnoed}) we rewrite $ \widehat j$-th equation
in (\ref{goraAAA}) as
$$
b_{\hat j}(x,D)\mbox{\bf v}
= -\sum^n_{j=1, j\ne \widehat j}\mbox{sign}(j-\hat j)\partial_{x_j}
v_{\hat j, j}-\frac{\lambda+2\mu}{\mu}\partial_{x_{\hat j}}v_2=g_{\hat j}\quad
\mbox{on}\,\,\Sigma.
$$
The construction of the operator $B_1(x,D)$ is complete. Now we construct
matrix $B_2(x').$
Next $v_{k,j}(x)=\nu_j(x')\partial_{\vec \nu} u_k-\nu_i(x')\partial_{\vec \nu}
u_j, 1\le k<j\le n.$ By the boundary condition (\ref{gora2A})
$$
v_{k,j}(y^*)=0 \quad \mbox{for} \,\,1\le k<j<n,\quad v_{j,n}(y^*)
= -\partial_{x_n} u_j(y^*),\quad v_2=-\partial_{x_n} u_j(y^*).
$$
Set $\widetilde{\mbox{\bf v}}=(v_{1,n},\dots, v_{2,n},\dots v_{n-1,n}, v_2).$
Obviously in the small neighborhood of  $y^*$ there exists a  smooth matrix
$B_3(x')$ such that
$$
(\partial_{x_n} u_1,\dots, \partial_{x_n} u_n)=B_3(x') \widetilde{\mbox{\bf v}}
\quad \forall x\in \Sigma\cap B(y^*,\delta).
$$
Then $B_2(x')\mbox{\bf v}=\mbox{\bf v}-(\nu_2\partial_{x_2} u_1-\nu_1
\partial_{x_1} u_2, \dots , \nu_n \partial_{x_n} u_1-\nu_1 \partial_{x_1} u_n,
\dots , \nu_n \partial_{x_n} u_{n-1}-\nu_{n-1}\partial_{x_{n-1}}
u_{n},$\newline $ \sum_{j=1}^n \nu_j\partial_{x_j} u_j)=0.$
The proof of Proposition \ref{loskut} is complete.
$\blacksquare$

By Proposition \ref{loskut} for the function $\mbox{\bf v}$ given by formulae (\ref{victory}) the Carleman estimate (\ref{3.2'1'}) holds true.
Estimating the right hand side of (\ref{3.2'1'}) using the inequalities (\ref{nina}) and (\ref{robot}) we have
\begin{eqnarray}\label{zmpolit1}
\Vert (\partial_{\vec \nu}\mbox{\bf v},\mbox{\bf  v})e^{s\varphi}\Vert_{H^{\frac 14,\frac 12, \widetilde s}(\Sigma)\times H^{\frac 12,\frac 32,\widetilde s}(\Sigma)} + \root\of{\int_{Q}\sum_{\vert \alpha\vert\le 2}
(s\widetilde\varphi)^{3-2\vert\alpha\vert}\vert
\partial^\alpha\mbox{\bf v}\vert^2e^{2s\varphi}dx}
\\
\le C_{25}\left(\Vert \mbox{\bf F}e^{s\varphi}\Vert_{H^{\frac 14, \frac 12, \widetilde s}(\Sigma)}
+ \Vert(\partial_{\vec \nu}
\mbox{\bf u})e^{s\varphi}\Vert_{H^{\frac 14, \frac 12, \widetilde s}(\Sigma)}
+\Vert (\partial_{\vec \nu}\mbox{\bf v},\mbox{\bf  v})e^{s\varphi}\Vert_{H^{\frac 14,\frac 12, \widetilde s}(\widetilde\Sigma)\times H^{\frac 12,\frac 32,\widetilde s}(\widetilde\Sigma)}
\right.\nonumber\\
+\Vert(\nabla'
d\omega_{\mbox{\bf u}}) e^{s\varphi}\Vert
_{L^2(Q)}
+\Vert(\nabla'  \mbox{div}\, \mbox{\bf u})e^{s\varphi}\Vert_{L^2(Q)}+\Vert (d\omega
_{\mbox{\bf u}})e^{s\varphi}\Vert_{L^2(Q)}\nonumber\\
\left.+\Vert  (\mbox{div}\,\mbox{\bf u})e^{s\varphi}\Vert_{L^2(Q)}+\Vert (d\omega_{\mbox{\bf F}})e^{s\varphi}
\Vert_{L^2(Q)}+\Vert  (\mbox{div}\, \mbox{\bf F})e^{s\varphi}\Vert_{L^2(Q)}+\Vert \mbox{\bf F}e^{s\varphi}\Vert_{L^2(Q)}\right ).\nonumber
\end{eqnarray}
By  (\ref{pinok22}) and (\ref{qupol})
%\begin{equation}\label{qupol}
% \Vert e^{s\varphi}\int_0^{x_0} g(\widetilde x_0)d\widetilde x_0\Vert_{L^2(-T,T)}\le o(1)\Vert e^{s\varphi} g\Vert_{L^2(-T,T)}\quad \mbox{as}\,\, %s\rightarrow +\infty.
%\end{equation}%
%This inequality imply that
 for all sufficiently large $s$
\begin{equation}\label{kom}
\int_Q \sum_{\vert \alpha\vert\le 2}(s\widetilde\varphi)^{3-2\vert\alpha\vert}\vert
\partial^\alpha(d \omega_{\mbox {\bf u}},\mbox{div}\,\mbox {\bf u})\vert^2e^{2s\varphi}dx\le C_{26}
 \int_{Q}\sum_{\vert \alpha\vert\le 2}
(s\widetilde\varphi)^{3-2\vert\alpha\vert}\vert
\partial^\alpha\mbox{\bf v}\vert^2e^{2s\varphi}dx.
\end{equation}
Next we prove the following:
\begin{proposition}\label{govno} Let positive $\delta$ be sufficiently small.
There exist $s_0>0$ such that for all $s\ge s_0$  there exists a constant $C_{27}>0$ independent of $s$
such that
\begin{equation}\label{govno1}
\Vert\partial^2_{\vec \nu} \text{\bf u}
e^{s\varphi}\Vert_{H^{\frac 14,\frac 12,\widetilde s}( \Sigma)}
+\Vert \partial_{\vec \nu}
\text{\bf u}
e^{s\varphi}\Vert_{H^{\frac 34,\frac 32,\widetilde s}( \Sigma)}\le C_{27} (\Vert \mbox{\bf v}
e^{s\varphi}\Vert_{H^{\frac 34,\frac 32, \widetilde s}(\Sigma)}
+\Vert \F e^{s\varphi}\Vert_{H^{\frac 14,\frac 12, \widetilde s}(\Sigma)}).
\end{equation}
\end{proposition}

{\bf Proof.} By (\ref{qupol}) there exists a constant $C_{28}$ independent of $s$ such that
\begin{eqnarray}\label{zina1}
\Vert (d \omega_{\mbox {\bf u}},\mbox{div}\,\mbox {\bf u})
e^{s\varphi}\Vert_{H^{\frac 34,\frac 32, \widetilde s}(\Sigma)}
\le C_{28}
\Vert  \mbox{\bf v}
e^{s\varphi}\Vert_{H^{\frac 34,\frac 32, \widetilde s}(\Sigma)}.
\end{eqnarray}
Thanks to the zero Dirichlet boundary conditions on $ \Sigma$ there exists a smooth matrix $\mathcal V(x)$ such that
$\partial_{\vec \nu} \text{\bf u}=\mathcal V(x)(d \omega_{\mbox {\bf u}},\mbox{div}\,\mbox {\bf u}).$
Therefore
\begin{eqnarray}\label{zina}
\Vert\partial_{\vec \nu}
\text{\bf u}
e^{s\varphi}\Vert_{H^{\frac 34,\frac 32, \widetilde s}(\Sigma)}\le C_{29}
\Vert  (d \omega_{\mbox {\bf u}},\mbox{div}\,\mbox {\bf u})
e^{s\varphi}\Vert_{H^{\frac 34,\frac 32, \widetilde s}(\Sigma)}
.
\end{eqnarray}
From equation (\ref{gora1A}) on $\Sigma$ we have
\begin{equation}\label{eb1}
\partial_{\vec \nu}^2  \text{\bf u}=\widetilde A(x, D')\partial_{\vec \nu}
\text{\bf u} +\widetilde B(x)(\text{\bf F}+\int_0^{x_0}\widetilde L_{\widetilde \lambda,
\widetilde\mu}(x, \widetilde x_0, D')\mbox{\bf u}(\widetilde x_0, x')
d\widetilde x_0),
\end{equation}
where $\widetilde A(x, D')$ is a first order differential operator on $\partial\Omega$ and $\widetilde B$ is
a $C^1-$ matrix function. From (\ref{eb1}) and (\ref{zina}) we have
\begin{eqnarray}\label{ep}
\Vert \partial^2_{\vec \nu} \text{\bf u}
e^{s\varphi}\Vert_{H^{\frac 14,\frac 12, \widetilde s}(\Sigma)}
 \le C_{30}(\Vert  \widetilde A(x, D')\partial_{\vec \nu}
\text{\bf u}\Vert_{H^{\frac 14,\frac 12, \widetilde s}(\Sigma)}+\Vert \F e^{s\varphi}\Vert_{H^{\frac 14,\frac 12, \widetilde s}(\Sigma)}\nonumber\\
 +\Vert e^{s\varphi} B\int_0^{x_0}\widetilde L_{\widetilde \lambda,
\widetilde\mu}(x, \widetilde x_0, D')\mbox{\bf u}(\widetilde x_0, x')
d\widetilde x_0)\Vert_{H^{\frac 14,\frac 12, \widetilde s}(\Sigma)}) \le C_{31}(\Vert  (d \omega_{\mbox {\bf u}},\mbox{div}\,\mbox {\bf u})
e^{s\varphi}\Vert_{H^{\frac 34,\frac 32, \widetilde s}(\Sigma)}
\nonumber\\+\sum_{j=1}^n\Vert e^{s\varphi} \int_0^{x_0}p_j(\widetilde x_0,x)\partial^2_{\vec \nu}u_j(\widetilde x_0, x')
d\widetilde x_0)\Vert_{H^{\frac 14,\frac 12, \widetilde s}(\Sigma)}+\sum_{j=1}^n\Vert e^{s\varphi} \int_0^{x_0}\widetilde p_j(\widetilde x_0,x)\partial_{\vec \nu} u_j(\widetilde x_0, x')
d\widetilde x_0)\Vert_{H^{\frac 34,\frac 32, \widetilde s}(\Sigma)}\nonumber\\
+\Vert \F e^{s\varphi}\Vert_{H^{\frac 14,\frac 12, \widetilde s}(\Sigma)}).
\end{eqnarray} In order to get the last inequality we used (\ref{zina}).
By (\ref{granata}) for any $s\ge 1$ we have
\begin{eqnarray}\label{eep}
\sum_{j=1}^n\Vert e^{s\varphi} \int_0^{x_0}p_j(\widetilde x_0,x)\partial^2_{\vec \nu}u_j(\widetilde x_0, x')
d\widetilde x_0)\Vert_{H^{\frac 14,\frac 12, \widetilde s}(\Sigma)}+\sum_{j=1}^n\Vert e^{s\varphi} \int_0^{x_0}\widetilde p_j(\widetilde x_0,x)\partial^2_{\vec \nu} u_j(\widetilde x_0, x')
d\widetilde x_0)\Vert_{H^{\frac 14,\frac 32, \widetilde s}(\Sigma)}\nonumber\\
\le \frac{C_{32}}{s}(\Vert \partial^2_{\vec \nu} \text{\bf u}
e^{s\varphi}\Vert_{H^{\frac 14,\frac 12, \widetilde s}(\Sigma)}+\Vert \partial_{\vec \nu} \text{\bf u}
e^{s\varphi}\Vert_{H^{\frac 34,\frac 32, \widetilde s}(\Sigma)}).
\end{eqnarray}  Form (\ref{ep}), (\ref{eep}) and (\ref{zina}) we obtain (\ref{govno1}).
The proof of Proposition \ref{govno} is complete. $\blacksquare$

Next we prove

\begin{proposition}\label{golos} Let $\mbox{\bf u}\in H^1(Q), \mbox{\bf u}\vert_{\Sigma}
=0.$
There exists   $s_0>1$  such that
for all
$s\ge s_0$
\begin{eqnarray}\label{3.2'}
\sum_{\vert\alpha\vert\le 2}\int_Q (s\widetilde\varphi)^{4-2\vert\alpha\vert}\vert\partial^\alpha \mbox{\bf u}\vert^2
e^{2s\varphi}dx
\le C_{33}(\Vert (s\widetilde\varphi)^\frac 12 \nabla' d\omega_{\mbox{\bf u}}\,
e^{s\varphi}\Vert^2_{L^2(Q)} +\Vert (s\widetilde\varphi)^{-\frac 12}(\partial_{x_0}\mbox{div}\,\mbox{\bf u})
e^{s\varphi}\Vert^2_{L^2(Q)}\nonumber\\+\Vert (s\widetilde\varphi)^{-\frac 12} \partial_{x_0} d\omega_{\mbox{\bf u}}\,
e^{s\varphi}\Vert^2_{L^2(Q)}+
\Vert (s\widetilde\varphi)^\frac 12(\nabla'\mbox{div}\,\mbox{\bf u})
e^{s\varphi}\Vert^2_{L^2(Q)}\nonumber\\
+s^2\Vert\widetilde\varphi \partial_{x_0}\partial_{\vec \nu} \text{\bf u}
e^{s\varphi}\Vert_{ L^2(\widetilde\Sigma)}^2
+ s^2\Vert \widetilde\varphi\partial_{\vec \nu}  \text{\bf u}e^{s\varphi}
\Vert^2_{{ L}^2(\widetilde\Sigma)}),
\end{eqnarray}
where $C_{32}$ is independent of $s$.
\end{proposition}

{\bf Proof.} Let $x_0\in (-T,T)$ be an arbitrary but fixed.
 For any index $\widehat j\in\{1,\dots, n\}$ the short computations imply
\begin{eqnarray}\label{gavnoed}
-\Delta u_{\widehat j} =\sum_{j=1, j\ne \widehat j}^n
-\partial_{x_j}(\partial_{x_j} u_{\widehat j}-\partial_{x_{\widehat j}} u_j)
-\partial^2_{x_{\widehat j}} u_{\widehat i} -\sum_{j=1,j\ne \widehat j}
\partial_{x_j}\partial_{x_{\widehat j} }u_j\nonumber\\
=-\sum_{j=1, j\ne \widehat j}\partial_{x_j}(\partial_{x_j} u_{\widehat j}
-\partial_{x_{\widehat j}} u_j)-\partial_{x_{\widehat j} }\text{div}\,
\mbox{\bf u }\quad \text{in}\,\,\Omega.
\end{eqnarray}
Then the Carleman estimate with boundary  for the Laplace operator implies
\begin{eqnarray}\label{R}
\int_Q(\sum_{j,k=1}^n \vert \partial^2_{x_kx_j}\mbox{\bf u}\vert^2+ s^2\widetilde\varphi^2 \vert \nabla '\mbox{\bf u}\vert^2+s^4\widetilde\varphi^4
\vert \mbox{\bf u}\vert^2)e^{2s\varphi}dx\le C_{34}(\Vert s^\frac 12\widetilde\varphi^
\frac 12\nabla ' \mbox{div}\, \mbox{\bf u}e^{s\varphi}\Vert^2_{L^2(Q)}
\nonumber\\
+\Vert s^\frac 12\widetilde\varphi^\frac 12 \nabla' d\omega_{\mbox{\bf u}}\,
e^{s\varphi}\Vert^2_{L^2(Q)}+\int_{\widetilde \Sigma}s^2\widetilde\varphi^2\vert
\partial_{\vec\nu}\mbox{\bf u}\vert^2 e^{2s\varphi}d\Sigma).
\end{eqnarray}
We differentiate both sides of equation (\ref{gavnoed})
with respect to the variable $x_0$ and take $H^{-1}$ Carleman estimate
by authors in \cite{IY7}:

\begin{eqnarray}\label{RR}
\int_Q s^2\widetilde\varphi^2\vert \partial_{x_0} \mbox{\bf u}\vert^2e^{2s\varphi}
dx\le C_{35}(\Vert s^\frac 12\widetilde\varphi^\frac 12 \mbox{div}\, \partial_{x_0}
\mbox{\bf u}e^{s\varphi}\Vert^2_{L^2(Q)}\nonumber\\+\Vert s^\frac 12\widetilde
\varphi^\frac 12  \partial_{x_0} d\omega_{\mbox{\bf u}}\, e^{s\varphi}\Vert^2
_{L^2(Q)}+\int_{\widetilde \Sigma}s^2\widetilde\varphi^2\vert \partial
_{\vec\nu}\partial_{x_0}\mbox{\bf u}\vert^2 e^{2s\varphi}d\Sigma).
\end{eqnarray}
Combination of (\ref{R}) and (\ref{RR}) implies (\ref{3.2'}).
The proof of the proposition is complete.
                                    $\blacksquare$

By
Proposition \ref{govno} and Proposition \ref{golos}  from (\ref{zmpolit1}), (\ref{kom}) we obtain the estimate
\begin{eqnarray}\label{zmpolit1o}
\Vert\partial^2_{\vec \nu} \text{\bf u}
e^{s\varphi}\Vert_{H^{\frac 14,\frac 12,\widetilde s}( \Sigma_0)}^2
+\Vert \partial_{\vec \nu}
\text{\bf u}
e^{s\varphi}\Vert_{H^{\frac 34,\frac 32,\widetilde s}( \Sigma_0)}^2\\
+\int_Q \sum_{\vert\alpha\vert\le 2} (s\widetilde\varphi)^{4-2\vert\alpha\vert}\vert\partial^\alpha \mbox{\bf u}\vert^2
e^{2s\varphi}dx + \int_{Q}\sum_{\vert\alpha\vert\le 1,\alpha_0=0}(s\widetilde\varphi)^{3-2\vert\alpha\vert}\vert\partial^\alpha (d\omega_{\mbox{\bf u}},\mbox{div}\,\text{\bf u})\vert^2
e^{2s\varphi} dx\nonumber\\
\le C_{36}\biggl(\Vert
\mbox{\bf F}e^{s\varphi}\Vert^2_{H^{\frac 14,\frac 12,\widetilde s}({\Sigma})}\nonumber+\Vert  (\mbox{div}\, \mbox{\bf F})e^{s\varphi}\Vert^2_{L^2(Q)}+\Vert \mbox{\bf F}e^{s\varphi}\Vert^2_{L^2(Q)}\\ +\Vert (d \omega_{\mbox {\bf u}},\mbox{div}\,\mbox {\bf u})
e^{s\varphi}\Vert_{H^{\frac 34,\frac 32, \widetilde s}(\widetilde\Sigma)}^2+\Vert \partial_{\vec \nu}(d \omega_{\mbox {\bf u}},\mbox{div}\,\mbox {\bf u})
e^{s\varphi}\Vert_{H^{\frac 14,\frac 12, \widetilde s}(\widetilde\Sigma)}^2+\Vert\widetilde\varphi \partial_{x_0}\partial_{\vec \nu} \text{\bf u}
e^{s\varphi}\Vert_{ L^2(\widetilde\Sigma)}^2)
        \nonumber
\end{eqnarray}
for all $s\ge s_0.$  By (\ref{nina}) we obtain from (\ref{zmpolit1o}) the estimate (\ref{2.9'}).
Thus the proof of Theorem \ref{opa3} is finished. $\blacksquare$

\section{Proof of Theorem \ref{theorem 1.2}.}\label{Q1}

We differentiate equations (\ref{(A1.1)})-(\ref{(A1.3)}) respect to $x_0$:
\begin{equation}\label{(1.1!)}
\rho\ppp^3_{x_0}\mbox{\bf u}
= L_{\lambda,\mu}(x',D') \ppp^2_{x_0}\mbox{\bf u} +  \LLLL(x,D')
\ppp_{x_0}\mbox{\bf u}+L_{\partial_{x_0}\widetilde\lambda,\partial_{x_0}\widetilde\mu}(x,D')
\mbox{\bf u} + \partial_{x_0}R(x)\f(x') \quad
\mbox{in $(0,T)\times \OOO$}.
\end{equation}
From (\ref{(A1.1)}) and (\ref{(A1.222)}) we obtain
\begin{equation}\label{(1.222!)}
\ppp^2_{x_0}\mbox{\bf u}(\eta,\cdot) =\frac{1}{\rho}( \LLL(x',D')\mathbf{b}+L_{\widetilde\mu(\eta,\cdot),\widetilde\lambda(\eta,\cdot)}(\eta,x',D')\mathbf{a}+R(\eta,\cdot)\mbox{\bf f}) , \quad\ppp_{x_0} \mbox{\bf u}(\eta,\cdot) = \mbox{\bf b},
\end{equation}

\begin{equation}\label{(1.3!)}
\ppp_{x_0}\mbox{\bf u}\vert_{(0,T)\times \ppp\OOO} = 0.
\end{equation}

We set
$\y = \ppp^2_{x_0}\uu$.  Then
$$
\partial_{x_0}\uu(x_0,x') = \int^{x_0}_{\eta} \y(\wwwx,x') d\wwwx + \mathbf{b}(x').
$$

Using this equality we rewrite (\ref{(1.1!)})-(\ref{(1.3!)}) in terms of the unknown function $\y:$
\begin{eqnarray}
 \rho\ppp_{x_0} \y = \LLL(x',D') \y
+ \LLLL(x,D')\int^{x_0}_{\eta} \y(\wwwx,x') d\wwwx +\mbox{\bf F},
\quad   \mbox{in $\Omega\times (0,T)$}, \nonumber\\
 \y(\eta,\cdot) = \frac{1}{\rho}( \LLL(x',D')\mathbf{b}+L_{\widetilde\mu(\eta,\cdot),\widetilde\lambda(\eta,\cdot)}(\eta,x',D')\mathbf{a}+R(\eta,\cdot)\mbox{\bf f}) \quad \mbox{in $\Omega$}, \quad
 \y\vert_{\ppp\Omega\times (0,T)} = 0, \nonumber
\end{eqnarray} where
$\mbox{\bf F}(x)=L_{\partial_{x_0}\widetilde\lambda,\partial_{x_0}\widetilde\mu}(x,D')
\mbox{\bf u}+\LLLL(x',D') \mathbf{b}(x')
+ \partial_{x_0} R(x)f(x').
$

Let $\widetilde \y(x)=\y (x_0+\eta,x'),\widetilde {\mbox{\bf u}}(x)={\mbox{\bf u}} (x_0+\eta,x'), \widetilde {\mbox{\bf F}}=\mbox{\bf F}(x_0+\eta,x').$ Then
\begin{equation}\label{viri}
 \rho\ppp_{x_0}\widetilde \y = \LLL(x',D')\widetilde \y
+ \LLLL(x',D')\int^{x_0}_{0} \widetilde\y(\wwwx,x') d\wwwx +\widetilde{\mbox{\bf F}},
\quad  \mbox{in $\Omega\times (-\eta,T-\eta)$}
\end{equation}
\begin{equation}\label{viri2}
\widetilde\y(0,\cdot) =\frac{1}{\rho}( \LLL(x',D')\mathbf{b}+L_{\widetilde\mu(\eta,\cdot),\widetilde\lambda(\eta,\cdot)}(\eta,x',D')\mathbf{a}+R(\eta,x')\mbox{\bf f}) \quad \mbox{in $\Omega$}, \quad
 \widetilde\y\vert_{\ppp\Omega\times (-\eta,T-\eta)} = 0.
\end{equation}

Denote $Q_{\hat T}=(-\hat T,\hat T)\times \Omega,  \widetilde\Sigma_{\hat T}=(-\hat T,\hat T)\times \widetilde\Gamma$ where $\hat T=\frac 12\mbox{max}\{\eta, T-\eta\}.$ Let $\varphi(x)=\frac{e^{\lambda \psi}-e^{2\lambda\Vert\psi\Vert_{C^0(\overline\Omega)}}}{(\hat
T-x_0)^3(\hat T+x_0)^3} $ where function $\psi$ is constructed in \cite{IM}, $\lambda$ is sufficiently large positive parameter.  We apply Carleman estimate (\ref{2.9'}) to the system (\ref{viri}), (\ref{viri2}):
\begin{eqnarray}\label{2.9'l}
\Vert \widetilde \y\Vert_{\mathcal B(\varphi,s,Q_{\hat T})}
\le C_1(\Vert \text{\bf F}
e^{s\varphi}\Vert_{\mathcal Y(\varphi,s,Q_{\hat T})}+ \Vert (d \omega_{\widetilde{\mbox {\bf y}}},\mbox{div}\,\widetilde{\mbox {\bf y}})
e^{s\varphi}\Vert_{H^{\frac 34,\frac 32,\widetilde s}(\widetilde\Sigma_{\hat T})}\\+\Vert \partial_{\vec \nu}(d \omega_{\widetilde{\mbox {\bf y}}},\mbox{div}\,\widetilde{\mbox {\bf y}})
e^{s\varphi}\Vert_{H^{\frac 14,\frac 12,\widetilde s}(\widetilde\Sigma_{\hat T})}+\Vert\widetilde\varphi \partial_{x_0}\partial_{\vec \nu} \widetilde{\text{\bf y}}
e^{s\varphi}\Vert_{ L^2(\widetilde\Sigma_{\hat T})})\quad\forall s\ge s_0.
\nonumber
\end{eqnarray}
By the stationary phase  argument (see e.g. \cite{Stein}) for all $s\ge 1$
\begin{eqnarray}\label{nokia4}
\Vert \text{\bf F}
e^{s\varphi}\Vert_{\mathcal Y(\varphi,s,Q_{\hat T})}\le\frac{ C_2}{\root\of s}(\sum_{\vert \alpha\vert\le 2}\Vert \partial^\alpha\mathbf{b} e^{s\varphi(0,\cdot)}\Vert_{H^{1,s}(\Omega)}\\+\Vert  (d \omega_{\widetilde{\mbox {\bf u}}},\mbox{div}\,\widetilde{\mbox {\bf u}}) e^{s\varphi}\Vert_{L^2(-\hat T,\hat T;H^{2,s}(\Omega))}+\Vert\widetilde{\mbox {\bf u}} e^{s\varphi}\Vert_{L^2(-\hat T,\hat T;H^{2,s}(\Omega))}+ \Vert \f e^{s\varphi(0,\cdot)}\Vert_{H^{1,s}(\Omega)} ) \nonumber.
\end{eqnarray}

On the other hand by (\ref{viri2}) and (\ref{(1.18)}) for all positive $s$ we have
\begin{eqnarray}\label{nokia}
\Vert \f e^{s\varphi(0,\cdot)}\Vert_{H^{1,s}(\Omega)}\le C_3 \Vert R(\eta,\cdot)\f e^{s\varphi(0,\cdot)}\Vert_{H^{1,s}(\Omega)}\nonumber \\
\le C_4(\Vert \widetilde\y(0,\cdot)e^{s\varphi(0,\cdot)}\Vert_{H^{1,s}(\Omega)}+\Vert \mathbf{b}e^{s\varphi(0,\cdot)}\Vert_{H^{3, s}(\Omega)}+\Vert \mathbf{a}e^{s\varphi(0,\cdot)}\Vert_{H^{3,s}(\Omega)}).
\end{eqnarray}
Observe that by (\ref{gopak}) there exists a constant $C_5$ independent  of $s$ such that
\begin{equation}\label{nokia1}
\Vert \widetilde\y(0,\cdot)e^{s\varphi(\cdot,0)}\Vert_{H^{1,s}(\Omega)}
\le C_5
\Vert \widetilde \y\Vert_{\mathcal B(\varphi,s,Q_{\hat T})}\quad\forall s\ge 1.
\end{equation}
Using (\ref{nokia1}) to estimate the  first term in the right hand of (\ref{nokia}) and applying (\ref{2.9'l}) we obtain
\begin{eqnarray}\label{nokia3}
\Vert \f e^{s\varphi(0,\cdot)}\Vert_{H^{1,s}(\Omega)}\le C_6 \Vert R(\eta,\cdot)\f e^{s\varphi(0,\cdot)}\Vert_{H^{1,s}(\Omega)}\nonumber \\
\le C_7(\Vert \mathbf{b}e^{s\varphi(0,\cdot)}\Vert_{H^{3, s}(\Omega)}+\Vert \mathbf{a}e^{s\varphi(0,\cdot)}\Vert_{H^{3,s}(\Omega)}+\Vert \widetilde \y\Vert_{\mathcal B(\varphi,s,Q_{\hat T})})\nonumber\\
\le
C_8(\Vert \mathbf{b}e^{s\varphi(0,\cdot)}\Vert_{H^{3, s}(\Omega)}+\Vert \mathbf{a}e^{s\varphi(0,\cdot)}\Vert_{H^{3, s}(\Omega)}\nonumber\\
+\Vert \text{\bf F}
e^{s\varphi}\Vert_{\mathcal Y(\varphi,s,Q_{\hat T})}+ \Vert (d \omega_{\widetilde{\mbox {\bf y}}},\mbox{div}\,\widetilde{\mbox {\bf y}})
e^{s\varphi}\Vert_{H^{\frac 34,\frac 32, \widetilde s}(\widetilde\Sigma_{\hat T})}\nonumber\\
+\Vert \partial_{\vec \nu}(d \omega_{\widetilde{\mbox {\bf y}}},\mbox{div}\,\widetilde{\mbox {\bf y}})
e^{s\varphi}\Vert_{H^{\frac 14,\frac 12, \widetilde s}(\widetilde\Sigma)}+\Vert\widetilde\varphi \partial_{x_0}\partial_{\vec \nu} \widetilde{\text{\bf y}}
e^{s\varphi}\Vert_{ L^2(\widetilde\Sigma_{\hat T})})\quad\forall s\ge s_1.
\end{eqnarray}
From (\ref{nokia3}), (\ref{nokia4}) we have
\begin{eqnarray}\label{nokia5}
\Vert \f e^{s\varphi(0,\cdot)}\Vert_{H^{1,s}(\Omega)}
\le
C_9(\Vert \mathbf{b}e^{s\varphi(\cdot,0)}\Vert_{H^{3,\widetilde s}(\Omega)}+\Vert \mathbf{a}e^{s\varphi(0,\cdot)}\Vert_{H^{3, s}(\Omega)}\nonumber\\
+\frac{1}{\root\of{s}}(\sum_{\vert \alpha\vert\le 2}\Vert \partial^\alpha\mathbf{b}e^{s\varphi(0,\cdot)}\Vert_{H^{1,s}(\Omega)}+ \Vert \f e^{s\varphi(0,\cdot)}\Vert_{H^{1,s}(\Omega)}\nonumber\\+\Vert  (d \omega_{\widetilde{\mbox {\bf u}}},\mbox{div}\,\widetilde{\mbox {\bf u}}) e^{s\varphi}\Vert_{L^2(-\hat T,\hat T;H^{2,s}(\Omega))}+\Vert\widetilde{\mbox {\bf u}} e^{s\varphi}\Vert_{L^2(-\hat T,\hat T;H^{2,s}(\Omega))}  )\nonumber\\+ \Vert (d \omega_{\widetilde{\mbox {\bf y}}},\mbox{div}\,\widetilde{\mbox {\bf y}})
e^{s\varphi}\Vert_{H^{\frac 34,\frac 32, \widetilde s}(\widetilde\Sigma_{\hat T})}\nonumber\\
+\Vert \partial_{\vec \nu}(d \omega_{\widetilde{\mbox {\bf y}}},\mbox{div}\,\widetilde{\mbox {\bf y}})
e^{s\varphi}\Vert_{H^{\frac 14,\frac 12, \widetilde s}(\widetilde\Sigma)}+\Vert\widetilde\varphi \partial_{x_0}\partial_{\vec \nu} \widetilde{\text{\bf y}}
e^{s\varphi}\Vert_{ L^2(\widetilde\Sigma_{\hat T})})\quad\forall s\ge s_1.
\end{eqnarray}
 We introduce the functions $\varphi_*(x)=\frac{e^{\lambda \psi}-e^{2\lambda\Vert\psi\Vert_{C^0(\overline\Omega)}}}{\ell_*(x_0)},\widetilde\varphi_*(x_0)=\frac{1}{(2\hat T-x_0)^3(2\hat T+x_0)^3} $ where $\ell_*\in C^3[-2\hat T,2\hat T]$ , $\ell_*(\pm 2\hat T)=\ell_*'(\pm 2\hat T)=0$, $\ell_*(x_0)=(\hat
T-x_0)^3(\hat T+x_0)^3$ on $[-\hat T/4,\hat T/4]$ and $\ell_*(x_0)\ge (\hat
T-x_0)^3(\hat T+x_0)^3$ on $[-\hat T,\hat T].$
Therefore
\begin{equation}\label{gosha}
\varphi_*(x)= \varphi(x)\quad \forall x\in Q_{\hat T}.
\end{equation}
To estimate the norm of  $\partial^\alpha\widetilde{\mbox {\bf u}}$ in the right hand side of  (\ref{nokia5}) we apply to equations for $\widetilde{\mbox{\bf u}}$   Carleman estimate (\ref{2.9'}) to the system  (\ref{(A1.1)})-(\ref{(A1.3)}) and using the stationary phase argument  we obtain:
\begin{eqnarray}\label{2.9'lA}
\Vert \widetilde {\mbox {\bf u}}\Vert_{\mathcal B(\varphi_*,s,Q_{2\hat T})}
\le C_{10}(\Vert R\text{\bf f}
e^{s\varphi_*}\Vert_{\mathcal Y(\varphi_*,s,Q_{2\hat T})}+ \Vert (d \omega_{\widetilde{\mbox {\bf u}}},\mbox{div}\,\widetilde{\mbox {\bf u}})
e^{s\varphi_*}\Vert_{H^{\frac 34,\frac 32,\widetilde s}(\widetilde\Sigma_{2\hat T})}\\+\Vert \partial_{\vec \nu}(d \omega_{\widetilde{\mbox {\bf u}}},\mbox{div}\,\widetilde{\mbox {\bf u}})
e^{s\varphi_*}\Vert_{H^{\frac 14,\frac 12,\widetilde s}(\widetilde\Sigma_{2\hat T})}+\Vert\widetilde\varphi_* \partial_{x_0}\partial_{\vec \nu} \widetilde{\text{\bf u}}
e^{s\varphi_*}\Vert_{ L^2(\widetilde\Sigma_{2\hat T})})\nonumber\\
\le C_{11}(s^{-\frac 12}\Vert \text{\bf f}
e^{s\varphi(0,\cdot)}\Vert_{H^{1,s}(\Omega)}+ \Vert (d \omega_{\widetilde{\mbox {\bf u}}},\mbox{div}\,\widetilde{\mbox {\bf u}})
e^{s\varphi_*}\Vert_{H^{\frac 34,\frac 32,\widetilde s}(\widetilde\Sigma_{2\hat T})}\nonumber\\+\Vert \partial_{\vec \nu}(d \omega_{\widetilde{\mbox {\bf u}}},\mbox{div}\,\widetilde{\mbox {\bf u}}) \nonumber
e^{s\varphi_*}\Vert_{H^{\frac 14,\frac 12,\widetilde s}(\widetilde\Sigma_{2\hat T})}+\Vert\widetilde\varphi_* \partial_{x_0}\partial_{\vec \nu} \widetilde{\text{\bf u}}
e^{s\varphi_*}\Vert_{ L^2(\widetilde\Sigma_{2\hat T})})\quad\forall s\ge s_2.
\nonumber
\end{eqnarray}
By (\ref{gopak})
\begin{equation}\label{grob}
\left (\int_{Q_{2\hat T}} \sum_{\vert \alpha\vert\le 2,\alpha_0=0}\left (\frac{1}{s\widetilde \varphi_*}\vert\partial^\alpha (d \omega_{\widetilde{\mbox {\bf u}}},\mbox{div}\,\widetilde{\mbox {\bf u}})\vert^2+\vert \partial^\alpha \widetilde{\mbox {\bf u}}\vert^2\right )e^{2s\varphi_*}dx\right )^\frac 12\le C_{12}\Vert \widetilde {\mbox {\bf u}}\Vert_{\mathcal B(\varphi_*,s,Q_{2\hat T})}.
\end{equation}
By (\ref{2.9'l}) and (\ref{grob})
\begin{eqnarray}\label{PPP}
\left (\int_{Q_{2\hat T}} \sum_{\vert \alpha\vert\le 2,\alpha_0=0}\left (\frac{1}{s\widetilde \varphi_*}\vert\partial^\alpha (d \omega_{\widetilde{\mbox {\bf u}}},\mbox{div}\,\widetilde{\mbox {\bf u}})\vert^2+\vert \partial^\alpha \widetilde{\mbox {\bf u}}\vert^2\right ) e^{2s\varphi_*}dx\right )^\frac 12\\
\le C_{13}(s^{-\frac 12}\Vert \text{\bf f}
e^{s\varphi_*(0,\cdot)}\Vert_{H^{1,s}(\Omega)}+ \Vert (d \omega_{\widetilde{\mbox {\bf u}}},\mbox{div}\,\widetilde{\mbox {\bf u}})
e^{s\varphi_*}\Vert_{H^{\frac 34,\frac 32,\widetilde s}(\widetilde\Sigma_{2\hat T})}\nonumber\\+\Vert \partial_{\vec \nu}(d \omega_{\widetilde{\mbox {\bf u}}},\mbox{div}\,\widetilde{\mbox {\bf u}}) \nonumber
e^{s\varphi_*}\Vert_{H^{\frac 14,\frac 12,\widetilde s}(\widetilde\Sigma_{2\hat T})}+\Vert\widetilde\varphi_* \partial_{x_0}\partial_{\vec \nu} \widetilde{\text{\bf u}}
e^{s\varphi_*}\Vert_{ L^2(\widetilde\Sigma_{2\hat T})})\quad\forall s\ge s_2.
\nonumber
\end{eqnarray}
By (\ref{gosha}) and
estimates (\ref{nokia5}) and (\ref{PPP})  for sufficiently large $s$ implies (\ref{(1.19)}).
 The proof of Theorem \ref{theorem 1.2} is completed. $\blacksquare$

 {\bf Proof of Theorem \ref{theorem 1.3}.} Let functions $(\rho_j, \mbox{\bf v}_j), j\in \{1,2\}$ satisfy the equations
 \begin{equation}\label{Z1}
\partial_{x_0}\rho_j+\mbox{div}(\mbox{\bf v}_j\rho_j)=0\quad \mbox{in} \,Q,
\end{equation}
\begin{equation}\label{Z2}
\rho_j\partial_{x_0}\mbox{\bf v}_j-L_{\lambda,\mu}(x',D')\mbox{\bf v}_j +\rho_j(\mbox{\bf v}_j,\nabla')\mbox{\bf v}_j+h(\rho_j)\nabla' \rho_j=R\mbox{\bf f}_j\quad \mbox{in} \,Q,
\end{equation}
\begin{equation}\label{Z3}
\mbox{\bf v}_j\vert_{\Sigma}=0.
\end{equation}
We set $\rho=\rho
_1-\rho_1,\mbox{\bf v}=\mbox{\bf v}_1-\mbox{\bf v}_2, \mbox{\bf f}=\mbox{\bf f}_1-\mbox{\bf f}_2.$ Then, from equations (\ref{Z1})-(\ref{Z3}) we have
\begin{equation}\label{ZZ1}
\partial_{x_0}\rho+\mbox{div}(\mbox{\bf v}_1\rho)=-\mbox{div}(\mbox{\bf v}\rho_2)\quad \mbox{in} \,Q,
\end{equation}
\begin{eqnarray}\label{ZZ2}
\rho_1\partial_{x_0}\mbox{\bf v}+\rho\partial_{x_0}\mbox{\bf v}_2-L_{\lambda,\mu}(x',D')\mbox{\bf v}+\rho(\mbox{\bf v}_1,\nabla')\mbox{\bf v}_1 +\rho_2(\mbox{\bf v}_1,\nabla')\mbox{\bf v}+\rho_2(\mbox{\bf v},\nabla')\mbox{\bf v}_2\nonumber\\+h(\rho_1)\nabla' \rho+(h(\rho_1)-h(\rho_2))\nabla' \rho_2=R\mbox{\bf f}\quad \mbox{in} \,Q,
\end{eqnarray}
\begin{equation}\label{ZZ3}
\mbox{\bf v}\vert_{\Sigma}=0.
\end{equation}

%Let $B(0,\hat r)=\{x '\in \Bbb R^n\vert \vert x\vert <\hat r\}$ be a ball in $\Bbb R^n$ of sufficiently large radius $\hat r$ such that %$\Omega\subset B(0,\hat r).$ We extend the function $\mbox{\bf v}_1$ up to $C^3([-T,T]\times  B(0,\hat r))$ function such that %$\mbox{supp}\,\mbox{\bf v}_1\subset [-T,T]\times B(0,\hat r).$
Let $y(t)=(y_1(t),\dots ,y_n(t))$ be a solution to the system of ordinary differential equations
$$
\frac{d y}{dt}=\mbox{\bf v}_1(t, y(t)).
$$ The curve $(t,y(t))$ is the characteristic curve  for the hyperbolic operator
$$
L_{\mbox{\bf v}_1}(x,D)w=\partial_{x_0}w+(\mbox{\bf v}_1(x),\nabla' w).
$$
Consider the initial value problem
\begin{equation}\label{oop}
L_{\mbox{\bf v}_1}(x,D)\psi=0\quad\mbox{on}\,\, Q,\quad \psi(0,x')=\psi_0(x')\quad x'\in \Omega
\end{equation}
such that the initial condition $\psi_0\in C^3(\overline\Omega)$  satisfies
\begin{equation}\label{pop}
\partial_{\vec \nu}\psi_0\vert_{\partial\Omega\setminus\widetilde \Gamma}<-C_{14}<0,\quad \vert \nabla' \psi_0(x)\vert>C_{15}>0\quad\forall x\in \Omega .
\end{equation}
 Existence of such  a function $\psi_0$  is proved in \cite{IM}. Thanks to (\ref{ZZ3})
solution $\psi \in C^3(\overline Q)$ to the Cauchy problem (\ref{oop}) exist, unique and satisfies the estimate
\begin{equation}\label{gora1}
\Vert \psi\Vert_{C^3(\overline Q)}\le C_{16}\Vert \psi_0\Vert_{C^3(\overline\Omega)}.
\end{equation}   Then, by (\ref{gora1}) and (\ref{pop}) there exist a positive $\delta$ such that

%\begin{equation}\label{begemot}\Vert \psi\Vert_{C^3(\overline\Omega)}\le C_{17},
%\end{equation}
\begin{equation}\label{zim}
 \partial_{\vec \nu}\psi\vert_{[-\delta,\delta]\times (\partial\Omega\setminus\widetilde \Gamma)}<-C_{17}<0,\quad \vert \nabla' \psi(x)\vert>C_{18}>0\quad\forall x\in Q_\delta\triangleq (-\delta,\delta)\times \Omega,
\end{equation}
\begin{equation}\label{coppp}
-\partial_{\vec \nu}\psi(x)>\frac {1}{\root\of 2}\,\root\of{\frac{\mu(x)}{(\lambda+2\mu)(x)}} \vert \partial_{\vec\tau}\psi(x)\vert \quad \forall x\in [-\delta,\delta]\times \overline \Gamma_0, \quad \forall \vec\tau\in T(\Gamma_0), \,\, \vert \vec\tau\vert=1.
\end{equation}

 Using such a function $\psi $ we construct the weight function $\varphi$ by formula
\begin{equation}\label{gold}\varphi(x)=\frac{e^{\lambda \psi(x)}-e^{2\lambda\Vert\psi\Vert_{C^0(\overline Q_\delta)}}}{{\widetilde\ell}(x_0)},\quad \widetilde \varphi(x_0)=\frac{1}{(\delta-x_0)^3(\delta+x_0)^3},
\end{equation}
where $\lambda$ is the large positive parameter and $\delta$ is a small positive parameter, ${\widetilde\ell}(t)=1-\vert x_0\vert$ on $[-\delta/2,\delta/2],$
$\partial_{x_0}{\widetilde\ell}(x_0) <0 $ on $[0,\delta]$ and $\partial_{x_0}{\widetilde\ell}(x_0) >0 $ on $[-\delta,0]$, ${\widetilde\ell}(x_0)>0$ on $(-\delta,\delta)$, $\partial_{x_0}^k{\widetilde\ell}(\pm\delta)=0$ for all $k\in \{0,1,2\},$ $\partial_{x_0}^3{\widetilde\ell}(\pm\delta)\ne 0,$ ${\widetilde\ell}\in W^1_\infty[-\delta,\delta]\cap C^2[0,\delta]\cap C^2[-\delta,0].$

%$\quad {\widetilde\ell}(x_0)\in C^\infty[-\delta,\delta],$ be a function strictly positive on $(-\delta,\delta)$  such that $\partial_{x_0}^2{\widetilde\ell}(0)< 0, %\partial_{x_0}{\widetilde\ell}(t)<0 \,\,\mbox{on}\, \,\,(0,\delta], $ \newline$ \partial_{x_0}{\widetilde\ell}(t)>0 \,\,\mbox{on}\,\, [-\delta,0), %{\widetilde\ell}(t)=(\delta^2-x_0^2)^3\,\, \mbox{on}\,\,([-\delta,-\delta/2]\cap [\delta/2,\delta])\times\Omega.
%$

We have

\begin{proposition}Let $\lambda>1$.  There exist a positive $\delta_0$  independent of $\lambda$ such that for all $\delta\in (0,\delta_0)$ function $\varphi$ given by (\ref{gold}) satisfies (\ref{pinok220}), (\ref{rep}), (\ref{giga1})
- (\ref{lomka}), Condition \ref{A1} and Condition \ref{A2} with $\beta = \mu$ and $\beta=\lambda+2\mu$. Moreover  for some positive constants $C_{19}, C_{20}$
\begin{equation}\label{vobla}
\frac{d\varphi(t,y(t))}{dt} <-C_{19}< 0 \quad \mbox{for}\,\, t\in [0,\delta)\quad \mbox{and}\quad \frac{d\varphi(t,y(t))}{dt}> C_{20}>0 \quad \mbox{for}\,\, t\in (-\delta,0].
\end{equation}
%and
%\begin{equation}\label{vobla1}
%\frac{d^2\varphi(t,y(t))}{dt^2}\vert_{t=0} < 0.
%\end{equation}
Moreover there exists a positive $\lambda_0$  and positive  $C_{21}$ such that  for all $\lambda>\lambda_0$
\begin{equation}\label{victoryA}
\partial_{x_0}\varphi(x)< -C_{21} <0\quad \mbox{on}\, \, [0,\delta)\times \overline\Omega, \quad \partial_{x_0}\varphi(x)>C_{21} >0\quad \mbox{on}\, \, (-\delta,0]\times \overline\Omega.
\end{equation}
\end{proposition}

{\bf Proof.} Formula (\ref{gold}), equation (\ref{oop}) and short computations imply
$$
\frac{d\varphi(t,y(t))}{dt}=\frac{d }{dt}\left(\frac{e^{\lambda \psi(t,y(t))}-e^{2\lambda\Vert\psi\Vert_{C^0(\overline Q_\delta)}}}{{\widetilde\ell}(t)} \right)=
$$
$$
\frac{e^{\lambda \psi(t,y(t))}(L_{\mbox{\bf v}_1}(x,D)\psi)(t,y(t))}{{\widetilde\ell}(t)}-\frac{e^{\lambda \psi(t,y(t))}-e^{2\lambda\Vert\psi\Vert_{C^0(\overline Q_\delta)}}}{{\widetilde\ell}^2(t)}{\widetilde\ell}'(t)=
$$
$$
-\frac{e^{\lambda \psi(t,y(t))}-e^{2\lambda\Vert\psi\Vert_{C^0(\overline Q_\delta)}}}{{\widetilde\ell}^2(t)}{\widetilde\ell}'(t).
$$
Let $\lambda$ be sufficiently large, since  ${\widetilde\ell}'(t)<0$ on $[0,\delta]$ and  ${\widetilde\ell}'(t)>0$ on $[-\delta,0]$  from the above formula we have (\ref{vobla}).
%Taking the second derivative of the function $\varphi(t,y(t))$ we have
%$$
%\frac{d^2\varphi(t,y(t))}{dt^2}=-\frac{d }{dt}\left( \frac{e^{\lambda \psi(t,y(t))}-e^{2\lambda\Vert\psi\Vert_{C^0(\overline %Q_\delta)}}}{{\widetilde\ell}^2(t)}{\widetilde\ell}'(t) \right)=
%$$
%$$ -\frac{e^{\lambda \psi(t,y(t))}-e^{2\lambda\Vert\psi\Vert_{C^0(\overline Q_\delta)}}}{{\widetilde\ell}^2(t)}{\widetilde\ell}''(t) +2 \frac{e^{\lambda %\psi(t,y(t))}-e^{2\lambda\Vert\psi\Vert_{C^0(\overline Q_\delta)}}}{{\widetilde\ell}^3(t)}({\widetilde\ell}'(t))^2.
%$$
%Since ${\widetilde\ell}'(0)=0$ this formula implies (\ref{vobla1}).
%Since $\nabla'\psi\ne 0$ on $Q_\delta$ pseudoconvexity condition \ref{A2} holds true provided that $\lambda$ is sufficiently large.
The inequality (\ref{pinok220}) follows from (\ref{zim}) and (\ref{gold}). Inequality (\ref{coppp}) implies (\ref{rep}).
Formula (\ref{gold}) immediately implies  (\ref{1})-(\ref{lomka}).
Now we check Condition \ref{A1}.
For simplicity of notations we denote $p(x,\xi)=p_{\rho,\beta}(x,\xi)=i\rho\xi_0+\sum_{k,j=1}^na_{k,j}\xi_k\xi_j.$ We remind that $\zeta$ is given by formula (\ref{propoganda1}).
 Observe that
$$
\partial_{ x_m}{\overline
p}(x,\xi_0,\overline{\widetilde\zeta})=\overline p_{m}(x,\xi_0,\overline{\widetilde\zeta})-i\vert
s\vert\sum_{k=1}^np^{(k)}(x,\xi_0,\overline{\widetilde\zeta})\partial^2_{x_k x_m}\varphi.
$$
Then
\begin{eqnarray}
\mbox{Im}\{\overline p(x,\xi_0,\overline{\widetilde\zeta}),{p(x,\xi_0,{\widetilde\zeta})}\}=\nonumber \\
\mbox{Im}\left(\sum_{k=0}^n\overline p^{(k)}(x,\xi_0,\overline{\widetilde\zeta})
\partial_{x_k}p(x,\xi_0,\widetilde\zeta)-\partial_{x_k}\overline p(x,\xi_0,\overline{\widetilde\zeta})
p^{(k)}(x,\xi_0,\widetilde\zeta)\right).\nonumber
\end{eqnarray}
Simple computations provide the following formulas:
\begin{eqnarray}
\mbox{Im}\left(\partial_{ \xi_0}\overline
p(x,\xi_0,\overline{\widetilde\zeta}) \partial_{
x_0}p(x,\xi_0,\widetilde\zeta)-\partial_{ x_0}\overline
p(x,\xi_0,\overline{\widetilde\zeta})\partial_{\xi_0}
p(x,\xi_0,\widetilde\zeta)\right)\nonumber
\\=\mbox{Im}\left((-i\rho)(p_{0}(x,\xi_0,\widetilde\zeta)+i\vert s\vert\sum_{m=1}^n p^{(m)}(x,\xi_0,\widetilde\zeta)
\partial^2_{ x_m x_0}\varphi)\right.\nonumber\\
\left. -i\rho(\overline p_{0}(x,\xi_0,\overline{\widetilde\zeta})-i\vert s\vert\sum_{m=1}^n
p^{(m)}(x,\xi_0,\overline{\widetilde\zeta})\partial^2_{x_m x_0}\varphi)\right\}=\nonumber\\=-2\rho\mbox{Re}\,p_{0}(x,\xi)+2\rho
s^2\sum_{m=1}^n
p^{(m)}(x,\nabla \varphi)\partial^2_{x_m x_0}\varphi+2s^2\rho a_{0}(x,\nabla'\varphi,\nabla'\varphi),\nonumber
\end{eqnarray}
where $a_{0}(x,\widetilde\eta,\widetilde\eta)=\sum_{k,j=1}^n\frac{\partial
a_{kj}}{\partial x_0} \eta_k\eta_j$ and for any $k\in \{1,\dots, n\}$
\begin{eqnarray}
\mbox{Im}\left(\partial_{\xi_k}\overline
p(x,\xi_0,\overline{\widetilde\zeta})\partial_{ x_k}
p(x,\xi_0,\widetilde{\zeta})-\partial_{ x_k}\overline
p(x,\xi_0,\overline{\widetilde\zeta})\partial_{\xi_k}
p(x,\xi_0,\widetilde\zeta)\right)\nonumber
\\=\mbox{Im}\left(\overline p^{(k)}(x,\xi_0,\overline{\widetilde\zeta})
(p_{k}(x,\xi_0,\widetilde\zeta)
+ i\vert s\vert\sum_{m=1}^n p^{(m)}(x,\xi_0,\widetilde\zeta)
\partial^2_{x_k x_m}\varphi)\right.\nonumber\\
-\left.
p^{(k)}(x,\xi_0,\widetilde\zeta)(\overline p_{k}(x,\xi_0,\overline{\widetilde\zeta})-i\vert
s\vert\sum_{m=1}^n p^{(m)}(x,\xi_0,\overline{\widetilde\zeta}
)\partial^2_{ x_k x_m}\varphi)\right
)=\nonumber\\-p^{(k)}(x,\vert
s\vert\nabla'\varphi)\mbox{Re}\,p_{k}(x,\xi_0,\overline \zeta)+p^{(k)}(x,\xi)\mbox{Im}\,p_{k}
(x,\xi_0,\widetilde\zeta)                                     \nonumber\\
+\vert s\vert p^{(k)}(x,\xi)\sum_{m=1}^n
p^{(m)}(x,\xi)\partial^2_{ x_k
x_m}\varphi+\vert s\vert p^{(k)}(x,\vert
s\vert\nabla'\varphi)\sum_{m=1}^n p^{(m)}(x,\vert
s\vert\nabla'\varphi)\partial^2_{
x_k
x_m}\varphi\nonumber\\
 -p^{(k)}(x,\vert
s\vert\nabla'\varphi)\mbox{Re}\,\overline p_{k}(x,\xi_0,\overline \zeta)-p^{(k)}(x,\xi)\mbox{Im}\,
p_{k}(x,\xi_0,\overline{\widetilde\zeta})\nonumber\\
 +\vert s\vert p^{(k)}(x,\xi)\sum_{m=1}^n p^{(m)}(x,\xi)
\partial^2_{
x_k x_m}\varphi+\vert s\vert p^{(k)}(x,\vert
s\vert\nabla'\varphi)\sum_{m=1}^n p^{(m)}(x,\vert
s\vert\nabla'\varphi)\partial^2_{
x_k x_m}\varphi.\nonumber
\end{eqnarray}
Therefore
\begin{eqnarray}
&&\frac 12\mbox{Im}\{\overline p(x,\xi_0,\overline{\widetilde\zeta}),
{p(x,\xi_0,\widetilde\zeta)}\} =                         \nonumber \\
& &\frac 12\mbox{Im}\left(\sum_{k=0}^n\partial_{
\xi_k}\overline p(x,\xi_0,\overline{\widetilde\zeta})\partial_{ x_k}
p(x,\xi_0,\widetilde\zeta)-\partial_{
x_k}\overline p(x,\xi_0,\overline{\widetilde\zeta})\partial_{ \xi_k}
p(x,\xi_0,\widetilde\zeta)\right)\nonumber\\
& &=-\rho\mbox{Re}\,p_{0}(x,\xi)+\rho
s^2\sum_{m=1}^n
p^{(m)}(x,\nabla \varphi)\partial^2_{x_m x_0}\varphi+s^2\rho a_{0}(x,\nabla'\varphi,\nabla'\varphi)
                                \nonumber\\
& &\sum_{k=1}^n( -p^{(k)}(x,\vert
s\vert\nabla'\varphi)( \mbox{Re}\,p_{k}(x,\xi)-p^{(k)}(x,\vert
s\vert\nabla'\varphi))+\frac 12p^{(k)}(x,\xi)\mbox{Im}\, p_{k}(x,\xi_0,\widetilde\zeta)\nonumber\nonumber\\
& &-\frac 12 p^{(k)}(x,\xi)\mbox{Im}\, p_{k}(x,\xi_0,\overline{\widetilde\zeta}))\nonumber\\
& & +\sum_{m,k=1}^n(\vert s\vert
p^{(k)}(x,\xi)p^{(m)}(x,\xi)+\vert s\vert p^{(k)}(x,\vert
s\vert\nabla'\varphi)p^{(m)}(x,\vert
s\vert\nabla'\varphi))\partial^2_{
x_k x_m}\varphi .
\end{eqnarray}
Observing that $\partial^2_{ x_k
x_m}\varphi=(\lambda^2\partial_{x_k}\psi \partial_{x_m}\psi+\lambda\partial^2_{x_kx_m}\psi)
\frac{e^{\lambda\psi}}{{\widetilde\ell}}$ for any $k,m\in \{1,\dots, n\}$,
we have
\begin{eqnarray}
&&I=\sum_{m,k=1}^n(\vert s\vert
p^{(k)}(x,\xi)p^{(m)}(x,\xi)+\vert s\vert p^{(k)}(x,\vert
s\vert\nabla'\varphi)p^{(m)}(x,\vert
s\vert\nabla'\varphi))\partial^2_{
x_k
x_m}\varphi\nonumber\\
& & =\lambda^2\vert s\vert
(a(x,\xi',\nabla'\psi)^2+s^2\frac{e^{2\lambda\psi(x)}}
{{\widetilde\ell}^{2}}a(x,\nabla'\psi,\nabla'\psi)^2)
\frac{e^{\lambda\psi}}{{\widetilde\ell}}+\nonumber\\
& &\sum_{m,k=1}^n(\vert s\vert
p^{(k)}(x,\xi)p^{(m)}(x,\xi)+\vert s\vert p^{(k)}(x,\vert
s\vert\nabla'\varphi)p^{(m)}(x,\vert
s\vert\nabla'\varphi))\lambda\psi_{x_kx_m}\frac{e^{\lambda\psi}}{{\widetilde\ell}}.\nonumber
\end{eqnarray}
Since $(x,\xi,s)\in \mathcal S$ providing that $\lambda$ is sufficiently large the following inequality holds
$$
a(x,\xi', \xi')=s^2a(x,\nabla '\varphi,\nabla '\varphi)\ge
 C_{22}\lambda ^2 \vert (\xi',s\frac{e^{\lambda\psi}}{{\widetilde\ell}})\vert^2.
$$
Taking $\hat \lambda$ sufficiently large, for all $\lambda\ge
\hat\lambda$ we have
\begin{equation}\label{mi}
I\ge \frac {\lambda^4}{2} C_{23}\vert
s\vert\frac{e^{\lambda\psi}}{{\widetilde\ell}}
\vert(\xi',s\frac{e^{\lambda\psi}}{{\widetilde\ell}})\vert ^2\quad\forall
(x,\xi,s)\in \mathcal S,
\end{equation}
where positive constant  $C_{23}$ is independent of $(\lambda,x,\xi,s).$

\noindent Finally observing that  $$ \vert\xi_0\vert\le \vert
a(x,\xi',\vert s\vert\nabla'\varphi)\vert/\rho(x) \quad\forall
(x,\xi,s)\in\mathcal S$$  from  (\ref{mi}) we find
\begin{equation}\label{mi1}
I\ge  C_{24}\frac{\lambda^4e^{\lambda\psi}}{2{\widetilde\ell}}\vert s\vert
M^2\left (\xi,s\frac{e^{\lambda\psi}}{{\widetilde\ell}
}\right )\quad\forall (x,\xi,s)\in \mathcal S,
\end{equation}
where positive constant  $C_{24}$ is independent of $(\lambda,x,\xi,s).$ On the other
hand
\begin{eqnarray}\label{mo}
&&\vert -\rho\mbox{Re}\,p_{0}(x,\xi)
+ \rho s^2a_{0}(x,\nabla'\varphi,
\nabla'\varphi)
+ \rho s^2\sum_{m=1}^n
p^{(m)}(x,\nabla'\varphi)\partial^2_{
x_k x_0}\varphi\nonumber\\
& & \sum_{m=1}^n-p^{(m)}(x,\vert
s\vert\nabla'\varphi)(p_{m}(x,\xi)-p_{m}(x,\vert
s\vert\nabla'\varphi))+\frac 12p^{(m)}(x,\xi)\mbox{Im}\,p_{m}(x,\xi_0,\widetilde\zeta)\nonumber\nonumber\\
&&-\frac 12p^{(m)}(x,\xi)\mbox{Im}\,\overline p_{m}(x,\xi_0,\overline{\widetilde\zeta})\vert\le C_{25}\vert s\vert \lambda^2\frac{e^{\lambda\psi}}{{\widetilde\ell}}
M^2\left (\xi,s\frac{e^{\lambda\psi}}{{\widetilde\ell}}\right).
\end{eqnarray}
Inequalities (\ref{mo}), (\ref{mi1}) imply (\ref{ma}). Proof of the fact   that function  $\varphi$ satisfies Condition \ref{A2} is same.
In order to prove inequalities  (\ref{victoryA}) we differentiate function $\varphi$ on $(0,\delta/2)$:
$$
\partial_{x_0}\varphi=\frac{(e^{\lambda\psi}-e^{2\lambda\Vert \psi\Vert_{C^0(\overline Q_\delta)}})}{{\widetilde\ell}^2(x_0)}+\frac{\lambda e^{\lambda\psi}\partial_{x_0}\psi}{{\widetilde\ell}(x_0)}=\frac{(e^{\lambda\psi}-e^{2\lambda\Vert \psi\Vert_{C^0(\overline Q_\delta)}})}{{\widetilde\ell}^2(x_0)}-\frac{\lambda (\mbox{\bf v}_1,\nabla'\psi)e^{\lambda\psi}}{{\widetilde\ell}(x_0)}
$$
$$\le \frac{(e^{\lambda\psi}-e^{2\lambda\Vert \psi\Vert_{C^0(\overline Q_\delta)}})}{{\widetilde\ell}^2(x_0)}-\frac{\lambda \Vert \mbox{\bf v}_1\Vert_{C^0(\overline Q_\delta)}\Vert \psi_0\Vert_{C^3(\overline \Omega)}e^{\lambda\psi}}{{\widetilde\ell}(x_0)}.$$ Hence for all $\lambda$ sufficiently large we proved the first inequality in (\ref{victoryA}) on $(0,\delta/2)$ . Taking the derivative of the function $\varphi$ on $(\delta/2, \delta]$ we have
$$
\partial_{x_0}\varphi=\frac{{\widetilde\ell}'(x_0)}{{\widetilde\ell}^2(x_0)}(e^{\lambda\psi}-e^{2\lambda\Vert \psi\Vert_{C^0(\overline Q_\delta)}})+\frac{\lambda e^{\lambda\psi}\partial_{x_0}\psi}{{\widetilde\ell}(x_0)}=\frac{-{\widetilde\ell}'(x_0)}{{\widetilde\ell}^2(x_0)}(e^{2\lambda\Vert \psi\Vert_{C^0(\overline Q_\delta)}}-e^{\lambda\psi})-\frac{\lambda (\mbox{\bf v}_1,\nabla'\psi)e^{\lambda\psi}}{{\widetilde\ell}(x_0)}
$$
$$\le \frac{\mbox{inf}_{x_0\in [\delta/2,\delta]}(-{\widetilde\ell}'(x_0))}{{\widetilde\ell}^2(x_0)}(e^{2\lambda\Vert \psi\Vert_{C^0(\overline Q_\delta)}}-e^{\lambda\psi})-\frac{\lambda\Vert \mbox{\bf v}_1\Vert_{C^0(\overline Q_\delta)}\Vert \psi_0\Vert_{C^3(\overline \Omega)}e^{\lambda\psi}}{{\widetilde\ell}(x_0)}.$$
 Hence for all $\lambda$ sufficiently large we proved the first inequality in (\ref{victoryA}) on $(\delta/2,\delta]$ .
The proof of the second inequality in (\ref{victoryA}) is the same.
Proof of proposition is complete.
$\blacksquare$

%Let $x$ be an arbitrary point of $[0,\delta]\time\Omega$ and $(t,y(t))$ be characteristic curve which pass through this point at moment $x_0$: %$(x_0, y(x_0))=x.$
%Then  for $x_0>0$  using the method of  characteristic we solve equation (\ref{ZZ1}):
%$$\rho(x)= \rho(0,y(0))e^{\int_0^{x_0}div \mbox{\bf v}_1(t,y(t)) dt}+e^{\int_0^{x_0}div \mbox{\bf v}_1(t,y(t)) dt}\int_0^{x_0}\mbox{div}(\mbox{\bf %v}\rho_2)(t,y(t)) e^{-\int_0^{t}div \mbox{\bf v}_1(s,y(s)) ds}dt.
%$$
%Multiplying the above formula by $e^{s\varphi}$ we have
%\begin{eqnarray}\label{globus}
%e^{s\varphi(x)}\rho(x)= e^{s\varphi(x)}\rho(0,y(0))e^{\int_0^{x_0}div \mbox{\bf v}_1(t,y(t)) dt}\nonumber\\+e^{s\varphi(x)}e^{\int_0^{x_0}div %\mbox{\bf v}(t,y(t)) dt}\int_0^{x_0}\mbox{div}(\mbox{\bf v}\rho_2)(t,y(t)) e^{-\int_0^{t}div \mbox{\bf v}_1(s,y(s)) ds}dt.
%\end{eqnarray}
%So, by (\ref{pinok220}), for all positive $s$ from (\ref{globus}) we obtain the inequality
%\begin{eqnarray}\label{globus1}
%\vert e^{s\varphi(x)}\rho(x)\vert\le e^{s\varphi(0,y(0))}\rho(0,y(0))\vert e^{\int_0^{x_0}div \mbox{\bf v}_1(t,y(t)) dt}\vert\\
%+e^{\int_0^{x_0}div \mbox{\bf v}_1(t,y(t)) dt}e^{s\varphi(x)}\int_0^{x_0} \vert\mbox{div}(\mbox{\bf v}\rho_2)(t,y(t))\vert e^{-\int_0^{t}div %\mbox{\bf v}_1(s,y(s)) ds} dt.\nonumber
%\end{eqnarray}
%Integrating the inequality (\ref{globus1}) over $Q_\delta$ using (\ref{pinok22}) and Proposition \ref{gonka} for all $s\ge 1$ we have

Next we prove
\begin{proposition}\label{estimate} Let function $p\in L^2(Q_\delta),\rho_0\in L^2(\Omega)$ and $\widetilde \rho\in C^0([-\delta,\delta];L^2(\Omega))$ be solution to the initial value problem
\begin{equation}\label{LAZZ1}
-L_{\mbox{\bf v}_1}(x,D)\widetilde \rho =p\quad\mbox{in}\,\,Q_\delta, \quad \widetilde\rho(0,\cdot)=\rho_0.
\end{equation}
Then where exist a constant $C_{26}$ independent of $s$ such that for all $s\ge 1$ we have
\begin{equation}\label{zamok}
\Vert e^{s\varphi}\widetilde \rho\Vert_{L^2(Q_\delta)}\le C_{26}(\Vert e^{s\varphi(0,\cdot)}\rho_0\Vert_{L^2(\Omega)}+\Vert \frac{p}{{s\widetilde \varphi}} e^{s\varphi}\Vert_{L^2(Q_\delta)}).
\end{equation}
\end{proposition}
{\bf Proof.}
Then  for $x_0>0$  using the method of  characteristic we solve equation (\ref{LAZZ1}):
\begin{equation}\label{ggo}\widetilde \rho(x)= \rho_0(y(0))e^{\int_0^{x_0}\mbox{div}\, \mbox{\bf v}_1(t,y(t)) dt}+e^{\int_0^{x_0}div \mbox{\bf v}_1(t,y(t)) dt}\int_0^{x_0}p(t,y(t)) e^{-\int_0^t\mbox{div}\, \mbox{\bf v}_1(s,y(s)) ds}dt.
\end{equation}
Consider $\mathcal F_{x_0}$ the diffeomorphism of domain $\Omega$ into $\Omega$  defined in the following way:
$$
\mathcal F_{x_0}(x')=y(0)\quad \mbox{where function $y$ solves the Cauchy problem}\,\, \frac{d y}{dt}=\mbox{\bf v}_1(t,y),  \quad y(x_0)=x'.
$$
By (\ref{vobla})
$$
\Vert e^{s\varphi(x)}\rho_0(y(0))e^{\int_0^{x_0}div \mbox{\bf v}_1(t,y(t)) dt}\Vert_{L^2(\Omega)}\le \Vert e^{s\varphi(0,y(0))}\rho_0(y(0))e^{\int_0^{x_0}div \mbox{\bf v}_1(t,y(t)) dt}\Vert_{L^2(\Omega)}.
$$
The short computations imply
$$
\Vert e^{s\varphi(0,y(0))}\rho_0(y(0))e^{\int_0^{x_0}div \mbox{\bf v}_1(t,y(t)) dt}\Vert_{L^2(\Omega)}\le C_{27}
\Vert e^{s\varphi(0,y(0))}\rho_0(y(0))\Vert_{L^2(\Omega)}=
$$
$$
\Vert e^{s\varphi(0,\mathcal F_{x_0}(x'))}\rho_0(\mathcal F_{x_0}(x'))\Vert_{L^2(\Omega)}.
$$
Making the change of variables $\mathcal F_{x_0}(x')=z=(z_1,\dots,z_n)$ we have
$$
\Vert e^{s\varphi(0,\mathcal F_{x_0}(x'))}\rho_0(\mathcal F_{x_0}(x'))\Vert_{L^2(\Omega)}=\Vert e^{s\varphi(0,z)}\rho_0(z) \vert \mbox{det}\, (\mathcal F_{x_0}^{-1})'\vert^\frac 12\Vert_{L^2(\Omega)}\le
$$
$$
\le C_{28}\Vert e^{s\varphi(0,\cdot)}\rho_0\Vert_{L^2(\Omega)}.
$$
Therefore
\begin{equation}\label{oskol}
\int_{-\delta}^\delta\Vert e^{s\varphi(x)}\rho_0(y(0))e^{\int_0^{x_0}div \mbox{\bf v}_1(t,y(t)) dt}\Vert_{L^2(\Omega)}^2dx_0\le 2\delta C_{29} \Vert e^{s\varphi(0,\cdot)}\rho_0\Vert^2_{L^2(\Omega)}.
\end{equation}
Next we estimate the $L^2$ norm second term in the right hand side of (\ref{ggo}).
Then
$$
\int_{-\delta}^\delta\int_\Omega \left (e^{\int_0^{x_0}div \mbox{\bf v}_1(t,y(t)) dt}e^{s\varphi(x)}\int_0^{x_0} \vert p(t,y(t))\vert e^{-\int_0^{t}div \mbox{\bf v}_1(s,y(s)) ds} dt\right)^2dx_0dx'
$$
$$\le C_{30} \int_\Omega \int_{-\delta}^\delta\left (e^{s\varphi(x)}\int_0^{x_0} p(t,y(t)) dt\right)^2dx_0dx'
$$
%$$
%\le C_{29} \int_\Omega \int_{-\delta}^\delta\left (e^{s\varphi(x)}\int_0^{x_0} p(t,y(t)) dt\right)^2dx_0dx'\le C \int_\Omega %\int_{-\delta}^\delta\left (e^{s\varphi(x)}\int_0^{x_0}\vert g(t,y(t))\vert dt\right)^2dx_0dx'
%$$
$$
\le C_{31} \int_\Omega \int_{-\delta}^\delta \frac { (L_{\mbox{\bf v}_1}(x,D)e^{2s\varphi(x)}) }{2sL_{\mbox{\bf v}_1}(x,D) \varphi}  \left ( \int_0^{x_0}\ p(t,y(t))dt\right)^2 dx_0dx'.
$$
Observe that the function $\frac { 1}{L_{\mbox{\bf v}_1}(x,D) \varphi}\in W_\infty^1((0,\delta)\times \Omega)\cap W_\infty^1((-\delta,0)\times \Omega).$ Therefore
%$$
% \int_\Omega \int_{-\delta}^\delta e^{2s\varphi(x)}\left\vert\frac { x_0}{L_{\mbox{\bf v}_1}(x,D) \varphi} L_{\mbox{\bf v}_1}(x,D) \varphi %\right\vert \int_0^{x_0}\vert p(t,y(t))\vert^2 dt\vert dx_0dx'
% $$
% $$\le C_{32} \int_\Omega \int_{-\delta}^\delta e^{2s\varphi(x)}\left \vert L_{\mbox{\bf v}_1}(x,D) \varphi \right \vert \int_0^{x_0}\vert % %p(t,y(t))\vert^2 dt\vert dx_0dx'=
% $$
 $$\mathcal I=\int_\Omega \int_{Q_\delta} \frac { (L_{\mbox{\bf v}_1}(x,D)e^{2s\varphi(x)}) }{2sL_{\mbox{\bf v}_1}(x,D) \varphi}  \left ( \int_0^{x_0}\ p(t,y(t))dt\right)^2 dx_0dx'=
 $$
 $$-\int_{Q_\delta} \frac{p}{sL_{\mbox{\bf v}_1}(x,D)\varphi}  \left ( \int_0^{x_0} p(t,y(t))dt\right)e^{2s\varphi} dx
 $$
 $$+\frac{1}{2s}\int_{Q_\delta}e^{2s\varphi}L_{\mbox{\bf v}_1}(x,D)^*\left( \frac{1}{L_{\mbox{\bf v}_1}(x,D)\varphi}\right) \left ( \int_0^{x_0} p(t,y(t))dt\right)^2 dx_0dx'
 $$
 $$-\frac{1}{2s}\int_{Q_\delta}e^{2s\varphi}\mbox{div}\,\mbox{\bf v}_1\left( \frac{1}{L_{\mbox{\bf v}_1}(x,D)\varphi}\right) \left ( \int_0^{x_0} p(t,y(t))dt\right)^2 dx_0dx'=
 \sum_{k=1}^3\mathcal I_{k}.
 $$
 The simple computations imply
 $$
 \left\vert  \frac{1}{L_{\mbox{\bf v}_1}(x,D)\varphi}\right \vert \le \frac{C_{32}}{\widetilde \varphi}\quad \mbox{and} \quad \left \vert L_{\mbox{\bf v}_1}(x,D)\left( \frac{1}{L_{\mbox{\bf v}_1}(x,D)\varphi}\right)\right \vert\le \frac{C_{33}}{\widetilde \varphi}\quad \mbox{on}\,\, Q_\delta .
 $$
 Therefore
 $$
 \vert\mathcal I_{2}\vert+
 \vert\mathcal I_{3}\vert=\left\vert \frac{1}{2s}\int_{Q_\delta}e^{2s\varphi}\mbox{div}\,\mbox{\bf v}_1\left( \frac{1}{L_{\mbox{\bf v}_1}(x,D)^*\varphi}\right) \left ( \int_0^{x_0} p(t,y(t))dt\right)^2 dx_0dx'\right\vert
 $$
 $$+\left \vert\frac{1}{2s}\int_{Q_\delta}e^{2s\varphi}L_{\mbox{\bf v}_1}(x,D)^*\left( \frac{1}{L_{\mbox{\bf v}_1}(x,D)\varphi}\right) \left ( \int_0^{x_0} p(t,y(t))dt\right)^2 dx_0dx'\right\vert
 $$
 \begin{equation}\label{sharmanka}
 \le C_{34} \int_{Q_\delta}e^{2s\varphi}\frac{1}{s\widetilde \varphi}\left ( \int_0^{x_0} p(t,y(t))dt\right)^2  dx_0dx'.
 \end{equation}

 By (\ref{vobla}) there exists a constant  $C_{35}$ independent of $s$ such that
 $$
 \vert\mathcal I_{1}\vert\le \frac{ \mathcal I}{2}+C_{35}\Vert pe^{s\varphi}/s\widetilde\varphi\Vert^2_{L^2(Q)}.
 $$

 This inequality, (\ref{sharmanka}) and (\ref{oskol}) imply (\ref{zamok}). Proof of the proposition is complete. $\blacksquare$

 Applying the Proposition \ref{estimate} to equation (\ref{ZZ1}) we have
\begin{equation}\label{gora}
\Vert e^{s\varphi}\rho\Vert_{L^2(Q_\delta)}\le C_{36}(\Vert e^{s\varphi(0,\cdot)}\rho(0,\cdot)\Vert_{L^2(\Omega)}+\Vert \frac{1}{{s\widetilde \varphi}}\mbox{div}\,\mbox{\bf v}e^{s\varphi}\Vert_{L^2(Q_\delta)}+\Vert \frac{1}{{s\widetilde \varphi}} \mbox{\bf v}e^{s\varphi}\Vert_{L^2(Q_\delta)}) \quad \forall s\ge s_0,
\end{equation}
where $s_0$ is sufficiently large.
From (\ref{ZZ1}) we have
\begin{equation}\label{ZIMA}
L_{\mbox{\bf v}_1}(x,D)\partial^\alpha \rho=f_\alpha\quad \mbox{in} \,Q,
\end{equation}  $f_\alpha= [L_{\mbox{\bf v}_1},\partial^\alpha ]\rho-\partial^\alpha\mbox{div}\,(\mbox{\bf v}\rho_2)-\partial^\alpha(\mbox{div}\,\mbox{\bf v}_1\rho).$
%Repeating arguments  (\ref{globus}), (\ref{globus1}) we obtain

%$$
%\vert e^{s\varphi(x)}\partial^\alpha\rho(x)\vert\le e^{s\varphi(0,y(0))}\vert\partial^\alpha\rho(0,y(0))\vert e^{\int_0^{x_0}div \mbox{\bf %v}_1(t,y(t)) dt}
%$$
%$$+e^{\int_0^{x_0}div \mbox{\bf v}_1(t,y(t)) dt}\int_0^{x_0} e^{s\varphi(t,y(t))}\vert f_\alpha (t,y(t))\vert e^{-\int_0^{t}div \mbox{\bf %v}_1(s,y(s)) ds} dt.
%$$
%So, from the  above equality, applying (\ref{pinok22}) and Proposition \ref{gonka} we have
%\begin{equation}\label{gora1}
%\Vert e^{s\varphi}\partial^\alpha\rho\Vert_{L^2(Q_\delta)}\le C_{12}(\Vert %e^{s\varphi(0,\cdot)}\partial^\alpha\rho(0,\cdot)\Vert_{L^2(\Omega)}+\frac{1}{\root\of s}\Vert f_\alpha e^{s\varphi}\Vert_{L^2(Q_\delta)}).
%\end{equation}
Let $\vert\alpha\vert=1,\alpha_0=0.$  Applying to equation (\ref{ZIMA}) Proposition \ref{estimate}  and using (\ref{gora}) we obtain for all $s\ge s_1$
%\begin{eqnarray}\label{gora2}
%\Vert f_\alpha e^{s\varphi}\Vert_{L^2(Q_\delta)}\le C_{13}(\Vert \rho e^{s\varphi}\Vert_{L^2(Q_\delta)}+\sum_{\vert\beta\vert\le 1,\beta_0=0} %\Vert  e^{s\varphi} \partial^\beta \mbox{\bf v}\Vert_{L^2(Q_\delta)}+\Vert e^{s\varphi} \nabla'\mbox{div}\, \mbox{\bf %v}\Vert_{L^2(Q_\delta)})\nonumber\\
%\le C_{14}(\Vert e^{s\varphi(0,\cdot)}\rho(0,\cdot)\Vert_{L^2(\Omega)}+\Vert  e^{s\varphi}\nabla'\mbox{div}\, \mbox{\bf %v}\Vert_{L^2(Q_\delta)}+\sum_{\vert\beta\vert\le 1,\beta_0=0} \Vert  e^{s\varphi}\partial^\beta \mbox{\bf v}\Vert_{L^2(Q_\delta)})
%\end{eqnarray} for all $s\ge 1.$
%From (\ref{gora1}) and (\ref{gora2}) for all $s\ge 1$ we have
\begin{eqnarray}\label{gora3}
\Vert e^{s\varphi}\nabla '\rho\Vert_{L^2(Q_\delta)}\le C_{37}(\sum_{\vert\alpha\vert\le 1,\alpha_0=0}\Vert e^{s\varphi(0,\cdot)}\partial^\alpha\rho(0,\cdot)\Vert_{L^2(\Omega)}\nonumber\\+\Vert e^{s\varphi}  \frac{1}{{s\widetilde \varphi}}\nabla'\mbox{div}\, \mbox{\bf v}\Vert_{L^2(Q_\delta)}+\sum_{\vert\alpha\vert\le 1,\alpha_0=0} \Vert  \frac{1}{{s\widetilde \varphi}} e^{s\varphi}\partial^\alpha \mbox{\bf v}\Vert_{L^2(Q_\delta)}).
\end{eqnarray}
Let $\vert\alpha\vert=2,\alpha_0=0.$ Applying to equation (\ref{ZIMA}) Proposition \ref{estimate}  and using (\ref{gora}), (\ref{gora3}) we obtain
%\begin{eqnarray}\label{gora4}
%\Vert f_\alpha e^{s\varphi}\Vert_{L^2(Q_\delta)}\le C_{16}(\Vert (\vert \nabla '\rho \vert +\rho)e^{s\varphi}\Vert_{L^2(Q_\delta)}+s^{-\frac %12}\sum_{\vert\beta\vert\le 2,\beta_0=0} \Vert \partial^\beta \mbox{\bf v}\Vert_{L^2(Q_\delta)}\nonumber\\+s^{-\frac 12}\sum_{\vert %\alpha\vert=2,\alpha_0=0}\Vert \partial^\alpha\mbox{div}\, \mbox{\bf v}\Vert_{L^2(Q_\delta)})
%\le C_{17}(\Vert e^{s\varphi(0,\cdot)}\sum_{\vert\beta\vert\le %2}\partial^\beta\rho(0,\cdot)\Vert_{L^2(\Omega)}\nonumber\\+s^{-\frac{1}{2}}\sum_{\vert\beta\vert\le 2,\beta_0=0} \Vert  %e^{s\varphi}\partial^\beta \mbox{\bf v}\Vert_{L^2(Q_\delta)}+s^{-\frac 12}\sum_{\vert \beta\vert=2,\beta_0=0}\Vert  %e^{s\varphi}\partial^\beta\mbox{div}\, \mbox{\bf v}\Vert_{L^2(Q_\delta)}).
%\end{eqnarray}
%Inequalities (\ref{gora1}) and (\ref{gora4}) imply
\begin{eqnarray}\label{gora5}
\sum_{\vert\alpha\vert\le 2,\alpha_0=0}\Vert e^{s\varphi}\partial^\alpha\rho\Vert_{L^2(Q_\delta)}\le C_{38}(\sum_{\vert\alpha\vert\le 2,\alpha_0=0}\Vert e^{s\varphi(0,\cdot)}\partial^\alpha\rho(0,\cdot)\Vert_{L^2(\Omega)}+\nonumber\\+\sum_{\vert \alpha\vert\le 2,\alpha_0=0}\Vert e^{s\varphi}\frac{1}{{s\widetilde \varphi}}\partial^\alpha\mbox{div}\, \mbox{\bf v}\Vert_{L^2(Q_\delta)}+\sum_{\vert\alpha\vert\le 2,\alpha_0=0} \Vert e^{s\varphi} \frac{1}{{s\widetilde \varphi}}\partial^\alpha \mbox{\bf v}\Vert_{L^2(Q_\delta)})\quad \forall s\ge s_2.
\end{eqnarray}
From (\ref{ZZ1}) and (\ref{gora5})  for all $s\ge s_3$ we have
\begin{eqnarray}\label{okorok}
\Vert \nabla '\partial_{x_0}\rho e^{s\varphi}\Vert_{L^2(Q_\delta)}+\Vert \partial_{x_0}\rho e^{s\varphi}\Vert_{L^2(Q_\delta)}\le C_{39} (\sum_{\vert\alpha\vert\le 2,\alpha_0=0}\Vert e^{s\varphi}\partial^\alpha\rho\Vert_{L^2(Q_\delta)}\\+\sum_{\vert \alpha\vert\le 2,\alpha_0=0}\Vert e^{s\varphi}\frac{1}{{s\widetilde \varphi}}\partial^\alpha\mbox{div}\, \mbox{\bf v}\Vert_{L^2(Q_\delta)}+\sum_{\vert\alpha\vert\le 2,\alpha_0=0}\Vert e^{s\varphi}\frac{1}{{s\widetilde \varphi}}\partial^\alpha\mbox{\bf v}\Vert_{L^2(Q_\delta)})\nonumber\\
 \le C_{40}(\Vert e^{s\varphi(0,\cdot)}\rho(0,\cdot)\Vert_{H^{2,s}(\Omega)}+\sum_{\vert \alpha\vert\le 2,\alpha_0=0}\Vert e^{s\varphi}\frac{1}{{s\widetilde \varphi}}\partial^\alpha\mbox{div}\, \mbox{\bf v}\Vert_{L^2(Q_\delta)}+\sum_{\vert\alpha\vert\le 2,\alpha_0=0} \Vert e^{s\varphi} \frac{1}{{s\widetilde \varphi}}\partial^\alpha \mbox{\bf v}\Vert_{L^2(Q_\delta)}).\nonumber
\end{eqnarray}
We differentiate (\ref{ZZ1}) respect to variable $x_0$
\begin{equation}\label{ZZZ1}
L_{\mbox{\bf v}_1}(x,D)^*\partial^\alpha\partial_{x_0} \rho=\widetilde f_\alpha\quad \mbox{in} \,Q,
\end{equation}  where $\widetilde f_\alpha= [L_{\mbox{\bf v}_1}^*,\partial^\alpha\partial_{x_0} ]\rho-\partial^\alpha\partial_{x_0}\mbox{div}(\mbox{\bf v}\rho_2)-\partial^\alpha\partial_{x_0}(\mbox{div}\,\mbox{\bf v}_1\rho).$
%Repeating arguments  (\ref{globus}), (\ref{globus1}) we obtain

From (\ref{ZZZ1}) and (\ref{okorok}) for ass $s\ge s_4$
\begin{eqnarray}\label{okorokZ}
\sum_{\vert\alpha\vert\le 2,\alpha_0=0}\Vert \partial^\alpha\partial_{x_0}\rho e^{s\varphi}\Vert_{L^2(Q_\delta)}
 \le C_{41}(\sum_{\vert\alpha\vert\le 2,\alpha_0=0}\Vert e^{s\varphi(0,\cdot)}\partial^\alpha\partial_{x_0}\rho(0,\cdot)\Vert_{L^2(\Omega)}\nonumber\\+\sum_{\vert \alpha\vert\le 2,\alpha_0=0}
 \Vert e^{s\varphi}\frac{1}{{s\widetilde \varphi}}\partial^\alpha\widetilde f_\alpha\Vert_{L^2(Q_\delta)})\nonumber\\
  \le C_{42}(\Vert e^{s\varphi(0,\cdot)}\rho(0,\cdot)\Vert_{H^{3,s}(\Omega)}+\Vert e^{s\varphi(0,\cdot)}\mbox{\bf v}(0,\cdot)\Vert_{H^{3,s}(\Omega)}\\+\sum_{\vert \alpha\vert\le 2,\alpha_0=0}\Vert e^{s\varphi}\frac{1}{{s\widetilde \varphi}}\partial^\alpha\mbox{div}\, \mbox{\bf v}\Vert_{L^2(Q_\delta)}+\sum_{\vert\alpha\vert\le 2,\alpha_0=0} \Vert e^{s\varphi} \frac{1}{{s\widetilde \varphi}}\partial^\alpha \mbox{\bf v}\Vert_{L^2(Q_\delta)}).\nonumber
\end{eqnarray}

Observe that the function $\varphi$ defined by (\ref{gold}) with function $\psi$ which verifies    (\ref{oop}) - (\ref{coppp}) satisfies  (\ref{giga1})
- (\ref{lomka}), Condition \ref{A1} and Condition \ref{A2} provided that parameter $\lambda$ is large enough.

Applying to equation (\ref{ZZ2}) Carleman estimate (\ref{2.9'}) for all $s\ge s_5$ we have
\begin{eqnarray}\label{K2.9'!}
\Vert \mbox{\bf v}\Vert_{\mathcal B(\varphi,s,Q_\delta)}
+ \Vert\partial^2_{\vec \nu} \text{\bf v}
e^{s\varphi}\Vert_{H^{\frac 14,\frac 12,\widetilde s}( \Sigma_{0,\delta})}
+\Vert \partial_{\vec \nu}
\text{\bf v}
e^{s\varphi}\Vert_{H^{\frac 34,\frac 32,\widetilde s}( \Sigma_{0,\delta})}\le C_{43}(\Vert \text{\bf F}
e^{s\varphi}\Vert_{\mathcal Y(\varphi,s,Q_\delta)}
\\
+\Vert (d \omega_{\mbox {\bf v}},\mbox{div}\,\mbox {\bf v})
e^{s\varphi}\Vert_{H^{\frac 34,\frac 32,\widetilde s}(\widetilde\Sigma_\delta)}+ \Vert \partial_{\vec \nu}(d \omega_{\mbox {\bf v}},\mbox{div}\,\mbox {\bf v})
e^{s\varphi}\Vert_{H^{\frac 14,\frac 12,\widetilde s}(\widetilde\Sigma_\delta)}+\Vert\widetilde\varphi \partial_{x_0}\partial_{\vec \nu} \text{\bf v}
e^{s\varphi}\Vert_{ L^2(\widetilde\Sigma_\delta)}),
\nonumber
\end{eqnarray}
where
$\mbox{\bf F}=-\rho\partial_{x_0}\mbox{\bf v}_2-(\mbox{\bf v},\nabla)\mbox{\bf v}_2-h(\rho_1)\nabla' \rho-(h(\rho_1)-h(\rho_2))\nabla \rho_2+R\mbox{\bf f},\widetilde\Sigma_\delta=[-\delta,\delta]\times\widetilde\Gamma, \widetilde\Sigma_{0,\delta}=[-\delta,\delta]\times\partial\Omega\setminus\widetilde\Gamma.$

From (\ref{gora5}) and (\ref{K2.9'!})  for all $s\ge s_6$ we obtain
\begin{eqnarray}\label{KK2.9'}
\Vert \mbox{\bf v}\Vert_{\mathcal B(\varphi,s,Q_\delta)}
+ \Vert\partial^2_{\vec \nu} \text{\bf v}
e^{s\varphi}\Vert_{H^{\frac 14,\frac 12,\widetilde s}( \Sigma_{0,\delta})}
+\Vert \partial_{\vec \nu}
\text{\bf v}
e^{s\varphi}\Vert_{H^{\frac 34,\frac 32,\widetilde s}( \Sigma_{0,\delta})}\le C_{44}(\Vert R \text{\bf f}
e^{s\varphi}\Vert_{\mathcal Y(\varphi,s,Q_\delta)}
\nonumber\\
+\Vert (d \omega_{\mbox {\bf v}},\mbox{div}\,\mbox {\bf v})
e^{s\varphi}\Vert_{H^{\frac 34,\frac 32,\widetilde s}(\widetilde\Sigma_\delta)}+ \Vert \partial_{\vec \nu}(d \omega_{\mbox {\bf v}},\mbox{div}\,\mbox {\bf v})
e^{s\varphi}\Vert_{H^{\frac 14,\frac 12,\widetilde s}(\widetilde\Sigma_\delta)}+\Vert\widetilde\varphi \partial_{x_0}\partial_{\vec \nu} \text{\bf v}
e^{s\varphi}\Vert_{ L^2(\widetilde\Sigma_\delta)}\nonumber\\
+\Vert e^{s\varphi(0,\cdot)}\rho(0,\cdot)\Vert_{H^{2,s}(\Omega)}).
\end{eqnarray}
Next we differentiate equations (\ref{ZZ1})-(\ref{ZZ3}) respect to variable $x_0.$ Setting $\widetilde \rho=\partial_{x_0}\rho$ and $\mbox{\bf u}=\partial_{x_0}\mbox{\bf v}$ we obtain

%\begin{equation}\label{AZZ1}
%\partial_{x_0}\widetilde \rho+\mbox{div}(\mbox{\bf v}_1\widetilde \rho)=-\mbox{div}(\mbox{\bf u}\rho_2)-\mbox{div}(\mbox{\bf v}\partial_{x_0}\rho_2)-
%\mbox{div}(\partial_{x_0}\mbox{\bf v}_1 \rho)\quad \mbox{in} \,Q,
%\end{equation}
\begin{eqnarray}\label{AZZ2}
\rho_1\partial_{x_0}\mbox{\bf u}+\widetilde \rho\partial_{x_0}\mbox{\bf v}_2+L_{\lambda,\mu}(x',D')\mbox{\bf u} +(\mbox{\bf v}_1,\nabla')\mbox{\bf u}+(\mbox{\bf u},\nabla')\mbox{\bf v}_2+h(\rho_1)\nabla' \widetilde\rho\nonumber\\+\partial_{x_0}(h(\rho_1)-h(\rho_2))\nabla ' \rho_2=\partial_{x_0}R(x)\mbox{\bf f}+\mbox{\bf G}\quad \mbox{in} \,Q, \quad \mbox{\bf u}\vert_{\Sigma}=0,
\end{eqnarray}
%$-\rho_1\partial_{x_0}\mbox{\bf v}-(\partial_{x_0}\mbox{\bf v}_1,\nabla)\mbox{\bf v}$
%\begin{equation}\label{AZZ3}
%\mbox{\bf u}\vert_{\Sigma}=0,
%\end{equation}
where $\mbox{\bf G}=-\partial_{x_0}\rho_1\partial_{x_0}\mbox{\bf v}-\rho\partial^2_{x_0}\mbox{\bf v}_2-(\partial_{x_0}\mbox{\bf v}_1,\nabla')\mbox{\bf v}-(\mbox{\bf v},\nabla')\partial_{x_0}\mbox{\bf v}_2-\partial_{x_0}h(\rho_1)\nabla' \rho-(h(\rho_1)-h(\rho_2))\nabla ' \partial_{x_0}\rho_2.$

Applying to equation (\ref{AZZ2}) Carleman estimate (\ref{2.9'}) for all $s\ge s_7$ we have
\begin{eqnarray}\label{LK2.9'}
\Vert \mbox{\bf u}\Vert_{\mathcal B(\varphi,s,Q_\delta)}
+ \Vert\partial^2_{\vec \nu} \text{\bf u}
e^{s\varphi}\Vert_{H^{\frac 14,\frac 12,\widetilde s}( \Sigma_{0,\delta})}
+\Vert \partial_{\vec \nu}
\text{\bf u}
e^{s\varphi}\Vert_{H^{\frac 34,\frac 32,\widetilde s}( \Sigma_{0,\delta})}\\\le C_{45}(\Vert \text{\bf G}
e^{s\varphi}\Vert_{\mathcal Y(\varphi,s,Q_\delta)}
+
\Vert (d \omega_{\mbox {\bf u}},\mbox{div}\,\mbox {\bf u})
e^{s\varphi}\Vert_{H^{\frac 34,\frac 32,\widetilde s}(\widetilde\Sigma_\delta)}\nonumber\\+ \Vert \partial_{\vec \nu}(d \omega_{\mbox {\bf u}},\mbox{div}\,\mbox {\bf u})
e^{s\varphi}\Vert_{H^{\frac 14,\frac 12,\widetilde s}(\widetilde\Sigma_\delta)}+\Vert\widetilde\varphi \partial_{x_0}\partial_{\vec \nu} \text{\bf u}
e^{s\varphi}\Vert_{ L^2(\widetilde\Sigma_\delta)}).
\nonumber
\end{eqnarray}
 From (\ref{LK2.9'}) using (\ref{okorokZ}) and   (\ref{KK2.9'}) we obtain for all $s\ge s_8$
\begin{eqnarray}\label{KKK2.9'}
\sum_{k=0}^1(\Vert\partial^k_{x_0} \mbox{\bf v}\Vert_{\mathcal B(\varphi,s,Q_\delta)}
+ \Vert\partial^2_{\vec \nu}\partial^k_{x_0} \text{\bf v}
e^{s\varphi}\Vert_{H^{\frac 14,\frac 12,\widetilde s}( \Sigma_{0,\delta})}
+\Vert \partial_{\vec \nu}\partial^k_{x_0}
\text{\bf v}
e^{s\varphi}\Vert_{H^{\frac 34,\frac 32,\widetilde s}( \Sigma_{0,\delta})})\le C_{46}(\Vert \text{\bf f}
e^{s\varphi}\Vert_{\mathcal Y(\varphi,s,Q_\delta)}
+\nonumber\\
\sum_{k=0}^1(\Vert \partial^k_{x_0}(d \omega_{\mbox {\bf v}},\mbox{div}\,\mbox {\bf v})
e^{s\varphi}\Vert_{H^{\frac 34,\frac 32,\widetilde s}(\widetilde\Sigma_\delta)}+ \Vert\partial^k_{x_0} \partial_{\vec \nu}(d \omega_{\mbox {\bf v}},\mbox{div}\,\mbox {\bf v})
e^{s\varphi}\Vert_{H^{\frac 14,\frac 12,\widetilde s}(\widetilde\Sigma_\delta)}+\Vert\widetilde\varphi \partial^{k+1}_{x_0}\partial_{\vec \nu} \text{\bf v}
e^{s\varphi}\Vert_{ L^2(\widetilde\Sigma_\delta)}\nonumber\\
+\Vert e^{s\varphi(0,\cdot)}\partial^k_{x_0}\rho(0,\cdot)\Vert_{H^{2,s}(\Omega)}).
\end{eqnarray}

By (\ref{victoryA}) for any $\epsilon >0$ there exist $s_0(\epsilon)$ such that for all $s\ge s_9(\epsilon)$ we have
\begin{equation}\label{Lnokia4}
\Vert \text{\bf f}
e^{s\varphi}\Vert_{\mathcal Y(\varphi,s,Q_\delta)}\le C_{47}\Vert \text{\bf f}
e^{s\varphi}\Vert_{L^2(-\delta,\delta;H^{1,\widetilde s}(\Omega))} \le \epsilon\Vert \f e^{\varphi(0,\cdot)}\Vert_{H^{1,s}(\Omega)}.
\end{equation}
On the other hand there exists a constant $C_{48}$ independent of $s$ such that
\begin{eqnarray}\label{compot3}
\Vert e^{s\varphi(0,\cdot)}\partial_{x_0}\mbox{\bf v}(0,\cdot)\Vert^2_{H^{1,s}(\Omega)}\le C_{48}\sum_{k=0}^1\int_{Q_\delta}\left(\frac{1}{s^2\widetilde\varphi^2}\vert\partial^{k+1}_{x_0}\nabla' \mbox{\bf v}\vert^2+s^2\widetilde\varphi^2\vert\partial^k_{x_0}\nabla' \mbox{\bf v}\vert^2 \right)e^{2s\varphi}dx\nonumber\\ \le C_{49} \sum_{k=0}^1\Vert\partial^k_{x_0} \mbox{\bf v}\Vert^2_{\mathcal B(\varphi,s,Q_\delta)}\quad \forall s\ge 1.
\end{eqnarray}
By (\ref{ZZ2}) for all $s\ge 1$ the following estimate is true

\begin{eqnarray}\label{compot1}
\Vert e^{s\varphi(0,\cdot)}\rho_1(0,\cdot)\partial_{x_0}\mbox{\bf v}(0,\cdot)\Vert_{H^{1, s}(\Omega)}\ge \Vert e^{s\varphi(0,\cdot)}R(0,\cdot) \mbox{\bf f}\Vert_{H^{1, s}(\Omega)}\\-C_{50}(\Vert e^{s\varphi(0,\cdot)}(\mbox{\bf v}_1(0,\cdot)-\mbox{\bf v}_2(0,\cdot))\Vert_{H^{3, s}(\Omega)}+\Vert e^{s\varphi(0,\cdot)}( \rho_1(0,\cdot)-\rho_2(0,\cdot))\Vert_{H^{2, s}(\Omega)}).\nonumber
\end{eqnarray}

By (\ref{compot}) we obtain from (\ref{compot1}) that there exist a positive constant $C_{51}$ such that
\begin{eqnarray}\label{compot2}
\Vert e^{s\varphi(0,\cdot)}\rho_1(0,\cdot)\partial_{x_0}\mbox{\bf v}(0,\cdot)\Vert_{H^{1,s}(\Omega)}\ge C_{51}\Vert e^{s\varphi(0,\cdot)} \mbox{\bf f}\Vert_{H^{1,s}(\Omega)}\\-C_{52}(\Vert e^{s\varphi(0,\cdot)}(\mbox{\bf v}_1(0,\cdot)-\mbox{\bf v}_2(0,\cdot))\Vert_{H^{3, s}(\Omega)}\nonumber+\Vert e^{s\varphi(0,\cdot)}( \rho_1(0,\cdot)-\rho_2(0,\cdot))\Vert_{H^{2,s}(\Omega)})\quad\forall s\ge 1.
\end{eqnarray}
Then (\ref{compot2}), (\ref{compot3}) imply the estimate
\begin{eqnarray}\label{compot22}
\Vert e^{s\varphi(0,\cdot)} \mbox{\bf f}\Vert_{H^{1,s}(\Omega)}\le C_{53}\left (\sum_{k=0}^1\Vert\partial^k_{x_0} \mbox{\bf v}\Vert_{\mathcal B(\varphi,s,Q_\delta)}+\Vert e^{s\varphi(0,\cdot)}(\mbox{\bf v}_1(0,\cdot)-\mbox{\bf v}_2(0,\cdot))\Vert_{H^{3, s}(\Omega)}\nonumber\right.\\
+\Vert e^{s\varphi(0,\cdot)}( \rho_1(0,\cdot)-\rho_2(0,\cdot))\Vert_{H^{2,s}(\Omega)}\bigg)\quad \forall s\ge 1.\nonumber
\end{eqnarray}
This estimate and (\ref{KKK2.9'}) imply
\begin{eqnarray}\label{compot4}
\Vert e^{s\varphi(0,\cdot)} \mbox{\bf f}\Vert_{H^{1,s}(\Omega)}\le C_{54}\left (\Vert e^{s\varphi(0,\cdot)}(\mbox{\bf v}_1(0,\cdot)-\mbox{\bf v}_2(0,\cdot))\Vert_{H^{3, s}(\Omega)}\nonumber\right.\\+\Vert e^{s\varphi(0,\cdot)}( \rho_1(0,\cdot)-\rho_2(0,\cdot))\Vert_{H^{3,\widetilde s}(\Omega)}
+ \Vert \text{\bf f}
e^{s\varphi}\Vert_{\mathcal Y(\varphi,s,Q_\delta)}
+\nonumber\\
\sum_{k=0}^1(\Vert \partial^k_{x_0}(d \omega_{\mbox {\bf v}},\mbox{div}\,\mbox {\bf v})
e^{s\varphi}\Vert_{H^{\frac 34,\frac 32,\widetilde s}(\widetilde\Sigma_\delta)}+ \Vert\partial^k_{x_0} \partial_{\vec \nu}(d \omega_{\mbox {\bf v}},\mbox{div}\,\mbox {\bf v})
e^{s\varphi}\Vert_{H^{\frac 14,\frac 12,\widetilde s}(\widetilde\Sigma_\delta)}\nonumber\\
+\Vert\widetilde\varphi \partial^{k+1}_{x_0}\partial_{\vec \nu} \text{\bf v}
e^{s\varphi}\Vert_{ L^2(\widetilde\Sigma_\delta)}+\Vert e^{s\varphi(0,\cdot)}\partial^k_{x_0}\rho(0,\cdot)\Vert_{H^{2,s}(\Omega)}))\quad\forall s\ge s_{10}.
\end{eqnarray}
Using (\ref{Lnokia4})  to estimate the norm of function $\text{\bf f}
e^{s\varphi(0,\cdot)}$  in the right hand side of (\ref{compot4}) we have
\begin{eqnarray}\label{compot5}
\Vert e^{s\varphi(0,\cdot)} \mbox{\bf f}\Vert_{H^{1,s}(\Omega)}\le C_{55}\left (\Vert e^{s\varphi(0,\cdot)}(\mbox{\bf v}_1(0,\cdot)-\mbox{\bf v}_2(0,\cdot))\Vert_{H^{3,s}(\Omega)}\right.\nonumber\\+\Vert e^{s\varphi(0,\cdot)}( \rho_1(0,\cdot)-\rho_2(0,\cdot))\Vert_{H^{3,s}(\Omega)}
+ \epsilon \Vert \f e^{\varphi(0,\cdot)}\Vert_{H^{1,s}(\Omega)}
\nonumber\\
+\sum_{k=0}^1(\Vert \partial^k_{x_0}(d \omega_{\mbox {\bf v}},\mbox{div}\,\mbox {\bf v})
e^{s\varphi}\Vert_{H^{\frac 34,\frac 32,\widetilde s}(\widetilde\Sigma_\delta)}+ \Vert\partial^k_{x_0} \partial_{\vec \nu}(d \omega_{\mbox {\bf v}},\mbox{div}\,\mbox {\bf v})
e^{s\varphi}\Vert_{H^{\frac 14,\frac 12,\widetilde s}(\widetilde\Sigma_\delta)}\nonumber\\
 +\Vert\widetilde\varphi \partial^{k+1}_{x_0}\partial_{\vec \nu} \text{\bf v}
e^{s\varphi}\Vert_{ L^2(\widetilde\Sigma_\delta)}
+\Vert e^{s\varphi(0,\cdot)}\partial^k_{x_0}\rho(0,\cdot)\Vert_{H^{2,s}(\Omega)}))\quad\forall s\ge s_{11}.
\end{eqnarray}
Then taking in (\ref{compot5}) parameter $\epsilon$ sufficiently small we obtain (\ref{(1.19)}).  Proof of theorem is complete. $\blacksquare$
%\centerline{\bf Appendix .}

\section{Appendix}\label{Q!8}

In the appendix we formulate several lemmata which were used for the proof of
Theorem \ref{opa3}. Proof of these lemmata is the same as in  e.g. \cite{IPY} or \cite{IY2018}.

The following lemma allows us to extend the definition of the
operator $A$ on Sobolev spaces.
\begin{lemma}\label{Fops0}
Let $a(\widetilde y,\widetilde\xi,s)\in L^\infty_{cl}S^{1,s}({\mathcal O}).$ Then
$A\in {\mathcal L}(H^{\frac 12,1,s}_0({\mathcal O});L^2({\mathcal O}))$ and \newline $\Vert
A\Vert_{{\mathcal L}(H^{\frac 12, 1,s}_0({\mathcal O});L^2({\mathcal
O}))}\le C_1(\pi_{L^\infty}(a)). $
\end{lemma}

\begin{lemma}\label{Fops1}
Let $a(\widetilde y,\widetilde\xi,s)\in W^{\ell,\infty}_{cl}S^{\ell,s}({\mathcal
O}).$ Then $A(\widetilde y,\widetilde D,s)^*
=A^*(\widetilde y,
\widetilde D,s)+R$, where $A^*$ is the
pseudodifferential operator with symbol
$\overline{a(\widetilde y,\widetilde\xi,s)}$ and $R\in {\mathcal
L}(H^{\frac{\ell-1}{2},\ell-1, s}_0({\mathcal O}),L^2({\mathcal O}))$ satisfies
$$
\Vert R\Vert_{{\mathcal L}(H_0^{\frac{\ell-1}{2},\ell-1,s}({\mathcal O}),L^2({\mathcal O}))}\le
C_2\pi_{W^{\ell,\infty}(\mathcal O)}(a).
$$
\end{lemma}

\begin{lemma}\label{Fops3}
Let $a(\widetilde y,\widetilde\xi,s)\in W^{1,\infty}_{cl}S^{1,s}({\mathcal O})$ and
$b(\widetilde y,\widetilde\xi,s)\in W^{1,\infty}_{cl}S^{\mu,s}({\mathcal
O}).$ Then \newline $A(\widetilde y,\widetilde D,s)B(\widetilde y,\widetilde D,s)
=C(\widetilde y,\widetilde D,s)+R_0$ where $C(\widetilde y,\widetilde D,s)$
is the operator with symbol
$a(\widetilde y,\widetilde\xi,s)b(\widetilde y,\widetilde\xi,s)$ and $R_0\in {\mathcal
L}(H_0^{\frac 12,1,s}({\mathcal O}),L^2({\mathcal O})).$ Moreover we have
$$
\Vert R_{0}\Vert_{{\mathcal L}(H_0^{\frac {1}{2},1,s}({\mathcal
O}),L^2({\mathcal O}))}\le C_3(\pi_{W^{0,\infty}(\mathcal O)}(a)\pi_{W^{1,\infty}(\mathcal O)}(b)+\pi_{W^{1,\infty}(\mathcal O)}(a)\pi_{W^{0,\infty}(\mathcal O)}(b)).
$$
%$$
%\Vert R_{0}\Vert_{{\mathcal
%L}(H_0^{\frac{\mu+\widetilde s}{2},\mu+{\widetilde s},s}({\mathcal
%O}),H^{\frac{\widetilde s+1}{2},{\widetilde s}+1,s}({\mathcal O}))}\le
%C\pi_{W^{1,\infty}(\mathcal O)}(a)\pi_{W^{1,\infty}(\mathcal O)}(b)\quad \mbox{for}\,\,\widehat j
%=0.
%$$
\end{lemma}
\noindent

\noindent
The direct consequence of Lemma \ref{Fops3} is the following
commutator estimate.

\begin{lemma}\label{Fops2}
Let $a(\widetilde y,\widetilde\xi,s)\in W^{1,\infty}_{cl}S^{1,s}({\mathcal
O})$ and $b(\widetilde y,\widetilde \xi,s)\in W^{1,\infty}_{cl}S^{1,s}({\mathcal
O}).$ \\
Then the commutator $[A,B]$ belongs to the space $ {\mathcal L}(H^{\frac 12,1,s}(
{\mathcal O});L^2({\mathcal
O}))$ and
$$\Vert [A,B]\Vert_{{\mathcal L}(H^{\frac 12,1,s}({\mathcal O});L^2({\mathcal O}))}\le
C_4(\pi_{C^0(\mathcal O)}(a)\pi_{C^0(\mathcal O)}(b)+\pi_{C^0(\mathcal O)}(a)
\pi_{W^{1,\infty}(\mathcal O)}(b)+\pi_{W^{1,\infty}(\mathcal O)}(a)\pi_{C^0(\mathcal O)}(b)).
$$
\end{lemma}
\noindent

\begin{lemma}\label{Fops4}
Let $a(\widetilde y,\widetilde \xi,s)\in W^{1,\infty}_{cl}S^{1,s}({\mathcal O})$ be a
symbol with compact support in ${\mathcal O}.$ Let  $\mathcal O_i
\subset\subset\mathcal O$ and $\overline{\mathcal O_1}\cap\overline{\mathcal O_2}=\emptyset .$
Suppose that $u\in H^{\frac 12,1,s}({\mathcal O})$ and $\mbox{supp}\,u
\subset {\mathcal O_1}.$  Then there exists a constant
$C_5$ such that
\begin{equation}\label{gabon}
\Vert A(\widetilde y,\widetilde D,s)u\Vert_{H^{\frac 12,1,s}({\mathcal O_2})}
\le \frac{C_5 \pi_{W^{1,\infty}}(a)}{\mbox{dist} (\mathcal O_2,\mathcal O_1)^{4n+3}}
\Vert u\Vert_{H^{\frac 12,1,s}({\mathcal O})}.
\end{equation}
\end{lemma}

\noindent We shall use the following variant of the G\aa rding
inequality:
\begin{lemma}\label{Fops5}
Let $p(\widetilde y,\widetilde\xi,s)\in W^{1,\infty}_{cl}S^{2,s}({\mathcal O})$ be a
symbol with compact support in ${\mathcal O}.$
Let $u\in H^{\frac 12,1,s}({\mathcal O})$ and $\mbox{supp}\,u\subset\mathcal O_1.$ Let
$\mathcal O_1\subset\subset\mathcal O_2\subset\subset\mathcal O_3
\subset\subset \mathcal O$  and $\widetilde\gamma\in C^\infty_0(\mathcal O_3)$
be a function
such that $\widetilde\gamma\vert_{\mathcal O_2}=1$  be such that
$\mbox{Re}\,p(\widetilde y,\widetilde\xi,s)> \hat C M^2(\widetilde \xi,s)$ for any
$\widetilde y\in\mathcal O_3$.  Then
\begin{eqnarray}\label{goblin}
\mbox{Re}(P(\widetilde y,\widetilde D,s)u,u)_{L^2(\mathcal O)}\ge \frac{\hat C}{2}
\Vert u\Vert^2_{H^{\frac 12,1,s}({\mathcal O})}\\-C_6[(\pi_{W^{0,\infty}(\mathcal O_3)}(p)+1)\pi_{C^{0}(\mathcal O_3)}
(\widetilde \gamma)+\sum_{k=0}^1(\pi_{W^{k,\infty}(\mathcal O_3)}(p)+1)\pi_{C^{1-k}(\mathcal O_3)}(\widetilde\gamma)+\pi_{C^{1}(\mathcal O_3)}(\widetilde\gamma)]^2\Vert u\Vert^2_{L^2({\mathcal O})}.\nonumber
\end{eqnarray}
\end{lemma}
{\bf Proof.}
Consider the
pseudodifferential operator $A(\widetilde y,\widetilde D,s)$ with symbol
$A(\widetilde y,\widetilde \xi,s)=$\newline$(\widetilde\gamma\mbox{Re}\, p(\widetilde y,\widetilde \xi,s)
-\widetilde\gamma\frac{\hat C}{2} M^2(\widetilde \xi,s))^\frac 12\in W^{1,\infty}_{cl}S^{1,s}
({\mathcal O}).$ Then, according to Lemma \ref{Fops3}, Lemma \ref{Fops0} and Lemma \ref{Fops1}
$$
A(\widetilde y,\widetilde D,s)^*A(\widetilde y,\widetilde D,s)
=\widetilde \gamma\mbox{Re}\, p(\widetilde y,\widetilde D,s)-\widetilde\gamma\frac{\hat
C}{2}M^2(\widetilde D,s)+R,
$$
where $R\in {\mathcal L}(H^{\frac 12,1,s}_0({\mathcal O});L^2({\mathcal O}))$ and
\begin{eqnarray}\label{sobol}
\Vert R\Vert_{{\mathcal L}(H^{\frac 12,1,s}({\mathcal O});L^2({\mathcal O}))}
\le C_{7}(\pi_{W^{1,\infty}(\mathcal O_3)}(a)\pi_{W^{0,\infty}(\mathcal O_3)}(a))\nonumber\\
\le C_{8}(\pi_{W^{0,\infty}(\mathcal O_3)}(p)+1)\pi_{C^{0}(\mathcal O_3)}+
\sum_{k=0}^1(\pi_{W^{k,\infty}(\mathcal O_3)}(p)+1)\pi_{C^{1-k}(\mathcal O_3)}(\widetilde\gamma)).
\end{eqnarray}

Therefore
\begin{eqnarray}
\mbox{Re}(P(\widetilde y,\widetilde D,s)u,u)_{L^2({\mathcal O})}
=\Vert A(\widetilde y,\widetilde D,s)\Vert^2_{L^2({\mathcal O})}
-((1-\widetilde\gamma)
M^2(\widetilde D,s)u,u)_{L^2({\Bbb R}^n)}\nonumber\\
+\frac{\hat C}{2}\Vert u\Vert^2_{H^{\frac 12,1,s}({\mathcal
O})}+(Ru,u)_{L^2({\mathcal O})}.\nonumber
\end{eqnarray}
Observing that by (\ref{sobol})\begin{eqnarray}\vert(Ru,u)_{L^2({\mathcal O})}\vert
\\
 \le C_{9}(\pi_{W^{0,\infty}(\mathcal O_3)}(p)+1)\pi_{C^{0}(\mathcal O_3)}
(\widetilde \gamma)+\sum_{k=0}^1(\pi_{W^{k,\infty}(\mathcal O_3)}(p)+1)\pi_{C^{1-k}(\mathcal O_3)}(\widetilde\gamma))\Vert u\Vert_{L^2({\mathcal O})}
\Vert u\Vert_{H^{\frac 12,1,s}({\mathcal O})},\nonumber
\end{eqnarray}
and since by Lemma \ref{Fops4}
 we have
$$
\vert((1-\widetilde\gamma)M^2(\widetilde D,s)u,u)
_{L^2({\Bbb R}^n)}\vert \le C_{10}\pi_{C^{1}(\mathcal O_3)}(\widetilde\gamma)\Vert u\Vert_{L^2({\mathcal O})}
\Vert u\Vert_{H^{\frac 12,1,s}({\mathcal O})},
$$
we obtain the statement of the lemma.
$\blacksquare$

\begin{lemma}\label{Fops5} Let $p\in \{\frac 12, \frac 32\}$. Then there exist a constant $C_{11}$ independent of $s$
\begin{equation}\label{norm}
\sum_{\ell=-\infty}^\infty\Vert \text{\bf z}_\ell\Vert^2_{H^{\frac p2,p,\widetilde s}(\Sigma)}\le C \Vert \text{\bf z}\Vert^2_{H^{\frac p2,p,\widetilde s}(\Sigma)}.
\end{equation}
\end{lemma}

{\bf Proof.}
Observe that
\begin{equation}\label{grom}
\Vert \text{\bf z}\Vert^2_{H^{\frac p2,p,\widetilde s}(\Sigma)}=(\Vert \text{\bf z}\Vert^2_{H^{0,p,\widetilde s}(\Sigma)}+\Vert \text{\bf z}\Vert^2_{H^{\frac p2,0}(\Sigma)})^\frac 12.
\end{equation}
Since $\kappa_\ell$ depends only on $x_0$  by (\ref{gorokn}) we have
\begin{equation}\label{grom1}
\sum_{\ell=-\infty}^\infty\Vert \text{\bf z}_\ell\Vert^2_{H^{0,p,\widetilde s}(\Sigma)}\le C_{12} \Vert \text{\bf z}\Vert^2_{H^{0,p,\widetilde s}(\Sigma)}.
\end{equation}
Now we estimate the first term in the right hand side of (\ref{grom}). Let $\widetilde {\mbox{\bf z}} =\mbox{\bf z}\circ F^{-1}$ and the mapping $F$ given by (\ref{kupol}). By definition of norm in Sobolev-Slobodetskii space we have
$$
\sum_{\ell=-\infty}^\infty\Vert \kappa_\ell \mbox{\bf z}\Vert^2_{H^{\frac p2,0}(\Sigma)}=
%\sum_{\ell=-\infty}^\infty\Vert \kappa_\ell  \mbox{\bf z}\Vert^2_{H^{\frac p2,0}([-T,T]\times \partial\Omega)}=
\sum_{\ell=-\infty}^\infty\Vert \kappa_\ell \widetilde{\mbox{\bf z}}\Vert_{H^{\frac p2,0}([-T,T]\times \Bbb R^{n-1})}^2=
$$
$$
\sum_{\ell=-\infty}^\infty\int_{\Bbb R^{n-1}}\int_{[-T,T]\times[-T,T]} \frac{\vert(\kappa_\ell \widetilde{\mbox{\bf z}})(y_0)-(\kappa_\ell \widetilde{\mbox{\bf z}})(x_0)\vert^2}{\vert y_0-x_0\vert^{p+1}} dx_0dy_0dx'
$$
$$
\le 4\sum_{\ell=-\infty}^\infty\int_{\Bbb R^{n-1}}\int_{[-T,T]\times[-T,T]} \frac{\vert \kappa_\ell(y_0)-\kappa_\ell (x_0)\vert^2\vert \widetilde{\mbox{\bf z}}(y_0,x')\vert^2}{\vert y_0-x_0\vert^{p+1}} dx_0dy_0dx'
$$
$$
+4\sum_{\ell=-\infty}^\infty\int_{\Bbb R^{n-1}}\int_{[-T,T]\times[-T,T]} \frac{\vert \widetilde{\mbox{\bf z}}(y_0,x')-\widetilde{\mbox{\bf z}}(x_0,x')\vert^2\vert \kappa_\ell(y_0)\vert^2}{\vert y_0-x_0\vert^{p+1}} dx_0dy_0dx'=I_1+I_2.
$$

We estimate terms $I_j$ separately.
By (\ref{gorokn}) we have
$$
I_2=\sum_{\ell=-\infty}^\infty\int_{[-T,T]}\vert \kappa_\ell(y_0)\vert^2\int_{\Bbb R^{n-1}}\int_{[-T,T]} \frac{\vert \widetilde{\mbox{\bf z}}(y_0,x')-\widetilde{\mbox{\bf z}}(x_0,x')\vert^2 }{\vert y_0-x_0\vert^{p+1}} dxdy_0
$$
%$$
%\le \sum_{\ell=-\infty}^\infty\int_{[-T,T]} \kappa_\ell(y_0)\int_{\Bbb R^{n-1}}\int_{[-T,T]} \frac{\vert \widetilde{\mbox{\bf %z}}(y_0,x')-\widetilde{\mbox{\bf z}}(x_0,x')\vert^2 }{\vert y_0-x_0\vert^{p+1}} dxdy_0
%$$
$$
\le C_{13}\int_{[-T,T]} \int_{\Bbb R^{n-1}}\int_{[-T,T]} \frac{\vert \widetilde{\mbox{\bf z}}(y_0,x')-\widetilde{\mbox{\bf z}}(x_0,x')\vert^2 }{\vert y_0-x_0\vert^{p+1}} dxdy_0= C_{14} \Vert \widetilde{\mbox{\bf z}}\Vert^2_{H^{\frac p2,0}([-T,T]\times\Bbb R^{n-1})}.
$$
We set $\mathcal Z_1=\{(x_0,y_0)\in [-T,T]\times [-T,T] \vert \frac{1}{\vert x_0-y_0\vert}\le \lambda \widetilde \varphi(y_0)\}$ where $\lambda$ is the large positive parameter.
The short computations imply
$$I_1=\sum_{\ell=-\infty}^\infty\int_{\Bbb R^{n-1}}\int_{[-T,T]\times[-T,T]} \frac{\vert \kappa_\ell(y_0)-\kappa_\ell (x_0)\vert^2\vert \widetilde{\mbox{\bf z}}(y_0,x')\vert^2}{\vert y_0-x_0\vert^{p+1}} dx_0dy_0dx'=
$$
$$=\sum_{\ell=-\infty}^\infty\int_{\Bbb R^{n-1}}\int_{\mathcal Z_1} \frac{\vert \kappa_\ell(y_0)-\kappa_\ell (x_0)\vert^2\vert \widetilde{\mbox{\bf z}}(y_0,x')\vert^2}{\vert y_0-x_0\vert^{p+1}} dx_0dy_0dx'+
$$
$$\sum_{\ell=-\infty}^\infty\int_{\Bbb R^{n-1}}\int_{[-T,T]\times[-T,T]\setminus \mathcal Z_1} \frac{\vert \kappa_\ell(y_0)-\kappa_\ell (x_0)\vert^2\vert \widetilde{\mbox{\bf z}}(y_0,x')\vert^2}{\vert y_0-x_0\vert^{p+1}} dx_0dy_0dx'= P_1+P_2.
$$
Using the definition of the set $\mathcal Z_1$ we have
$$
P_2\le C_{15}(\lambda)\sum_{\ell=-\infty}^\infty\int_{\Bbb R^{n-1}}\int_{\mathcal Z_1} \frac{\vert \kappa_\ell(y_0)-\kappa_\ell (x_0)\vert^2\vert \widetilde{\mbox{\bf z}}(y_0,x')\vert^2}{\vert y_0-x_0\vert^{1-\delta}\widetilde \varphi^\frac{p+\delta}{3}(y_0)} dx_0dy_0dx'
$$
$$\le 2C_{15}(\lambda) \sum_{\ell=-\infty}^\infty\int_{\Bbb R^{n-1}}\int_{[-T,T]\times[-T,T]} \frac{(\vert \kappa_\ell(y_0)\vert^2+\vert\kappa_\ell (x_0)\vert^2)\vert \widetilde{\mbox{\bf z}}(y_0,x')\vert^2}{\vert y_0-x_0\vert^{1-\delta}\widetilde \varphi^\frac{p+\delta}{3}(y_0)} dx_0dy_0dx'
$$
$$
\le C_{16}(\lambda)\int_{\Bbb R^{n-1}}\int_{[-T,T]}\frac{\vert \widetilde{\mbox{\bf z}}(y_0,x')\vert^2}{\widetilde \varphi^\frac{p+\delta}{3}(y_0)} sup_{y_0\in [-T,T]}\Vert \vert y- x\vert\Vert^2_{L^2[-T,T]} dy_0dx'
$$
$$+C_{17}(\lambda)\sum_{\ell=-\infty}^\infty\int_{\Bbb R^{n-1}}\int_{[-T,T]}\frac{\vert \kappa_\ell(y_0)\vert^2\vert \widetilde{\mbox{\bf z}}(y_0,x')\vert^2}{\widetilde \varphi^\frac{p+\delta}{3}(y_0)} sup_{y_0\in [-T,T]}\Vert \vert y-x\vert\Vert^2_{L^2[-T,T]} dy_0dx'
$$
$$
\le C_{18} \Vert \widetilde{\mbox{\bf z}}/\widetilde \varphi^\frac{p+\delta}{6}\Vert_{L^2([-T,T]\times \Bbb R^{n-1})}.
$$

Next we estimate $P_1$. For any positive $\epsilon_1$ there exists a constant $C_{19}(\epsilon)$ such that
$$
P_1\le C_{19}\bigg(\sum_{\ell=-\infty}^\infty\int_{([T,T-\epsilon_1]\times [T,T-\epsilon_1]\setminus \mathcal Z_1)\times \Bbb R^{n-1}} \frac{\vert \kappa_\ell(y_0)-\kappa_\ell (x_0)\vert^2\vert \widetilde{\mbox{\bf z}}(y_0,x')\vert^2}{\vert y_0-x_0\vert^{1+p}} dx_0dy_0dx'
$$
$$
+\sum_{\ell=-\infty}^\infty\int_{([-T,-T+\epsilon_1]\times [T,-T+\epsilon_1]\setminus \mathcal Z_1)\times \Bbb R^{n-1}} \frac{\vert \kappa_\ell(y_0)-\kappa_\ell (x_0)\vert^2\vert \widetilde{\mbox{\bf z}}(y_0,x')\vert^2}{\vert y_0-x_0\vert^{1+p}} dx_0dy_0dx'
$$
$$
+\Vert\widetilde{\mbox{\bf z}}\Vert^2_{L^2([-T,T]\times \Bbb R^{n-1})}\bigg ).
$$
We estimate the first term in this inequality. The estimate of the second one is the same

If $(x_0,y_0)\in [T,T-\epsilon_1]\times [T,T-\epsilon_1]\setminus \mathcal Z_1$ we have
$$
\vert x_0-y_0\vert\le 2\vert T-y_0\vert^\frac 3p.
$$
So
$$
\vert \kappa(x_0)-\kappa(y_0)\vert \le sup_{\lambda\in [0,1]} \frac{\vert y_0-x_0\vert}{\vert T-\lambda x_0-(1-\lambda) y_0\vert^\frac 3p}\le  \frac{C_{20}\vert y_0-x_0\vert}{(\vert T-y_0\vert -\vert x_0- y_0\vert)^\frac 94}\le \frac{C_{21}\vert y_0-x_0\vert}{\vert T-y_0\vert^\frac 94}.
$$
Using this inequality we obtain
$$
 \sum_{\ell=-\infty}^\infty\int_{([T,T-\epsilon_1]\times [T,T-\epsilon_1]\setminus \mathcal Z_1)\times \Bbb R^{n-1}} \frac{\vert \kappa_\ell(y_0)-\kappa_\ell (x_0)\vert^2\vert \widetilde{\mbox{\bf z}}(y_0,x')\vert^2}{\vert y_0-x_0\vert^{1+p}} dx_0dy_0dx'
$$
$$
 \le C_{22}\sum_{\ell=-\infty}^\infty\int_{([T,T-\epsilon_1]\times [T,T-\epsilon_1]\setminus \mathcal Z_1)\times \Bbb R^{n-1}} \frac{\vert \kappa_\ell(y_0)-\kappa_\ell (x_0)\vert^{2-p-\delta}\vert \widetilde{\mbox{\bf z}}(y_0,x')\vert^2}{\vert y_0-x_0\vert^{1-\delta}\widetilde\varphi(y_0)^{\frac 34(p+\delta)}} dx_0dy_0dx'
$$
$$
 \le C_{23}\sum_{\ell=-\infty}^\infty\int_{([T,T-\epsilon_1]\times [T,T-\epsilon_1]\setminus \mathcal Z_1)\times \Bbb R^{n-1}} \frac{(\vert \kappa_\ell(y_0)\vert^{2-p-\delta}+\vert\kappa_\ell (x_0)\vert^{2-p-\delta} )\vert \widetilde{\mbox{\bf z}}(y_0,x')\vert^2}{\vert y_0-x_0\vert^{1-\delta}\widetilde\varphi(y_0)^{\frac 34(p+\delta)}} dx_0dy_0dx'
$$
$$
\le C_{24} \Vert\frac{\widetilde{\mbox{\bf z}}}{\widetilde\varphi^{\frac 38(p+\delta)}}\Vert^2_{L^2([-T,T]\times \Bbb R^{n-1})}.
$$
Taking parameter $\delta$ sufficiently small we obtain (\ref{norm}). Proof of the proposition is complete.
$\blacksquare$

\end{document}